\documentclass[twoside,11pt]{article}
\usepackage{blindtext}
\usepackage{jmlr2e}
\usepackage{amsmath}
\usepackage{amsfonts}
\usepackage{amssymb}
\usepackage{mathabx}
\usepackage{graphicx,enumitem}
\usepackage{color}

\def\T{{\mathrm{\scriptscriptstyle T}}}
\newcommand{\bX}{{\bf X}}
\newcommand{\be}{{\bf e}}
\newcommand{\bbb}{{\bf b}}
\newcommand{\bU}{{\bf U}}
\newcommand{\bZ}{{\bf Z}}
\newcommand{\bS}{{\bf S}}
\newcommand{\bR}{{\bf R}}
\newcommand{\bB}{{\bf B}}
\newcommand{\bY}{{\bf Y}}
\newcommand{\bI}{{\bf I}}
\newcommand{\bJ}{{\bf J}}
\newcommand{\bL}{{\bf L}}
\newcommand{\cU}{{\cal U}}
\newcommand{\bmu}{\boldsymbol \mu}
\newcommand{\btheta}{\boldsymbol \theta}
\newcommand{\bphi}{\boldsymbol \phi}

\newcommand{\bpsi}{\boldsymbol \psi}
\newcommand{\bSigma}{\boldsymbol \Sigma}
\newcommand{\bXi}{\boldsymbol \Xi}
\newcommand{\bLambda}{\boldsymbol \Lambda}

\newcommand{\cE}{\mathbb{E}}
\newcommand{\cov}{\text{cov}}
\newcommand{\var}{\text{var}}
\newcommand{\pr}{\mathbb{P}}
\newcommand{\cS}{{\cal S}}
\newcommand{\diag}{\mbox{diag}}
\newcommand{\sSigma}{\scriptscriptstyle{\Sigma}}
\newcommand{\tF}{\text{F}}
\def\p{{P}}

\newtheorem{assumption}{Assumption}

\usepackage{lastpage}
\jmlrheading{26}{2025}{1-\pageref{LastPage}}{11/23; Revised 5/24}{1/25}{23-1578}{Shaojun Guo, Dong Li, Xinghao Qiao, Yizhu Wang} \ShortHeadings{Phase Transitions for Functional Data in High Dimensions}{Guo, Li, Qiao, and Wang}

\firstpageno{1}

\begin{document}

\title{From Sparse to Dense Functional Data in High Dimensions: Revisiting Phase Transitions from a Non-Asymptotic Perspective}
%\title{From sparse to dense functional data in high dimensions: Revisiting phase transitions from a non-asymptotic perspective}

\author{\name Shaojun Guo \email sjguo@ruc.edu.cn \\
       \addr Institute of Statistics and Big Data\\
       Renmin University of China\\
       Beijing 100872, China
       \AND
       \name Dong Li \email malidong@tsinghua.edu.cn \\
       \addr Department of Statistics and Data Science\\
       Tsinghua University\\
       Beijing 100084, China
       \AND
       \name Xinghao Qiao\thanks{The authors' names are sorted alphabetically, and the corresponding author is Xinghao Qiao.}  \email xinghaoq@hku.hk \\
       \addr Faculty of Business and Economics,\\
       The University of Hong Kong\\
       Hong Kong, China
       \AND
       \name Yizhu Wang \email wangyizh20@mails.tsinghua.edu.cn \\
       \addr Department of Mathematical Sciences\\
       Tsinghua University\\
       Beijing 100084, China}

\editor{Genevera Allen}

\maketitle

\begin{abstract}%   <- trailing '%' for backward compatibility of .sty file
%\blindtext
Nonparametric estimation of the mean and covariance functions is ubiquitous in functional data analysis and local linear smoothing techniques are most frequently used. \cite{zhang2016} explored different types of asymptotic properties of the estimation, which reveal interesting phase transition phenomena based on the relative order of the average sampling frequency per subject $T$ to the number of subjects $n$, partitioning the data into three categories: ``sparse'', ``semi-dense'', and ``ultra-dense''. In an increasingly available high-dimensional scenario,
where the number of functional variables $p$ is large in relation to $n$, we revisit this open problem from a non-asymptotic perspective by deriving comprehensive concentration inequalities for the local linear smoothers. Besides being of interest by themselves, our non-asymptotic results lead to elementwise maximum rates of $L_2$ convergence and uniform convergence serving as a fundamentally important tool for further convergence analysis when $p$ grows exponentially with $n$ and possibly $T$. With the presence of extra $\log p$ terms to account for the high-dimensional effect, we then investigate the scaled phase transitions and the corresponding elementwise maximum rates
from sparse to semi-dense to ultra-dense functional data in high dimensions. 
We also discuss a couple of applications of our theoretical results.
Finally, numerical studies are carried out to confirm the established theoretical properties.
\end{abstract}


\section{Introduction}
\label{sec:intro}
A fundamental issue in functional data analysis is the nonparametric estimation of the mean and covariance functions based on discretely sampled and noisy curves. Despite being of interest by themselves, the estimated quantities serve as building blocks for dimension reduction and subsequent modeling of functional data,
%for frequently used approaches to handle functional data,
such as functional principal component analysis (FPCA) \cite[]{yao2005a,li2010} and functional linear regression \cite[]{yao2005b,chen2022}.
Among candidate nonparametric smoothers, we focus on the most commonly-adopted local linear smoothing method due to its simplicity and attractive local and boundary correction properties.

In a typical functional data setting, we have $n$ random curves, representing $n$ subjects, observed with errors, at $T_i$ randomly sampled time points for the $i$th subject. The sampling frequency $T_i$ plays a pivotal role in the estimation, as it may affect the choice of the estimation procedure. The literature can be loosely divided into two categories.
The first category corresponds to dense functional data, where $T_i$'s are larger than some order of $n.$ A conventional approach to handle such data implements nonparametric smoothing to the observations from each subject to eliminate the noise, thus reconstructing each individual curve before subsequent analysis \cite[]{zhang2007}.
The second category referred to as sparse functional data, accords with bounded $T_i$'s. Under such a scenario, the pre-smoothing step is no longer applicable, an alternative pooling strategy considers pooling the data from all subjects to build strength across all observations \cite[]{yao2005a,li2010}.
More recently, \cite{zhang2016} provided a comprehensive analysis of phase transitions and the associated rates of convergence for three types of asymptotic properties: local asymptotic normality, $L_2$ convergence, and uniform convergence. They proposed to further partition dense functional data into new categories: ``semi-dense'' and ``ultra-dense'', depending on whether the root-$n$ rate is achieved with negligible asymptotic bias or not. However, these aforementioned asymptotic results are only suitable for handling univariate or low-dimensional multivariate functional data.

With recent advances in data collection technology, high-dimensional functional data sets become increasingly available. Examples include time-course gene expression data, and electroencephalography and functional magnetic resonance imaging data, where signals are measured over time at a large number of regions of interest \cite[]{zhu2016,li2018,zapata2022,fang2022}.
Those data can be represented as a $p$-vector of random functions
$\bX_{i}(\cdot)=\{X_{i1}(\cdot), \dots, X_{ip}(\cdot)\}^\T$ for $i=1, \dots,n$ defined on a compact set $\cal U,$ with the $p$-vector of mean functions $\bmu(\cdot)=\{\mu_1(\cdot),\dots, \mu_p(\cdot)\}^{\T}=\cE\{\bX(\cdot)\}$ and the $(p \times p)$-matrix of marginal- and cross-covariance functions 
\begin{equation}
\label{cov.fmat}
\bSigma(u,v)=\{\Sigma_{jk}(u,v)\}_{p \times p}, ~~\Sigma_{jk}(u,v)=\cov\{X_{ij}(u), X_{ik}(v)\}.
\end{equation}
In a high-dimensional regime, the dimension $p$ can be diverging with, or even larger than, the number of subjects $n$. In practice, each $X_{ij}(\cdot)$ is observed subject to error contamination at $T_{ij}$ random time points. See (\ref{model}) below. %In this paper, we propose a unified local linear smoothing approach to estimate the mean and covariance functions for both sparsely and densely sampled high-dimensional functional data.

The estimation of mean and covariance functions for high-dimensional functional data is not only interesting at its own right but also plays a foundational role as a critical starting point for subsequent analysis. Due to the infinite-dimensionality of functional data, it is common practice to truncate the basis expansion of each function $X_{ij}(\cdot)$ at a finite level, using either data-driven basis expansion via FPCA or pre-fixed basis expansion. This allows for the subsequent development of regularized methods for the sparse estimation to address the high-dimensionality based on estimated FPC scores or estimated basis coefficients. Importantly, the estimation of mean and covariance functions are implicitly involved in this procedure.
Applications of this type of approach to estimate sparse high-dimensional functional models include 
functional graphical models \cite[]{qiao2019,zhao2022,solea2022,zapata2022,lee2021},
functional additive regressions \cite[]{fan2014,fan2015,kong2016,luo2017,xue2021}, sparse FPCA \cite[]{hu2022} and functional linear discriminant analysis \cite[]{xue2023}.
Another line of applications considers the sparsity-induced estimation of the covariance matrix function without any basis expansion, such as \cite{fang2022,li2023} and \cite{leng2024}. Both types of applications call for the investigation of non-asymptotic properties of the componentwise mean and covariance estimators, which serves as the motivation of our paper when dealing with the practical scenario of partially observed curves in high dimensions. See Section~\ref{sec:app} for details of some applications.

Within the high-dimensional statistical learning framework, it is essential to conduct non-asymptotic analysis of the estimators by developing concentration inequalities under a given performance metric, which can lead to probabilistic error bounds in the elementwise maximum norm as a function of $n$, $p$, and possibly $T_{ij}$'s (depending whether they are diverging or bounded) under our setup.
Existing literature has mainly focused on fully observed functional data, based on which concentration inequalities for the estimated covariance functions were established in \cite{qiao2019} and \cite{zapata2022}.
In practical scenarios where curves are partially observed with errors, addressing dense functional data is achievable by applying the pre-smoothing technique to observations from each 
$i, j$ \cite[]{kong2016}. Alternatively, a unified pooling-type local linear smoothing approach can be employed for estimating the mean functions $\mu_j(\cdot)$'s and marginal- and (or) cross-covariance functions $\Sigma_{jk}(\cdot,\cdot)$'s across $j,k$
to handle both sparsely and densely observed functional data \cite[]{qiao2020,lee2021,fang2022}. 
%Under certain lower-dimensional structural assumptions, one can possibly develop the nonparametric smoothing method for the joint estimation of elements in $\bSigma(\cdot,\cdot)$, which, however, becomes challenging due to the observation of each $X_{ij}(\cdot)$ at different sets of points. 
Although such nonparametric smoothing approach suffers from high computational cost when $p$ is large, it can be substantially accelerated in a common practical scenario where each $X_{ij}(\cdot)$ is observed at the same set of points across $j \in [p]$ especially with the aid of linear binning \cite[]{fan1994}, resulting in an efficient estimation procedure. 
%Such elementwise approach is computationally feasible as it can be easily parallelized and 
Moreover, those commonly-adopted FPCA-based methods only necessitate the estimation of marginal- instead of cross-covariance functions across $j \in [p]$ \cite[]{qiao2019,solea2022}, and can be easily paralleled for fast computation. 
See Remark~\ref{rmk_comp}.

On the theory side, this approach entails dealing with the second-order $U$-statistics with complex dependence structures, posing a technically challenging task. 
\cite{qiao2020} made the first attempt to derive some sub-optimal concentration inequalities for local linear smoothers of marginal-covariance functions $\widehat \Sigma_{jj}(\cdot,\cdot)$'s, albeit under a restrictive finite-dimensional setting. 
\cite{lee2021} established the convergence of their proposed estimation of conditional functional graphical models under the assumption of elementwise maximum rate for the covariance smoothers:
\begin{equation}
\label{max.cfgm}
\max_{1\leq j,k \leq p}\|\widehat \Sigma_{jk} - \Sigma_{jk}\|_{\cS}=O_\p(\log p n^{-\tau}),
\end{equation}
where $\|\cdot\|_{\cS}$ denotes the Hilbert--Schmidt norm, and the parameter $\tau\in(0,1/2]$ reflects the average sampling frequency, with larger values yielding denser observational points. 
\cite{fang2022} developed the functional covariance estimation with theoretical guarantees by assuming generalized sub-Gaussian-type concentration inequalities for local linear  smoothers $\widehat\Sigma_{jk}(\cdot,\cdot)$'s, resulting in an improved elementwise maximum rate: 
\begin{equation}
\label{max.fcm}
\max_{1\leq j,k \leq p}\|\widehat \Sigma_{jk} - \Sigma_{jk}\|_{\cS}=O_\p\{(\log p)^{1/2} n^{-\tau} +h^2\},
\end{equation}
where $h>0$ is the bandwidth parameter. However, it remains of theoretical interest to ask: 
\begin{itemize}
\item What are the exact forms of such rates as functions of $n,$ $p,$ $T_{ij}$'s, and associated bandwidth parameters under cases with different sampling frequencies? 
\item Are these rates well-established in the sense of specifying the largest values of $\tau$ and, compared to \cite{zhang2016}, exhibiting any corresponding phase transition phenomena in the high-dimensional setting?
\end{itemize}

This paper aims to fill crucial theoretical gaps related to local linear smoothers frequently adopted in existing literature. Specifically, we
present a systematic and unified non-asymptotic analysis of local linear smoothers for the mean and covariance functions to accommodate both sparsely and densely observed functional data in high dimensions. While our focus is not to introduce new methodologies for handling high-dimensional partially observed functional data, we make three new contributions as follows.
%The new contribution of this paper is threefold.

%While our paper does not primarily aim to introduce new methodology for handling high-dimensional partially observed functional data, it focuses on addressing crucial theoretical gaps related to local linear smoothers frequently adopted in existing literature. Specifically, we
%present a systematic and unified non-asymptotic analysis of the local linear smoothers for the mean and covariance functions to accommodate both sparsely and densely observed functional data in high dimensions.
%it remains an important open problem to investigate the non-asymptotic proprieties of both the mean and covariance estimators when handling partially observed functional data.
%precise specifications of the largest values of $\tau$ under cases with different sampling frequencies is 
%developed their methodologies by simply assuming that the elementwise maximum rates of the covariance estimators are not faster than the parametric counterpart without specifying the exact forms. 
%This paper presents a systematic and unified non-asymptotic analysis of the local linear smoothers for the mean and covariance functions to accommodate both sparsely and densely observed functional data in high dimensions. 
%The new contribution of this paper is threefold.
\begin{itemize}
\item First, we develop generalized sub-Gaussian-type concentration inequalities for each functional element of the mean and covariance estimators in both $L_2$ norm and supremum norm. Compared to the asymptotic results in \cite{zhang2016}, our non-asymptotic error bounds lead to the same rates of $L_2$ convergence and uniform convergence, and reveal the same phase transition phenomena depending on the relative order of the average sampling frequency per subject to $n^{1/4}$ for dense functional data. See Remarks~\ref{rmk_asy_rates} and \ref{rmk_asy_phase}.

\item Second, we derive elementwise maximum rates of both $L_2$ and uniform convergence for the mean and covariance estimators. Notably, we fundamentally improve the rates (\ref{max.cfgm}) and (\ref{max.fcm}) assumed in existing literature in the sense of precisely specifying the largest values of $\tau$ under cases with different sampling frequencies. These established rates in Theorems~\ref{thm_mean_maxrate} and \ref{thm_cov_maxrate} serve as a foundational tool to provide theoretical guarantees for a set of 
aforementioned sparse high-dimensional functional models in the existing literature when dealing with the practical scenario of partially observed functional data in high dimensions.
%that can handle high-dimensional partially observed functional data.
%, such as functional graphical models \cite[]{li2018,qiao2019,zhao2022,solea2022,zapata2022,lee2021,zhao2023},functional additive regressions \cite[]{fan2014,fan2015,kong2016,luo2017,wang2022}, sparse FPCA \cite[]{hu2022} and functional covariance estimation \cite[]{fang2022,li2023}.

\item Third, with the presence of additional $\log p$ terms to account for the high-dimensional effect in our established elementwise maximum rates, the scaled phase transitions for high-dimensional dense functional data occur based on the ratios of the average sampling frequency per subject to $n^{1/4}(\log p)^{-1/4}.$ This leads to a further partition of dense functional data into categories of  ``semi-dense'' and ``ultra-dense'', depending on whether the parametric rate $(\log p)^{1/2}n^{-1/2}$ can be attained or not. With suitable choices of optimal bandwidths, we also present the optimal elementwise maximum rates from sparse to semi-dense to ultra-dense functional data, which correspondingly extend the optimal rates in \cite{zhang2016} to the high-dimensional setting. See Remarks~\ref{rmk_nonasy_phase} and \ref{rmk_nonasy_phase2}.
\end{itemize}

{\bf Outline of the paper}.
Section~\ref{sec.method} presents the nonparametric smoothing approach to estimate the mean and covariance functions.
In Section~\ref{sec.theory}, we investigate the non-asymptotic properties of the proposed local linear smoothers and discuss the associated phase transition phenomena.
In Section~\ref{sec:app}, we outline a couple of applications of the non-asymptotic theory for the local linear smoothers.
The established theoretical results are validated through simulations in Section~\ref{sec.sim}.
All technical proofs are relegated to the appendix.

{\bf Notation}.
We summarize here some notation to be used throughout the paper. For a positive integer $q$, we write $[q]=\{1, \dots, q\}$.
 For $x, y \in \mathbb{R},$ we write $x \vee y=\max(x,y)$ and $x \wedge y = \min(x,y)$. We use $I(\cdot)$ to denote an indicator function. %\textcolor{red}{and $\|\cdot\|$  to denote Euclidean norm for matrices and vectors.}
We use $\bI_p$ to denote the $p \times p$ identity matrix.
Let  $L_2(\cU)$ be a Hilbert space of square-integrable functions on a compact interval $\cU$ equipped with the inner product $\langle f, g\rangle =\int f(u) g(u) {\rm d}u$ for $f(\cdot)$, $g(\cdot) \in L_2(\cU)$ and the induced $L_2$ norm $\|\cdot\|_2=\langle \cdot, \cdot \rangle^{1/2}.$ For any bivariate function $\Phi(\cdot,\cdot)$ in $L_2(\cU \times \cU),$ we also use $\Phi$ to denote the linear  operator induced from the kernel function
$\Phi(\cdot,\cdot),$ that is, for any $f(\cdot) \in L_2(\cU),$ $\Phi(f)(\cdot)=\int \Phi(\cdot,v)f(v) {\rm d}v \in L_2(\cU),$ and denote
its Hilbert--Schmidt norm by $\{\int\int \Phi(u,v)^2 {\rm d}u {\rm d}v\}^{1/2}.$ 
For two positive sequences $\{a_n\}$ and $\{b_n\},$ we write $a_n \lesssim b_n$ or $b_n \gtrsim a_n$ if there exist a positive constant $c$ such that $\lim \sup_{n \rightarrow \infty} a_n/b_n \leq c.$ We write $a_n \asymp b_n$ if and only if $a_n \lesssim b_n$ and $b_n \lesssim a_n$ hold simultaneously.

\section{Methodology}\label{sec.method}
%\subsection{Estimation}\label{sec_est}
Let $\bX_{i}(\cdot)=\{X_{i1}(\cdot), \dots, X_{ip}(\cdot)\}^\T$ for $i \in [n]$ be  independently and identically distributed
copies of $\bX(\cdot)$ defined on $\cU$ 
with mean $\bmu(\cdot)$ and covariance $\bSigma(\cdot,\cdot).$
%with continuous $p$-vector of mean functions $\bmu(\cdot)=\{\mu_1(\cdot),\dots, \mu_p(\cdot)\}^{\T}=\cE\{\bX(\cdot)\}$ and $(p \times p)$-matrix of marginal- and cross-covariance functions $\bSigma(u,v)=\{\Sigma_{jk}(u,v)\}_{p \times p}$ whose $(j,k)$-th entry $\Sigma_{jk}(u,v)=\cov\{X_{ij}(u), X_{ik}(v)\}.$ 
For any $i \in [n]$ and $j \in [p],$ $X_{ij}(\cdot)$ is not directly observable in practice. Instead it is observed, with random errors, at $T_{ij}$ random time points, $U_{ij1}, \dots, U_{ijT_{ij}} \in {\cal U}.$ % where for dense measurement designs, all $T_{ij}$'s are larger than some order of $n$ and for sparse designs, all $T_{ij}$'s are bounded \cite[]{zhang2016,qiao2019}.
Let $Y_{ijt}$ be the observed value of $X_{ij}(U_{ijt})$ satisfying
\begin{equation}
\label{model}
Y_{ijt} = X_{ij}(U_{ijt}) + \varepsilon_{ijt},
\end{equation}
where the errors $\varepsilon_{ijt}$'s, independent of $X_{ij}$'s, are independently and identically distributed copies of $\varepsilon_j$ with $\cE(\varepsilon_j)=0$ and $\var(\varepsilon_j)=\sigma_j^2 <\infty.$

Based on the observed data $\{(U_{ijt},Y_{ijt}): i\in [n], j\in [p], t\in [T_{ij}]\},$ we present a unified procedure to estimate the mean functions $\mu_j(\cdot)$'s and the marginal- and cross-covariance functions $\Sigma_{jk}(\cdot,\cdot)$'s for both sparsely and densely observed functional data. In what follows, denote $K_h(\cdot)=h^{-1}K(\cdot/h)$ for a univariate kernel $K$ with bandwidth $h>0$. For each $j$, a local linear smoother is firstly applied to $\{(U_{ijt}, Y_{ijt}): i\in [n],
t\in [T_{ij}]\}$, and hence the estimated mean function is attained via $\hat\mu_j(u)=\hat b_0$, where
\begin{equation*}\label{mean_crit}
   (\hat b_0, \hat b_1)=\arg\min_{b_0,b_1}\sum_{i=1}^n v_{ij} \sum_{t=1}^{T_{ij}}\Big\{Y_{ijt} - b_0 - b_1(U_{ijt}-u)\Big\}^2 K_{h_{\mu,j}}(U_{ijt}-u).
\end{equation*}
The weight $v_{ij}$ is attached to each observation for the $i$th subject and the $j$th functional variable such that $\sum_{i=1}^n T_{ij}v_{ij}=1$ \cite[]{zhang2016}.

For each $i \in [n], j, k\in [p], t\in [T_{ij}]$ and $s\in [T_{ik}],$ once the mean functions are estimated, let $\Theta_{ijkts}=\{Y_{ijt}-\hat\mu_j(U_{ijt})\}\{Y_{iks}-\hat\mu_k(U_{iks})\}$ be the ``raw covariance" between $Y_{ijt}$ and $Y_{iks}.$
Notice that
$\cov(Y_{ijt},Y_{iks})=\Sigma_{jk}(U_{ijt},U_{iks})+\sigma_j^2 I(j=k)I(t=s)$. To estimate the marginal-covariance function $\Sigma_{jj}(\cdot,\cdot)$ for each $j$ or the cross-covariance function $\Sigma_{jk}(\cdot,\cdot)$ for each $j \neq k$, we employ local linear surface smoothers to the off-diagonals of the raw marginal-covariances $(\Theta_{ijjts})_{1 \leq t \neq s \leq T_{ij}}$ or to the raw cross-covariances $(\Theta_{ijkts})_{t \in [T_{ij}], s \in [T_{ik}]}$. Specifically, we
minimize
\begin{equation*}\label{mcov_crit}
  \sum_{i=1}^n w_{ijk} \underset{(t,s) \in {\cal T}}{\sum}
\Big\{\Theta_{ijkts} - \beta_0 - \beta_1(U_{ijt}-u)- \beta_2(U_{iks}-v)\Big\}^2K_{h_{\sSigma,jk}}(U_{ijt}-u) K_{h_{\sSigma,jk}}(U_{iks}-v)
\end{equation*}
with respect to $(\beta_0, \beta_1, \beta_2)$,
where the set ${\cal T}$ equals to $\{(t,s): t \in [T_{ij}], s \in [T_{ij}], t \neq s\}$ if $j=k$ or $\{(t,s): t \in [T_{ij}], s \in [T_{ik}]\}$ if $j\neq k,$ and the weight $w_{ijk}$ is assigned to each triplet $(i,j,k)$ such that $\sum_{i=1}^n T_{ij}\{T_{ik}-I(j=k)\}w_{ijk}=1$. See the weights to estimate marginal-covariance functions in \cite{zhang2016}. The resulting marginal- or cross-covariance estimator is $\widehat \Sigma_{jk}(u,v)=\hat \beta_0$. %We refer to \cite{yao2005a,li2010} and \cite{zhang2016} for details on different choices of the weights $v_{ij}$ and $w_{ij}$ \cite[]{zhang2016}.
For ease of presentation, we assume that the mean functions $\mu_j(\cdot)$'s are known in advance when discussing the concentration and convergence results related to the covariance estimators $\widehat \Sigma_{jk}(u,v)$'s in Section~\ref{sec.theory} below. However, it is noteworthy that all our discussions remain valid even when $\mu_j(\cdot)$'s are unknown as long as a few additional technical assumptions are imposed.

%To estimate the cross-covariance function $C_{jk}(\cdot,\cdot)$ for each $j \neq k,$ we finally implement a local linear surface smoother to the raw cross-covariances $(\Sigma_{ijkts})_{t \in [T_{ij}], s \in [T_{ik}]}$ by minimizing
%\begin{equation}
%\label{ccov_crit}
%  \sum_{i=1}^n \tilde w_{ijk} \sum_{t=1}^{T_{ij}} \sum_{s=1}^{T_{ik}}\Big\{\Sigma_{ijkts} - \beta_0 - \beta_1(U_{ijt}-u)- \beta_2(U_{iks}-v)\Big\}^2 K_{h_{Cjk}}(U_{ijt}-u) K_{h_{Cjk}}(U_{iks}-v).
%\end{equation}
%with respect to $(\beta_0,\beta_1, \beta_2),$ thus obtaining the cross-covariance estimator $\widehat C_{jk}(u,v)=\hat \beta_0.$
%The weight $\tilde w_{ijk}$ is attached to each $\Sigma_{ijkts}$ for the triplet $(i,j,k)$ such that $\sum_{i=1}^n T_{ij} T_{ik} \tilde w_{ijk}=1.$

Our estimation procedure allows general weighting schemes for $\{v_{i}\}_{i \in [n]},$ $\{w_{ijk}\}_{i \in [n], j,k \in [p]}$ such that two types of frequently-used schemes in existing literature are special cases of them. One type assigns the same weights to each observation \cite[]{yao2005a} with $v_{ij}=(\sum_{i=1}^n T_{ij})^{-1}$ and $w_{ijk}=[\sum_{i=1}^n T_{ij}\{T_{ik}-I(j=k)\}]^{-1},$ so a subject with a larger number of observations receives more weights in total.
The other type assigns the same weights to each subject \cite[]{li2010}, thus leading to $v_{ij}=(nT_{ij})^{-1}$ and $w_{ijk}=[nT_{ij}\{T_{ik}-I(j=k)\}]^{-1}.$

\begin{remark}
\label{rmk_comp}
(i) Suppose that the estimated mean and covariance functions are evaluated at a grid of $R \times R$ locations over $\cU^2.$ Under high-dimensional settings, it is apparent that our nonparametric smoothing approach suffers from high computational cost in kernel evaluations, particularly when estimating $p(p+1)/2$ marginal- and cross-covariance functions. In a common practical scenario, where each $X_{ij}(\cdot)$ is observed at the same set of time points $U_{1}, \dots, U_{T_i} \in \cU$ across $j \in [p],$ model~(\ref{model}) simplifies to
\begin{equation*}
\label{eq.kernel.simple}
Y_{ijt} = X_{ij}(U_{it}) + \varepsilon_{ijt}, ~~~ t=1, \dots, T_i,
\end{equation*}
which reduces the number of kernel evaluations from $O(\sum_{i=1}^n \sum_{j=1}^p T_{ij}R)$ to $O(\sum_{i=1}^n T_i R),$ substantially accelerating the computation in high-dimensional settings. 
To further speed up the computation, we can use the linear binning technique \cite[]{fan1994} to approximate the mean and covariance estimation. This would largely reduce the number of kernel evaluations to $O(R),$ and the number of operations (i.e. additions and multiplications) to $O(np^2R^2 + p^2R^4 + p \sum_{i=1}^n T_i),$ from $O(p^2R^2 \sum_{i=1}^n T_i^2).$
See the detailed implementation of binning and the associated computational complexity analysis in \cite{fang2022}. Our conducted numerical experiments show that such binned implementation offers significantly improved computational efficiency without sacrificing any estimation accuracy.

(ii) In a general scenario when $X_{ij}(\cdot)$'s are observed at different sets of time points, parallel estimation for $j, k \in [p]$ can be employed, resulting in a more efficient procedure. In contrast, under certain lower-dimensional structural assumptions, the possible development of a nonparametric smoothing method for the joint estimation of components in $\bSigma(\cdot,\cdot)$ becomes challenging in this general scenario, and is thus left for future research. 
 
(iii) Due to the infinite-dimensional nature of functional data, it is standard practice to employ FPCA as a first dimension reduction step before subsequent modeling in the commonly-adopted multi-step estimation when dealing with high-dimensional functional models. See the detailed discussion of this application in Section~\ref{sec_fpca}.
Specifically, implementing FPCA necessitates the nonparametric estimation of only marginal-covariance functions $\Sigma_{jj}(\cdot,\cdot)$ across $j \in [p],$ which can be easily paralleled for fast computation.
\end{remark}

%{\color{red} The weight $\tilde w_{ijk}$ is assigned to each $(i,j,k)$ such that $\sum_{i=1}^n T_{ij} T_{ik} \tilde w_{ijk}=1.$ In the weighing scheme, there are two frequently used schemes. One is the scheme which assigns the same weight to each observation, we have $ v_{ij} = 1 / \sum_{i=1}^{n} T_{ij} $, $ \omega_{ij} = 1 / \sum_{i=1}^{n} T_{ij}(T_{ij}-1) $ and $ \widetilde{\omega}_{ijk} = 1 / \sum_{i=1}^{n} T_{ij}T_{ik} $. And the other scheme which assigns the same weight to each subject, we have $ v_{ij} = 1 / n T_{ij} $, $ \omega_{ij} = 1 / n T_{ij}(T_{ij}-1) $ and $ \widetilde{\omega}_{ijk} = 1 / n T_{ij}T_{ik} $.}

%Note that we drop subscript $j$ of $h_{M,j}$ in (\ref{mcov_crit}) and subscripts $j,k$ of $h_{C,jk}$ in (\ref{ccov_crit})to simplify the notation in this section. However, we select different bandwidths $h_{C,jk}$ and $h_{M,j}$ across $j,k \in [p]$ in our empirical studies.

%\textcolor{red}{Consider a common practical situation, where binning \cite{fan1994,fang2022}}

\section{Theory}\label{sec.theory}
Before presenting the concentration and convergence results, we impose the following regularity assumptions.

\begin{assumption}\label{cond_subG}
For each $i\in [n]$ and $j \in [p],$ $X_{ij}(\cdot)$ is a sub-Gaussian random process and $\varepsilon_{ij}$ is a sub-Gaussian random variable, that is, there exists some positive constant $c$ such that $\cE\{\exp (\langle x, X_{ij} - \mu_j\rangle)\} \leq \exp \{2^{-1} c^2 \langle x, \Sigma_{jj}(x)\rangle \}$ for all $x(\cdot) \in L^2(\cU)$ and $\cE\{\exp(\varepsilon_{ij} z)\} \leq \exp(c^2 \sigma_j^2 z^2/2)$ for all $z \in {\mathbb R}.$
\end{assumption}

\begin{assumption}\label{cond_T}
For each $i \in [n]$ and $j \in [p],$ under the sparse design, $T_{ij} \leq T_0 <\infty,$ and, under the dense design, $T_{ij} \rightarrow \infty$ and there exists some positive constant $\bar c$ such that $\max_{i,j} T_{ij} (\min_{i,j} T_{ij})^{-1} \leq \bar c.$
%with $U_{ijt}$'s independent of $i.$
\end{assumption}

\begin{assumption}\label{cond_weight}
		Under the dense design, there exists some positive constant $c_0$ such that
  $\max_{i,j} v_{ij} (\min_{i,j} v_{ij})^{-1} \leq c_0$ and
%  $\max_{i,j} w_{ij} (\min_{i,j} w_{ij})^{-1} \leq c_0,$
  $\max_{i,j,k} w_{ijk} (\min_{i,j,k} w_{ijk})^{-1} \leq c_0.$
%		\begin{equation}
%			\begin{split}
%				& \frac{\max_{i.j} v_{ij}}{\min_{i,j} v_{i,j}} \leq K \\
%				& \frac{\max_{i.j} w_{ij}}{\min_{i,j} w_{i,j}} \leq K \\
%				& \frac{\max_{i.j,k} \tilde w_{ijk}}{\min_{i.j,k} \tilde w_{ijk}} \leq K \\
%				\nonumber
%			\end{split}
%		\end{equation}
	\end{assumption}

%\begin{condition}
%		\label{cond_timepoints}
%		Under the dense design, there exists some constant $\alpha \in (0, 1)$ such that
%		$\min_{1 \leq t \neq s \leq T_{ij}} |U_{ijt} - U_{ijs} | \geq \alpha %(T_{ij})^{-1}$ for any $i \in [n]$ and $j \in [p].$
%\end{condition}

\begin{assumption}
	\label{cond_T_bd}
	 (i) Let $\big\{U_{ijt}: i \in [n],j \in [p], t \in [T_{ij}]\big\}$ be independently and identically distributed copies of a random variable $U$ defined on $\cU.$ The density $f_{U}(\cdot)$ of $U$ satisfies $0<m_f \leq {\inf}_{\cU}f_{U}(u)\leq {\sup}_{\cU}f_{U}(u) \leq M_f<\infty$ for some positive constants $m_f$ and $M_f;$ (ii) $\bX(\cdot),$ $U$ and $\{\varepsilon_j\}_{j \in [p]}$ are mutually independent.
\end{assumption}

\begin{assumption}
		\label{cond_timepoints}
  Let $B_{jk} = [(k-1) \widetilde T_{j}^{-1},k \widetilde T_{j}^{-1}]$ for $k \in [\widetilde T_j]$ with $\widetilde T_{j} = \max_{i}T_{ij},$ there exists some constant $C>0$ such that the cardinality
$
\#\big\{U_{ijt}: U_{ijt} \in B_{jk}, t \in [T_{ij}]\big\} \le C
$
for each $i \in [n],j \in [p]$ and $k \in [\widetilde T_j].$
		%Under the dense design, there exists some constant $\alpha \in (0, 1)$ such that
		%$\min_{1 \leq t \neq s \leq T_{ij}} |U_{ijt} - U_{ijs} | \geq \alpha %(T_{ij})^{-1}$ for any $i \in [n]$ and $j \in [p].$
\end{assumption}

\begin{assumption}
	\label{cond_kernel}
	(i) $K(\cdot)$ is a symmetric probability density function defined on $[-1,1]$ with $\int u^2K(u)du <\infty$ and $\int K(u)^2 du < \infty.$
	(ii) $K(\cdot)$ is Lipschitz continuous: there exists some positive constant $L$ such that
	$|K(u) -K(v)| \leq L|u-v|$ for any $u,v \in [-1,1].$
\end{assumption}

%\begin{assumption}
%\label{cond_bandwidth}
%{\color{red}(i) $n \widebar T_{\mu,j} h_{\mu,j} \rightarrow \infty$ as $n\rightarrow \infty;$
%(ii) $n \widebar T_{\sSigma,jk}^2 h_{\mu,j}^2 \rightarrow \infty$ as $n\rightarrow \infty$}
%\end{assumption}

\begin{assumption}
	\label{cond_cov_der}
	(i) $\partial^2 \mu_{j}(u)/\partial u^2$ is uniformly bounded over $u \in \cU$ and $j \in [p];$
    (ii) $\partial^2 \Sigma_{jk}(u,v)\\/\partial u^2,$ $\partial^2 \Sigma_{jk}(u,v)/\partial u \partial v,$ and $\partial^2 \Sigma_{jk}(u,v)/\partial v^2$ are uniformly bounded over $(u,v) \in \cU^2$ and $(j,k) \in [p]^2.$
\end{assumption}

%\begin{condition}
%	\label{cond_cov_positivity}
%	(i) strictly positive definite {\color{red} maybe not needed}
%\end{condition}

The sub-Gaussianities in Assumption~\ref{cond_subG} for both Hilbert space-valued random elements $X_{ij}(\cdot)$'s and random errors $\varepsilon_{ij}$'s together imply that the observations $Y_{ijt}$'s in (\ref{model}) are sub-Gaussian, which plays a crucial role in deriving our subsequent concentration inequalities.
The dense case in Assumption~\ref{cond_T} corresponds to a common practical scenario, where the sampling frequencies $T_{ij}$'s are of the same order across $i \in [n]$ and $j \in [p]$. Under such a scenario, Assumption~\ref{cond_weight} is automatically satisfied by two frequently-used weighting schemes including ``equal weight per observation'' and ``equal weight per subject''.
%Assumption~\ref{cond_T_bd} is standard in functional data analysis literature \cite[]{yao2005a,qiao2020}.
Assumption~\ref{cond_timepoints} means that all observational time points are distributed in the sense of ``uniformly on $\cU$". This prevents the occurrence of an extreme case where a large number of time points are concentrated in some small areas while leaving too few points in other regions. %to avoid the extreme case where many time points are located in some neighborhood and few points are in other places. {\color{red} to be added by Shaojun} %{\color{red} To Shaojun: add some remark for Assumption~\ref{cond_timepoints}}.
Assumptions~\ref{cond_T_bd}, \ref{cond_kernel} and \ref{cond_cov_der}
are standard in the literature of local linear smoothing for functional data \cite[]{yao2005a,zhang2016} adaptable to the multivariate setting.
%{\color{blue} \begin{condition}\label{same_order}
%Under the dense design, there exists a constant $ K_5 >0 $, such that
%\begin{equation}
%\frac{\max_{i,j} T_{ij}}{\min_{i,j} T_{ij}} \lesssim K_5\nonumber,
%\end{equation}
%\begin{equation}
%n(1 \wedge T_{ij}h) \geq 1.
%			\nonumber
%e\end{equation}
%\end{condition}
%I think this condition can be implied, to be discussed with Yizhu.
%}

\begin{theorem}\label{thm_mean_coneq}
Suppose that Assumptions~\ref{cond_subG}-\ref{cond_kernel} hold. For each $j \in [p],$ let $\gamma_{n,T,h,j}=n (1 \wedge \widebar{T}_{\mu,j} h_{\mu,j})$ with the corresponding average sampling frequency per subject $\widebar{T}_{\mu,j} = n^{-1} \sum_{i=1}^n T_{ij},$ then there exist some positive constants $c_1, c_2$ (independent of $n, p, \widebar{T}_{\mu,j}$'s) and arbitrarily small $\epsilon_1>0$ such that for any $\delta \in (0, 1],$

%(i) under the sparse design,
%\begin{equation}
%\label{mean_coneq_L2_sparse}
%\pr\Big(\|\hat \mu_{j}-\tilde \mu_{j}\| \geq \delta\Big) \leq c_2 \exp\Big(-\frac{c_1 n h_{\mu, j}\delta^2}{1+\delta}\Big),
%\end{equation}

%\begin{equation}
%\label{mean_coneq_sup_sparse}
%\pr\Big\{\sup_{u \in \cU} \big|\hat \mu_{j}(u)-\tilde \mu_{j}(u)\big| \geq \delta\Big\} \leq \frac{c_2 n}{h_{\mu, j}^{1/2}} \exp\Big(-\frac{c_1 n h_{\mu, j} \delta^2}{1+\delta}\Big),
%\end{equation}

%(ii) and, under the dense design with $\nu_{n,T,h,j}=n \wedge \sum_{i=1}^n T_{ij} h_{\mu,j},$

\begin{equation}
\label{mean_coneq_L2}
\pr\big(\|\hat \mu_{j}-\tilde \mu_{j}\|_2 \geq \delta\big) \leq c_2 \exp\big(-{c_1 \gamma_{n,T,h,j} \delta^2}\big),
\end{equation}
\begin{equation}
\label{mean_coneq_sup}
\pr\Big\{\sup_{u \in \cU} \big|\hat \mu_{j}(u)-\tilde \mu_{j}(u)\big| \geq \delta\Big\} \leq
\frac{c_2  (n^{\epsilon_1} \gamma_{n,T,h,j})^{1/2}}{h_{\mu,j}^2}
\exp\big(-{c_1 \gamma_{n,T,h,j} \delta^2}\big),
\end{equation}
where $\tilde \mu_j(u)$ is a deterministic univariate function that converges to $\mu_j(u)$ as $h_{\mu,j} \rightarrow 0.$ See (\ref{tilde_mu}) in Appendix~\ref{ap.pf} for the exact form of $\tilde \mu_j(u).$
\end{theorem}

%Next, we provide the concentration behavior of $\widehat \Sigma_{jk}.$ To facilitate representation, we suppose that $\mu_j(\cdot)$ is known in advance.
\begin{theorem}
\label{thm_cov_coneq}
%Suppose that Conditions~\ref{cond_subG}-\ref{cond_kernel} hold.
%Suppose that the assumptions in Theorem~\ref{thm_mean_coneq} hold.
Suppose that Assumptions~\ref{cond_subG}-\ref{cond_kernel} hold.
For each $j,k \in [p],$  let $\nu_{n,T,h,jk}=n(1 \wedge \widebar T_{\sSigma,jk}^2 h_{\sSigma,jk}^2)$ with the corresponding average sampling frequency per subject  being $\widebar T_{\sSigma,jk}= [n^{-1} \sum_{i=1}^n T_{ij}\{T_{ik}-I(j=k)\}]^{1/2},$ then there exist some positive constants $c_3, c_4$ (independent of $n, p, \widebar T_{\sSigma,jk}$'s) and arbitrarily small $\epsilon_2>0$ such that for any $\delta \in (0,1],$

%(i) under the sparse design,
%\begin{equation}
%\label{cov_coneq_L2_sparse}
%\pr\Big(\|\widehat C_{jk}-\widetilde C_{jk}\|_{\cS} \geq \delta\Big) \leq c_4 \exp\Big(-\frac{c_3 n h_{C,jk}^2\delta^2}{1+\delta}\Big),
%\end{equation}

%\begin{equation}
%\label{cov_coneq_sup_sparse}
%\pr\Big\{\sup_{(u,v) \in \cU^2} \big|\widehat C_{jk}(u,v)-\widetilde C_{jk}(u,v)\big| \geq \delta\Big\} \leq c_4 n^2 \exp\Big(-\frac{c_3 n h_{C,jk}^2\delta^2}{1+\delta}\Big),
%\end{equation}

%(ii) and, under the dense design with $\gamma_{n,T,h,jk}=n \wedge \sum_{i=1}^n T_{ij} T_{ik} h_{C,jk}^2,$
\begin{equation}
\label{cov_coneq_L2}
\pr\Big(\|\widehat \Sigma_{jk}-\widetilde \Sigma_{jk}\|_{\cS} \geq \delta\Big) \leq c_4 \exp\big(-{c_3\nu_{n,T,h,jk}\delta^2}\big),
\end{equation}

\begin{equation}\label{cov_coneq_sup}
\pr\Big\{\sup_{(u,v) \in \cU^2} \big|\widehat \Sigma_{jk}(u,v)-\widetilde \Sigma_{jk}(u,v)\big| \geq \delta\Big\} \leq \frac{c_4 n^{\epsilon_2} \nu_{n,T,h,jk}}{h_{\sSigma,jk}^6}\exp\big(-{c_3\nu_{n,T,h,jk}\delta^2}\big),
\end{equation}
where $\widetilde \Sigma_{jk}(u,v)$ is a deterministic bivariate function that converges to $\Sigma_{jk}(u,v)$ as $h_{\sSigma,jk} \rightarrow 0.$ See (\ref{Sigma_tilde}) in Appendix~\ref{ap.pf} for the exact form of $\widetilde \Sigma_{jk}(u,v).$
\end{theorem}

\begin{remark}
\label{rmk_asy_rates}
The concentration inequalities in Theorems~\ref{thm_mean_coneq} and \ref{thm_cov_coneq} imply that $\hat \mu_j$ and $\widehat \Sigma_{jk}$ are nicely concentrated around $\tilde \mu_j$ and $\widetilde \Sigma_{jk}$, respectively, in both $L_2$ norm and supremum norm with generalized sub-Gaussian-type tail behaviors. %depending on which term in $\delta^2 \wedge \delta$ is dominant.
%The concentration inequalities in Theorems~\ref{thm_mean_coneq} and \ref{thm_cov_coneq} imply that the tails for $\hat \mu_j -\tilde \mu_j$ and $\widehat C_{jk} -\widetilde C_{jk}$ in $L_2$ and supremum norms behave in a generalized sub-Gaussian or sub-exponential way depending on which term in $\delta^2 \wedge \delta$ is dominant.
It is worth mentioning that such $L_2$ and uniform concentration results are derived based on the local concentration inequalities of $\hat \mu_j(u)$
and $\widehat \Sigma_{j,k}(u,v)$ for fixed interior points $u,v \in \cU,$ which enjoy the same tail behaviors as (\ref{mean_coneq_L2}) and (\ref{cov_coneq_L2}).
Besides being fundamental to derive elementwise maximum error bounds that are essential for further convergence analysis under high-dimensional settings, these non-asymptotic results lead to the same rates of %optimal minimax rates of
$L_2$ convergence and uniform convergence compared to those in \cite{zhang2016}. Specifically, under extra Assumption~\ref{cond_cov_der}, it holds that
$$\|\hat \mu_j - \mu_j\|_2=O_\p\big\{n^{-1/2} + (n\widebar T_{\mu,j} h_{\mu,j})^{-1/2}  + h_{\mu,j}^2\big\},$$
$$\sup_{u \in \cU}\big|\hat \mu_j(u) - \mu_j(u)\big|=O_\p\big[(\log n)^{1/2} n^{-1/2}\{1 + (\widebar T_{\mu,j} h_{\mu,j})^{-1/2}\} + h_{\mu,j}^2\big],$$
$$\|\widehat \Sigma_{jk} - \Sigma_{jk}\|_{\cS}=O_\p\big\{n^{-1/2} + (n\widebar T_{\sSigma,jk}^2 h_{\sSigma,jk}^2)^{-1/2} + h_{\sSigma,jk}^2\big\},$$
$$\sup_{(u,v) \in \cU^2}\big|\widehat \Sigma_{jk}(u,v) - \Sigma_{jk}(u,v)\big|=O_\p\big[(\log n)^{1/2}n^{-1/2}\{1 +(\widebar T_{\sSigma,jk}^2 h_{\sSigma,jk}^2)^{-1/2} \} + h_{\sSigma,jk}^2\big].$$
\end{remark}
%Provided thes rates of $\hat \mu_{j}$ (or $\widehat C_{jk}$),
\begin{remark}
\label{rmk_asy_phase}
The above rates of convergence reveal interesting phase transition phenomena depending on the ratio of the average sampling frequency per subject $\widebar T_{\mu,j}$ (or $\widebar T_{\sSigma,jk}$) to $n^{1/4}.$ In the following, we use different rates of $L_2$  convergence for $\hat \mu_{j}$ and $\widehat \Sigma_{jk}$ to illustrate a systematic partition of partially observed functional data into three categories:

\begin{enumerate}%[label=\Roman*]
\item\label{case_a} Under the sparse design, %with bounded $T_{ij}$'s,
when $h_{\mu,j} \asymp n^{-1/5},$ $$\|\hat \mu_{j} - \mu_{j}\|_2=O_\p(n^{-1/2} h_{\mu,j}^{-1/2} + h_{\mu,j}^2)=O_\p(n^{-2/5});$$ when $h_{\sSigma,jk} \asymp n^{-1/6},$
$$\|\widehat \Sigma_{jk} - \Sigma_{jk}\|_{\cS}=O_\p(n^{-1/2} h_{\sSigma,jk}^{-1} + h_{\sSigma,jk}^2)=O_\p(n^{-1/3}).$$

\item\label{case_b} Under the dense design,
when $\widebar T_{\mu,j} n^{-1/4} \rightarrow 0$ with $h_{\mu,j} \asymp (n \widebar T_{\mu,j})^{-1/5},$
$$\|\hat \mu_{j} - \mu_{j}\|_2=O_\p(n^{-1/2} \widebar T_{\mu,j}^{-1/2} h_{\mu,j}^{-1/2} + h_{\mu,j}^2)=O_\p\{(n \widebar T_{\mu,j})^{-2/5}\};$$
when $\widebar T_{\sSigma,jk} n^{-1/4} \rightarrow 0$ with $h_{\sSigma,j} \asymp (n \widebar T_{\sSigma,jk}^2)^{-1/6},$
$$\|\widehat \Sigma_{jk} - \Sigma_{jk}\|_{\cS}=O_\p(n^{-1/2} \widebar T_{\sSigma,jk}^{-1} h_{\sSigma,jk}^{-1} + h_{\sSigma,jk}^2)=O_\p\{(n \widebar T_{\sSigma,jk}^2)^{-1/3}\}.$$

\item\label{case_c} Under the dense design,
when $\widebar T_{\mu,j} n^{-1/4} \rightarrow \tilde c$ (some positive constant) %($c \in (0, \infty)$)
with $h_{\mu,j} \asymp n^{-1/4}$ or $\widebar T_{\mu,j} n^{-1/4} \rightarrow \infty$ with $h_{\mu,j} = o(n^{-1/4})$ and $\widebar T_{\mu,j}h_{\mu,j} \rightarrow \infty,$
$$\|\hat\mu_{j} - \mu_{j}\|_2=O_\p(n^{-1/2});$$
when $\widebar T_{\sSigma,jk} n^{-1/4} \rightarrow \tilde c$ with $h_{\sSigma,jk} \asymp n^{-1/4}$
or $\widebar T_{\sSigma,jk} n^{-1/4} \rightarrow \infty$ with $h_{\sSigma,jk}=o(n^{-1/4})$ and $\widebar T_{\sSigma,jk}h_{\sSigma,jk}  \rightarrow \infty,$
$$\|\widehat \Sigma_{jk} - \Sigma_{jk}\|_{\cS}=O_\p(n^{-1/2}).$$
\end{enumerate}

As $\widebar T_{\mu,j}$ and $\widebar T_{\sSigma,jk}$ grow very fast, case~\ref{case_c} results in the root-$n$ rate complying with the parametric rate for fully observed functional data.
As $\widebar T_{\mu,j}$ and $\widebar T_{\sSigma,jk}$ grow moderately fast, case~\ref{case_b} corresponds to the optimal minimax rates \cite[]{zhang2016}, which are slower than root-$n$ but faster than the counterparts for sparsely observed functional data.
Our established convergence rates in cases~\ref{case_a}, \ref{case_b} and \ref{case_c} allow free choices of $(j,k),$ and are respectively consistent to those of the mean and covariance estimators under categories of ``sparse'', ``semi-dense'' and ``ultra-dense'' univariate functional data introduced in \cite{zhang2016}. %In a similar fashion to the above discussion, the rates for estimated mean functions under these three categories are also respectively of the same order compared to those in \cite{zhang2016}.
\end{remark}

%\begin{condition}
%\label{cond_bandwidth}
%(i) $\log p \{\min_j (n \wedge \sum_{i=1}^n T_{ij} h_{\mu,j})\}^{-1} \rightarrow 0$ and $\max_j h_{\mu,j} \rightarrow 0;$
%(ii) $\log p (\min_{j,k} [n \wedge \sum_{i=1}^n T_{ij}\{T_{ik}-I(j=k)\}h_{C,jk}^2])^{-1} \rightarrow 0$ and
%$\max_{j,k} h_{C,jk} \rightarrow 0.$
%\end{condition}

\begin{theorem}\label{thm_mean_maxrate}
Suppose that the assumptions in Theorem~\ref{thm_mean_coneq} and Assumption~\ref{cond_cov_der}(i) hold, and
$ (\min_j \gamma_{n,T,h,j})^{-1}\log p \rightarrow 0,$ $\max_j h_{\mu,j} \rightarrow 0$ as $n,p \rightarrow \infty.$
%Let Assumptions~\ref{cond_subG}--\ref{cond_kernel}, \ref{cond_cov_der}(i) hold and $\log p \{\min_j (n \wedge \sum_{i=1}^n T_{ij} h_{\mu,j})\}^{-1} \rightarrow 0,$ $\max_j h_{\mu,j} \rightarrow 0.$
It then holds that

%(i) under the sparse design,

%\begin{equation}
%\label{mean_max_L2_sparse}
%  \max_{j \in [p]}\|\widehat \mu_{j}-\mu_{j}\| = O_\p\left\{\Big(\frac{\log p}{n \min_j h_{\mu,j}}\Big)^{1/2} + \max_j h_{\mu,j}^2\right\},
%\end{equation}

%\begin{equation}
%\label{mean_max_sup_sparse}
%  \max_{j \in [p]} \sup_{u \in \cU} \big|\hat \mu_{j}(u)-\mu_{j}(u)\big| = O_\p\left[\Big\{\frac{\log (p \vee n)}{n \min_j h_{\mu,j} }\Big\}^{1/2} + \max_j h_{\mu,j}^2\right]
%\end{equation}
%with $\min_j h_{\mu,j} \gtrsim n^{-\eta}$ for some constant $\eta>0,$\\

%(ii) and, under the dense design,

\begin{equation}
\label{mean_max_L2_dense}
  \max_{j\in [p]}\|\hat \mu_{j}-\mu_{j}\|_2 = O_\p\left\{\Big(\frac{\log p}{\min_{j}\gamma_{n,T,h,j}}\Big)^{1/2} + \max_j h_{\mu,j}^2\right\},
\end{equation}
and, if $\min_j h_{\mu,j} \asymp \{\log (p\vee n)/n\}^{\kappa_1}$ for some $\kappa_1 \in (0, 1/2],$
\begin{equation}
\label{mean_max_sup_dense}
  \max_{j \in [p]} \sup_{u \in \cU} \big|\hat \mu_{j}(u)-\mu_{j}(u)\big| = O_\p\left[\Big\{\frac{\log (p \vee n)}{\min_{j}\gamma_{n,T,h,j}}\Big\}^{1/2} + \max_j h_{\mu,j}^2\right].
\end{equation}
%with $\min_j h_{\mu,j} \gtrsim n^{-\eta}$ for some constant $\eta>0.$
\end{theorem}

\begin{theorem}\label{thm_cov_maxrate}
Suppose that the assumptions in Theorem~\ref{thm_cov_coneq} and Assumption~\ref{cond_cov_der}(ii) hold, and
%Let Assumptions~\ref{cond_subG}--\ref{cond_kernel}, \ref{cond_cov_der}(ii) hold and $\log p (\min_{j,k} [n \wedge \sum_{i=1}^n T_{ij}\{T_{ik}-I(j=k)\}h_{\sSigma,jk}^2])^{-1} \rightarrow 0,$
$ (\min_{j,k} \nu_{n,T,h,jk})^{-1}\log p \rightarrow 0,$
$\max_{j,k} h_{\sSigma,jk} \rightarrow 0$ as $n,p \rightarrow \infty.$ It then holds that

%(i) under the sparse design,
%\begin{equation}
%\label{cov_max_L2_sparse}
%  \max_{j,k \in [p]}\|\widehat C_{jk}-C_{jk}\|_{\cS} = O_\p\left\{\Big(\frac{\log p}{\min_{j,k} h_{C,jk}^2}\Big)^{1/2} + \max_{j,k} h_{C,jk}^2\right\},
%\end{equation}

%\begin{equation}
%\label{cov_max_sup_sparse}
%  \max_{j,k \in [p]} \sup_{(u,v) \in \cU^2} \big|\widehat C_{jk}(u,v)-C_{jk}(u,v)\big| = O_\p\left[\Big\{\frac{\log (p \vee n)}{\min_{j,k} h_{C,jk}^2}\Big\}^{1/2} + \max_{j,k} h_{C,jk}^2\right]
%\end{equation}
%with $\min_{j,k} h_{C,jk} \gtrsim n^{-\tau}$ for some constant $\tau>0,$\\

%(ii) under the dense design,
\begin{equation}
\label{cov_max_L2_dense}
  \max_{j,k \in [p]}\|\widehat \Sigma_{jk}-\Sigma_{jk}\|_{\cS} = O_\p\left\{\Big(\frac{\log p}{\min_{j,k} \nu_{n,T,h,jk}}\Big)^{1/2} + \max_{j,k} h_{\sSigma,jk}^2\right\},
\end{equation}
and, if $\min_{j,k} h_{\sSigma,jk}\asymp \{\log (p\vee n)/n\}^{\kappa_2}$ for some $\kappa_2 \in (0, 1/2],$
\begin{equation}
\label{cov_max_sup_dense}
  \max_{j,k \in [p]} \sup_{(u,v) \in \cU^2} \big|\widehat \Sigma_{jk}(u,v)-\Sigma_{jk}(u,v)\big| = O_\p\left[\Big\{\frac{\log (p \vee n)}{\min_{j,k} \nu_{n,T,h,jk}}\Big\}^{1/2} + \max_{j,k} h_{\sSigma,jk}^2\right].
\end{equation}
%with $\min_{j,k} h_{C,jk} \gtrsim n^{-\tau}$ for some constant $\tau>0.$
\end{theorem}

We observe that the elementwise maximum rates of $L_2$ convergence and uniform convergence are governed by both dimensionality parameters ($n, p, \{\widebar T_{\mu,j}\}_{j \in [p]},\{\widebar T_{\sSigma,jk}\}_{j,k \in [p]}$) and internal parameters ($\{h_{\mu,j}\}_{j \in [p]}, \{h_{\sSigma,jk}\}_{j,k \in [p]}$). Each convergence rate is composed of two terms reflecting our familiar variance-bias tradeoff in nonparametric statistics. It is easy to see that the variance terms are determined by the least frequently sampled and smoothed components, that is the smallest $\widebar T_{\mu,j}$ (or $\widebar T_{\sSigma,jk}$) and $h_{\mu,j}$ (or $h_{\sSigma,jk}$) across $j,k,$ whereas the highest level of smoothness with the largest $h_{\mu,j}$ (or $h_{\sSigma,jk}$) controls the bias terms.

\begin{remark}
\label{rmk_nonasy_phase}
To facilitate further discussion, we consider the simplified setting where $\widebar T_{\mu,j} \asymp \widebar T_{\mu}, h_{\mu,j} \asymp h_{\mu}$ and $\widebar T_{\sSigma,jk} \asymp \widebar T_{\sSigma}, h_{\sSigma,jk} \asymp h_{\sSigma}$ for each $j, k.$ Compared to cases \ref{case_a}--\ref{case_c} above, the corresponding elementwise maximum rates of convergence for $\{\hat \mu_j\}_{j \in [p]}$ (or $\{\widehat \Sigma_{jk}\}_{j,k \in [p]}$) in Theorem~\ref{thm_mean_maxrate} (or Theorem~\ref{thm_cov_maxrate}) reveal scaled phase transitions for dense functional data depending on the relative order of $\widebar T_{\mu}$ (or $\widebar T_{\sSigma}$) to $n^{1/4} (\log p)^{-1/4}$ instead of $n^{1/4}.$ In the following, we use elementwise maximum rates of $L_2$ convergence to illustrate the phase transition phenomena and the optimal estimation from sparse to dense functional data in high dimensions. In terms of uniform convergence, the same phenomena occur as long as $p \gtrsim n.$

\begin{enumerate}[label=(\roman*)]
\item \label{case_i} Under the sparse design, %with $\max_{i,j} T_{ij} \leq T_0 <\infty,$
when $h_{\mu} \asymp (\log p)^{1/5} n^{-1/5},$

\begin{equation*}\label{mean_max_sparse_optimal}
\max_{j}\|\hat \mu_{j} - \mu_{j}\|_2=O_\p\left\{ \Big(\frac{\log p}{n  h_{\mu}}\Big)^{1/2} + h_{\mu}^2\right\} = O_\p \left\{\Big(\frac{\log p}{n}\Big)^{2/5}\right\};
\end{equation*}
when $h_{\sSigma} \asymp (\log p)^{1/6} n^{-1/6},$
\begin{equation*}\label{cov_max_sparse_optimal}
\max_{j,k}\|\widehat \Sigma_{jk} - \Sigma_{jk}\|_{\cS}  =O_\p\left\{ \Big(\frac{\log p}{n h_{\sSigma}^2}\Big)^{1/2} + h_{\sSigma}^2\right\} = O_\p \left\{\Big(\frac{\log p}{n}\Big)^{1/3}\right\}.
\end{equation*}
\item\label{case_ii} Under the dense design, when $\widebar T_{\mu} (\log p)^{1/4} n^{-1/4} \rightarrow 0$ with $h_{\mu} \asymp (\log p)^{1/5} (n \widebar T_{\mu})^{-1/5},$
\begin{equation}
\label{mean_max_dense_optimal}
\max_{j}\|\hat \mu_{j} - \mu_{j}\|_2=O_\p\left\{ \Big(\frac{\log p}{n \widebar T_{\mu} h_{\mu}}\Big)^{1/2} + h_{\mu}^2\right\} = O_\p \left\{\Big(\frac{\log p}{n \widebar T_{\mu}}\Big)^{2/5}\right\};
\end{equation}
when $\widebar T_{\sSigma} (\log p)^{1/4} n^{-1/4} \rightarrow 0$ with $h_{\sSigma} \asymp (\log p)^{1/6} (n \widebar T_{\sSigma}^2)^{-1/6},$
\begin{equation}
\label{cov_max_dense_optimal}
\max_{j,k}\|\widehat \Sigma_{jk} - \Sigma_{jk}\|_{\cS}  =O_\p\left\{ \Big(\frac{\log p}{n \widebar T_{\sSigma}^2 h_{\sSigma}^2}\Big)^{1/2} + h_{\sSigma}^2\right\} = O_\p \left\{\Big(\frac{\log p}{n \widebar T_{\sSigma}^2}\Big)^{1/3}\right\}.
\end{equation}

\item\label{case_iii} Under the dense design, when $\widebar T_{\mu} (\log p)^{1/4} n^{-1/4} \rightarrow \tilde c$ with $h_{\mu} \asymp (\log p)^{1/4} n^{-1/4}$
or $\widebar T_{\mu} (\log p)^{1/4} n^{-1/4} \rightarrow \infty$ with $h_{\mu}=o\{(\log p)^{1/4} n^{-1/4}\}$ and $\widebar T_{\mu} h_{\mu} \rightarrow \infty,$ $$\max_j \|\hat \mu_j -\mu_j\|_2=O_\p\left\{ \Big(\frac{\log p}{n}\Big)^{1/2}\right\};$$
when $\widebar T_{\sSigma} (\log p)^{1/4} n^{-1/4} \rightarrow \tilde c$ with $h_{\sSigma} \asymp (\log p)^{1/4} n^{-1/4}$
or $\widebar T_{\sSigma} (\log p)^{1/4} n^{-1/4} \rightarrow \infty$ with $h_{\sSigma}=o\{(\log p)^{1/4}n^{-1/4}\}$ and $\widebar T_{\sSigma}h_{\sSigma}  \rightarrow \infty,$
$$\max_{j,k}\|\widehat \Sigma_{jk} - \Sigma_{jk}\|_{\cS} =O_\p\left\{ \Big(\frac{\log p}{n}\Big)^{1/2}\right\}.$$
\end{enumerate}
\end{remark}

\begin{remark}
\label{rmk_nonasy_phase2}
In a similar spirit to the partitioned three categories for univariate functional data (see cases~\ref{case_a}, \ref{case_b} and \ref{case_c} above), we can also term the high-dimensional partially observed functional data in cases~\ref{case_i}, \ref{case_ii}, and \ref{case_iii} as ``sparse", ``semi-dense", and ``ultra-dense", respectively. The main difference lies in the presence of
additional $\log p$ terms to account for the high-dimensional effect.
\begin{itemize}
\item As $\widebar T_{\mu}$ and $\widebar T_{\sSigma}$ grow at least in the order of $n^{1/4} (\log p)^{-1/4},$ the attained optimal rate $(\log p)^{1/2} n^{-1/2}$ is identical to that for the fully observed functional data \cite[]{zapata2022}, presenting that the theory for high-dimensional ultra-dense functional data falls in the parametric paradigm.
\item As $\widebar T_{\mu}$ and $\widebar T_{\sSigma}$ diverge slower than $n^{1/4} (\log p)^{-1/4},$ if we let $h_{\mu} \asymp (\log p)^{1/5} (n \widebar T_{\mu})^{-1/5}$ and $h_{\sSigma} \asymp (\log p)^{1/6} (n \widebar T_{\sSigma}^2)^{-1/6}$ to balance the corresponding variance and bias terms, the optimal rates for high-dimensional semi-dense functional data are respectively achieved in (\ref{mean_max_dense_optimal}) and (\ref{cov_max_dense_optimal}). These rates degenerate to the minimax rates in case~\ref{case_b} when $p$ is fixed.
With the choice of elementwise optimal bandwidths $h_{\mu} \asymp (n \widebar T_{\mu})^{-1/5}$ and $h_{\sSigma} \asymp (n \widebar T_{\sSigma}^2)^{-1/6},$ we obtain
$$\max_{j}\|\hat \mu_{j} - \mu_{j}\|_2=O_\p\{(\log p)^{1/2} (n \widebar T_{\mu})^{-2/5}\},~\max_{j,k}\|\widehat \Sigma_{jk} - \Sigma_{jk}\|=O_\p\{(\log p)^{1/2} (n \widebar T_{\sSigma})^{-1/3}\},$$ which are respectively slower than the optimal rates in (\ref{mean_max_dense_optimal}) and (\ref{cov_max_dense_optimal}).
Such discussion applies analogously to the sparse functional setting. See cases~\ref{case_a} and \ref{case_i}.
\item Compared to the asymptotic results for cases~\ref{case_a}, \ref{case_b} and \ref{case_c} under a fixed $p$ scenario, the high-dimensionality in cases~\ref{case_i}, \ref{case_ii} and \ref{case_iii} leads to the scaled phase transitions,
%$\log p$-based multiplicative shifts of the phase transitions,
optimal selected bandwidths, and corresponding optimal rates, each of which is up to a factor of $\log p$ at some polynomial order.
%a polynomial-of-$\log p$ factor.
%\item It is worth mentioning that \cite{lee2021} simply assumed the elementwise maximum rate $$\max_{j,k}\|\widehat \Sigma_{jk} - \Sigma_{jk}\|_{\cS}=O_\p(\log p n^{-\tau}),$$ which was further reduced to $O_\p\{(\log p)^{1/2} n^{-\tau}\}$ in \cite{fang2022}. The involved quantity $\tau \in (0,1/2]$ reflects the average sampling frequency with larger values yielding denser observational points. We fundamentally improve their rates in the sense of precisely specifying the largest values of $\tau$ under cases with different sampling frequencies.
\end{itemize}
\end{remark}

\section{Applications}
\label{sec:app}
In this section, we outline three applications of our established non-asymptotic results for the local linear smoothers under high dimensional settings. %when handling high-dimensional partially observed functional data.

%\begin{itemize}
%\item FPCA, score, convergence rate for $\hat X_{ij}-X_{ij}(u)$ and its application kernel covariance  (by Shaojun)
%\item Sparse FPCA (by Shaojun), Xue's paper (by Xinghao)
%\item Adaptive thresholding to be done by Xinghao
%\end{itemize}

\subsection{Estimation Under the FPCA Framework}
\label{sec_fpca}
A standard procedure towards the estimation of models involving high-dimensional functional data consists of two or three steps. 
Due to the infinite-dimensionality of functional data, the first step performs dimension reduction via, e.g., FPCA, to approximate each $X_{ij}(\cdot)$ by the $d_j$-dimensional truncation. This effectively transforms the problem of modeling the $p$-vector of functional variables into that of modeling the $(\sum_{j=1}^p d_j)$-vector of FPC scores. 
To overcome the difficulties caused by high-dimensionality, some functional sparsity assumptions are commonly imposed, which results in the estimation under block sparsity constraints in the second step. Examples include the group-lasso penalized least squares estimation in regression setups \cite[]{fan2015,kong2016,wang2022}, the group graphical lasso in functional graphical model estimation \cite[]{qiao2019,solea2022,zapata2022} and other related applications mentioned in Section~\ref{sec:intro}. 
Finally, the third step re-transforms block sparse estimates obtained in the second step to functional sparse estimates via estimated principal component functions obtained in the first step. In functional graphical model estimation, the third step is no longer required.
Building upon established theoretical results, this section presents some non-asymptotic results within the FPCA framework, which are crucial not only in their own right but also in providing the theoretical support for such multi-step estimation procedure.

For each $j \in [p],$ the standard dimension reduction method performs Karhunen-Lo\`eve expansion of each target trajectory $X_{ij}(\cdot)$ and truncates the expansion to the first $d_j$ terms, which serves as the foundation of FPCA:
\begin{equation}
\label{fpca.expand}
X_{ij}(\cdot) = \mu_j(\cdot) +\sum_{l=1}^{\infty}\xi_{ijl}\phi_{jl}(\cdot) \approx \mu_j(\cdot) 
 + \sum_{l=1}^{d_j}\xi_{ijl}\phi_{jl}(\cdot),
\end{equation}
where the coefficients $\xi_{ijl}=\langle X_{tj}-\mu_j, \phi_{jl}\rangle$ for $l\geq 1,$ namely FPC scores, correspond to a sequence of random variables with ${\mathbb E}(\xi_{ijl})=\lambda_{jl}$ and $\cov(\xi_{ijl}, \xi_{ijl'})=\lambda_{jl}I(l \neq l')$ and $\lambda_{j1} \geq \lambda_{j2} \geq \cdots >0$ are the eigenvalues of $\Sigma_{jj}(u,v)$ and $\phi_{j1}(\cdot), \phi_{j2}(\cdot), \dots$ are the corresponding eigenfunctions. To implement FPCA based on partially observed functional data, we carry out an eigenanalysis of the local linear smoothers $\widehat\Sigma_{jj}(u,v)$ that leads to estimated eigenvalue/eigenfunction pairs $(\hat \lambda_{jl}, \widehat \phi_{jl}(\cdot))$ for $j \in [p]$ and $l \in [d_j].$ Based on Theorem~\ref{thm_cov_maxrate}, we can establish the elementwise maximum rates for estimated eigenvalues and eigenfunctions in the following proposition.

\begin{proposition}
\label{thm_fpca_maxrate}
Suppose that the assumptions in Theorem~\ref{thm_cov_maxrate} hold and $\lambda_{j1} > \lambda_{j2} > \cdots > 0$ for each $j \in [p].$  Let $\delta_{jl} = {\min}_{k \in [l]}\{\lambda_{jk} - \lambda_{j(k+1)}\}$ for $j \in [p]$ and $l \in [d_j]$ . Then we have
\begin{equation*}
  \label{bd_max_lam_phi}
\max_{j \in [p], l \in [d_j]} \left\{\big|\hat \lambda_{jl} - \lambda_{jl}\big| 
+ \delta_{jl} \big\| \widehat \phi_{jl} - \phi_{jl}\big\|_2 \right\}
= O_\p\left\{\Big(\frac{\log p}{\min_{j,k} \nu_{n,T,h,jk}}\Big)^{1/2} + \max_{j,k} h_{\sSigma,jk}^2\right\}.
\end{equation*}
\end{proposition}

Proposition~\ref{thm_fpca_maxrate} can be used to provide the theoretical guarantees for the first and third steps under high-dimensional settings. In the second step, the main target is to implement the block regularized estimation based on the estimated FPC scores. 
Under the dense design with $T_{ij} \rightarrow \infty$, the estimated FPC scores $\langle  X_{ij} - \hat \mu_j, \hat \phi_{jl}\rangle$ can be well approximated by the numerical integration based on $\{U_{ijt}, Y_{ijt}, \widehat \phi_{jl}(U_{ijt})\}$ for $t \in [T_{ij}].$ To be specific, we can employ a Trapezoid rule-based numerical integration to estimate FPC scores
\begin{equation}
\label{est.score.dense}
\widehat \xi_{ijl}^{(1)}=\sum_{t=2}^{T_{ij}} \frac{\big\{Y_{ij(t-1)} - \hat \mu_j(U_{ij(t-1)})\big\}\widehat\phi_{jl}(U_{ij(t-1)}) + \big\{Y_{ijt} - \hat\mu_{j}(U_{ijt})\big\} \widehat \phi_{jl}(U_{ijt})}{2}\big|U_{ijt} - U_{ij(t-1)}\big|.
\end{equation}

However, such numerical integration approach fails under the sparse design with $T_{ij} \leq T_0 <\infty.$ We instead employ the principal components analysis through conditional expectation (PACE) method \cite[]{yao2005a} to estimate FPC scores. For each $j \in [p],$ under the assumption that $\xi_{ijl}$ and $\varepsilon_{ijt}$ in (\ref{model}) are jointly Gaussian, the PACE estimation of the FPC scores for the $i$th subject given the data
from the individual reduces to the estimated conditional expectation 
%{\color{red}Write $\widetilde \bX_{ij}=(X_{ij}(U_{ij1}), \dots, X_{ij}(U_{ijT_{ij}})^\T,$ $\widetilde \bY_{ij}=(Y_{ij1}, \dots, Y_{ijT_{ij}})^\T,$ $\bmu_{ij}=(\mu_{j}(U_{ij1}), \dots, \mu_{j}(U_{ijT_{ij}}))^\T,$ and $\bphi_{ijl}=(\phi_{jl}(U_{ij1}), \dots, \phi_{jl}(U_{ijT_{ij}}))^\T.$ For each $j \in [p],$ the best prediction of the FPC scores for the $i$th subject given the datafrom the individual is the conditional expectation 
%$$
%\widetilde \xi_{ijl} = {\mathbb E}(\xi_{ijl}| \widetilde\bY_{ij}) = \lambda_{jl} \bphi_{ijl}^\T \bSigma_{\bY_{ij}}^{-1}(\widetilde\bY_{ij}-\bmu_{ij}),
%$$
%where $\bSigma_{\bY_{ij}}=\cov(\widetilde\bY_{ij}, \widetilde\bY_{ij})=\cov(\widetilde\bX_{ij},\widetilde\bX_{ij})+\sigma_j^2 \bI_{T_{ij}}.$ By substituting the estimates of $\lambda_{jl}, \bphi_{ijl}, \bSigma_{Y_{ij}}$ and $\bmu_{ij},$} we can obtain the estimator for $\widetilde \xi_{ijl}$ by

\begin{equation}
\label{est.score.sparse}
\widehat \xi_{ijl}^{(2)} = {\widehat{\mathbb E}}(\xi_{ijl}|\widetilde \bY_{ij}) =\hat\lambda_{jl} \widehat\bphi_{ijl}^\T \widehat\bSigma_{\bY_{ij}}^{-1}(\widetilde\bY_{ij}-\widehat\bmu_{ij}),
\end{equation}
where we write $\widetilde \bY_{ij}=(Y_{ij1}, \dots, Y_{ijT_{ij}})^\T,$ $\widehat\bSigma_{\bY_{ij}}$ is a $T_{ij} \times T_{ij}$ matrix with its $(t,t')$-th entry $(\widehat\bSigma_{\bY_{ij}})_{t,t'}=\widehat \Sigma_{jj}(U_{ijt}, U_{ijt'}) + \hat \sigma_j^2 I(t=t'),$
$\widehat\bphi_{ijl}=\{\hat\phi_{jl}(U_{ij1}), \dots, \hat\phi_{jl}(U_{ijT_{ij}})\}^\T,$ and 
$\widehat\bmu_{ij}=(\hat\mu_{j}(U_{ij1}), \dots, \hat\mu_{j}(U_{ijT_{ij}}))^\T.$ See \cite{yao2005a} for details on the estimate $\hat \sigma_j^2$ of $\sigma_j.$

For each $j,k \in [p], l \in [d_j]$ and $m \in [d_k],$ let $\sigma_{jklm}^{(h)}={\mathbb E}(\xi_{ijl}^{(h)}\xi_{ikm}^{(h)})$  and its sample estimator be $\widehat \sigma_{jklm}^{(h)}=n^{-1}\sum_{i=1}^n \widehat\xi_{ijl}^{(h)} \widehat \xi_{ikm}^{(h)}$ for $h=1$ (dense case) and $2$ (sparse case).
To theoretically support the second step, it is essential to establish the elementwise maximum rates for $\{\widehat \sigma_{jklm}^{(h)}\}$ with $h=1,2,$ i.e. the convergence rate of
$$\max_{j,k \in [p], l \in [d_j], m \in [d_k]} \frac{\Big|\widehat \sigma_{jklm}^{(h)} - \sigma_{jklm}^{(h)}\Big|}{\Lambda_{jklm}^{(h)}},$$ where $\Lambda_{jklm}^{(h)}$ represents some normalization term that may depend on $\delta_{jl}$ and $\delta_{km}.$ Our established elementwise maximum rate in Theorem~\ref{thm_cov_maxrate} may still be applicable in this context. However, this topic is beyond the scope of the current paper and is thus left for future research.
%While this topic is beyond the scope of the current paper, we believe our established elementwise maximum rates may still be applicable in this context and warrant further investigation in future research.

In addition, building upon the expansion in (\ref{fpca.expand}) and estimated FPC scores in (\ref{est.score.dense}) or (\ref{est.score.sparse}), we can estimate the target trajectories for prediction or subsequent modeling
\begin{equation}
\label{fpca.exp.est}
\widehat X_{ij}^{(h)}(\cdot) = \widehat \mu_j(\cdot) + \sum_{l=1}^{d_j} \widehat \xi_{ijl}^{(h)}  \widehat \phi_{jl}(\cdot), ~~~i\in [n], ~j \in [p], ~h=1,2.
\end{equation}
We leave the development of non-asymptotic results for (\ref{fpca.exp.est}) as future research. Such results could possibly be useful for the application of nonparametric functional graphical model estimation \cite[]{li2018} when handling the practical partially observed functional scenario in high dimensions.

\subsection{Sparse FPCA}
\label{sec_sfpca}
While componentwise FPCA in Section~\ref{sec_fpca} may fail to model the correlation between components in $X_{i1}(\cdot), \dots, X_{ij}(\cdot),$ more effective dimension reduction can be achieved by leveraging the correlation information across different components, such as multivariate FPCA \cite[]{happ2018} (for fixed $p$) and sparse FPCA \cite[]{hu2022}, which incorporates the notation of sparsity in multivariate statistics into the functional setting to accommodate high-dimensional functional data. The sparsity structure motivates a thresholding rule that is easy to compute by exploiting the relationship between the univariate orthonormal basis representation for infinite-dimensional processes and multivariate Karhunen-Lo\`eve expansion in the form of 
\begin{equation*}
\label{mfpca}
\bX_{i}(\cdot) = \bmu(\cdot) +\sum_{k=1}^{\infty} \zeta_{ik} \bpsi_{k}(\cdot),
\end{equation*}
where $\zeta_{ik} = \sum_{j=1}^p \langle X_{ij} - \mu_j, \psi_{kj} \rangle,$ $\tilde{\lambda}_k$ and $\bpsi_k(\cdot) = \{\psi_{k1}(\cdot),\ldots,\psi_{kp}(\cdot)\}^\T$ for $k=1, \dots, \infty$ are the eigenvalues and the corresponding eigenfunctions of $\bSigma(u,v),$ respectively, satisfying $\int_{\cU} \bSigma(u,v) \bpsi_k(v) {\rm d} v = \tilde \lambda_k \bpsi_k(u).$
%An alternative approach to consider is multivariate FPCA, which aids in the interpretation of high-dimensional functional data. Recently, \cite{hu2022} introduced an innovative sparse FPCA framework that incorporates sparsity into functional variables. They clarified the relationship between the multivariate Karhunen-Lo$\grave{e}$ve (K-L) expansion and the univariate orthonormal basis representation for infinite-dimensional processes. Furthermore, they proposed a thresholding technique to pinpoint significant processes. 
Specifically, provided with a complete and orthonormal basis  $\{b_l(\cdot):l=1, \dots, \infty\}$, each random process is represented as $X_{ij}(\cdot)  = \mu_j(\cdot) + \sum_{l = 1}^{\infty} \theta_{ijl} b_l(\cdot),$ where the basis coefficients $\theta_{ijl} = \langle X_{ij} -\mu_j,b_l\rangle.$ Denote $\psi_{kj}(\cdot) = \sum_{l=1}^\infty \eta_{kjl} b_l(\cdot),$ where $\eta_{kjl}=\langle \psi_{kj}, b_l\rangle.$ According to Proposition 1 in \cite{hu2022}, we have
$$
\sum_{j'=1}^p \sum_{l' = 1}^\infty  \cov(\theta_{ijl},\theta_{ij'l'}) \eta_{kjl} = \tilde \lambda_k \eta_{kjl}, ~~i \in [n],~j \in [p],~k,l = 1,2, \dots. 
$$

Under the practical scenario (\ref{model}) with $T_{ij} \rightarrow \infty$ (dense case), \cite{hu2022} proposed to estimate $\theta_{ijl}$ by $\widehat \theta_{ijl} = T_{ij}^{-1}\sum_{t=1}^{T_{ij}} \{Y_{ijt} - \widehat \mu_j(U_{ijt})\}b_l(U_{ijt})$ and computed the sample variances of $\widehat \theta_{ijl},$ base on which performing the thresholding selection to encourage the sparsity.
However, this approach is not applicable to sparsely observed functional data with bounded $T_{ij},$ as $\theta_{ijl}$ can not be accurately estimated. We consider bridging the gap under the sparse case by applying the local linear smoothers to estimate the covariance functions $\Sigma_{jj'}(u,v)$ directly. 
Note that $\cov(\theta_{ijl},\theta_{ij'l'}) = \int  b_l(u) \Sigma_{jj'}(u,v) b_{l'}(v) {\rm d}u {\rm d}v.$ Instead of estimating $\theta_{ijl}$'s, we propose to directly estimate $\cov(\theta_{ijl},\theta_{ij'l'})$ by $\widehat{\cov}(\theta_{ijl},\theta_{ij'l'}) = \int_{\cU}  b_l(u) \widehat \Sigma_{jj'}(u,v) b_{l'}(v) {\rm d}u{\rm d}v.$ In the special case with $j=j'$ and $l=l',$ this estimate degenerates to the estimated variance of $\theta_{ijl},$ and hence the thresholding idea proposed in \cite{hu2022} can still be employed. To establish the convergence properties under high-dimensional settings, it is essential to develop concentration results for $\{\widehat{\cov}(\theta_{jl},\theta_{j'l'})\}.$  By the fact that
\begin{equation}
\label{cov.theta}
|\widehat{\cov}(\theta_{ijl},\theta_{ij'l'}) - \cov(\theta_{ijl},\theta_{ij'l'})| \leq \|\widehat \Sigma_{jj'} - \Sigma_{jj'}\|_{\cS},
\end{equation}
our derived rate (\ref{cov_max_L2_dense}) in Theorem~\ref{thm_cov_maxrate} becomes applicable. Specially, when $j=j'$ and $l=l',$ (\ref{cov.theta}), we can obtain the the elementwise maximum rate for $\widehat{\var}(\theta_{ijl}).$ It is noteworthy that \cite{hu2022} relied on existing concentration inequalities for $\chi_n^2$ to establish the concentration bound on the sample variance of $\widehat\theta_{ijt}$ under the dense design by assuming that $\theta_{ijl}$ and  $\varepsilon_{ijt}$ are jointly Gaussian. By comparison, our proposal can handle both sparsely and densely observed functional data without requiring the Gaussianity assumption.

It is also worth mentioning that the above discussion also applies to other pre-fixed basis expansion methods when fitting functional additive regression models \cite[]{fan2014,fan2015,xue2021} or implementing functional linear discriminant analysis \cite[]{xue2023}. Specifically, our proposed $\{\widehat{\cov}(\theta_{ijl},\theta_{ij'l'})\}$ are also involved in the estimation, making Theorem~7 useful in this context.

%in extremely sparse scenarios, where only a handful of measurements are available for each trajectory, as $\theta_{ijl}$ cannot be accurately estimated. We propose to bridge this gap by estimating the covariances directly. Note that $cov(\theta_{jl},\theta_{j'l'}) = \int  b_l(u) \Sigma_{jj'}(u,v) b_{l'}(v) dudv.$  Instead of estimating $\theta_{ijl}$'s, we suggest directly estimating that $cov(\theta_{jl},\theta_{j'l'})$ by $\widehat{cov}(\theta_{jl},\theta_{j'l'}) = \int  b_l(u) \widehat \Sigma_{jj'}(u,v) b_{l'}(v) dudv.$ Consequently, leveraging the estimated covariances $\widehat{cov}(\theta_{jl},\theta_{j'l'}) $'s, we can recover the sparse FPCA through the methodology outlined in \cite{hu2022}. Furthermore, since 
%$$
%|\widehat{cov}(\theta_{jl},\theta_{j'l'}) - cov(\theta_{jl},\theta_{j'l'})| \le \| \widehat \Sigma_{jj'} - \Sigma_{jj'}\|_{\cS},
%$$
%which can be derived directly from the definition of the Hilbert-Schmidt norm $\|\cdot\|_{\cS}$, the asymptotic properties associated with this approach can be theoretically guaranteed.

\subsection{Functional Thresholding}
\label{sec_ft}
Our third application involves estimating the covariance matrix function $\bSigma(\cdot,\cdot)$ in (\ref{cov.fmat}). Under the functional sparsity assumption that $\bSigma$ belongs to a class of ``approximately sparse'' covariance matrix functions, \cite{fang2022} proposed to perform adaptive functional thresholding on the local linear smoothers $\{\widehat\bSigma_{jk}(\cdot,\cdot)\}_{j,k \in [p]}$ using entry-dependent functional thresholds that automatically adapt to the variability of $\widehat \bSigma_{jk}(\cdot,\cdot)$'s. Specifically, their adaptive functional thresholding estimator  is defined as
$$\widetilde \bSigma=\big\{\widetilde \Sigma_{jk}(\cdot,\cdot)\big\}_{p\times p}~~\text{with}~~\widetilde \Sigma_{ij} = \Psi_{jk}^{1/2} \times s_{\lambda}\Big(\frac{\widehat \Sigma_{jk}}{\Psi_{jk}}\Big),$$
where, for any thresholding parameter $\lambda \geq 0$, $s_{\lambda}(\cdot)$ is the functional thresholding operator to enforce the functional sparsity with the aid of Hilbert--Schmidt norm, and $\Psi_{jk}(\cdot,\cdot)$ is a surrogate estimator for the asymptotic variance of $\widehat\Sigma_{jk}(\cdot,\cdot)$ \cite[]{fang2022}. Alternatively, one 
can achieve a universal functional thresholding estimator 
$$\widecheck \bSigma=\big\{\widecheck \Sigma_{jk}(\cdot,\cdot)\big\}_{p\times p}~~\text{with}~~\widecheck \Sigma_{jk}=s_{\lambda}(\widehat \Sigma_{jk}),$$ where a universal threshold level is used for all entries. To investigate the theoretical properties of both functional thresholding estimators $\widetilde \bSigma$ and $\widecheck \bSigma$ under high-dimension settings, it is crucial to make use of the elementwise maximum rate of $\max_{j,k \in [p]}\|\widehat\Sigma_{jk}-\Sigma_{jk}\|_{\cS},$ as presented in Theorem~\ref{thm_cov_maxrate}. In contrast, \cite{fang2022} assumed a rough rate (\ref{max.fcm}) instead of providing a proof to facilitate their technical analysis.

Such functional sparsity assumption is restrictive for many data sets, especially in finance and economics, where the functional variables exhibit high correlations. To address this issue, \cite{li2023} and \cite{leng2024} employed the functional factor models framework for $\bX_i(\cdot),$ which is decomposed as the sum of a common component driven by low-dimensional latent factors and an idiosyncratic component $\be_i(\cdot)$. Instead of imposing the functional sparsity assumption on $\bSigma,$
they imposed it on the covariance matrix function of $\be_i(\cdot),$ and proposed different estimators for $\bSigma$ by performing the associated eigenanalysis to estimate the common covariance and then applying (adaptive) functional thresholding to the residual covariance. To develop the convergence rates of their proposed estimators, our elementwise maximum rate in Theorem~\ref{thm_cov_maxrate} becomes applicable.

\section{Simulations}\label{sec.sim}
In this section, we examine the finite-sample performance of the local linear smoothers for the mean and covariance function estimation in high dimensions.

We generalize the simulated example for univariate functional data in \cite{zhang2016} to the multivariate setting by
generating $$X_{ij}(u) = \mu_j(u) + \bphi(u)^{\T} \btheta_{ij},~~~i \in [n], j \in [p], u \in \cU=[0,1],$$ %for $i \in [n], j \in [p]$ and $u \in \cU=[0,1],$
where the true mean function
$\mu_j(u)=1.5 \sin\{3 \pi(u+0.5)\} +2 u^3,$ the basis function
$\bphi(u)=\big\{\sqrt{2} \cos(2\pi u), \sqrt{2} \sin(2\pi u),\sqrt{2} \cos(4\pi u),\sqrt{2} \sin(4\pi u)\big\}^{\T}$ and the basis coefficient vector $\btheta_i = \big(\btheta_{i1}^{\T}, \dots, \btheta_{ip}^{\T}\big)^{\T} \in {\mathbb R}^{4p}$ is sampled independently from a mean zero multivariate Gaussian distribution with block covariance matrix $\bLambda \in {\mathbb R}^{4p \times 4p}$ whose $(j,k)$th block is given by $\bLambda_{jk}=\rho^{|j-k|}\diag\{2^{-2}, \dots, 5^{-2}\} \in {\mathbb R}^{4 \times 4}$ for $j,k \in [p].$ Hence the $(j,k)$th entry of the true covariance functions $\bSigma(\cdot,\cdot) = \{\Sigma_{jk}(\cdot,\cdot)\}_{p \times p}$ is $\Sigma_{jk}(u,v)=\bphi(u)^{\T} \bLambda_{jk} \bphi(v).$ We then generate the observed values $Y_{ijt} = X_{ij}(U_{ijt}) + \varepsilon_{ijt}$ for $t=1, \dots, T_{ij}=T,$ where the time points $U_{ijt}$'s and errors $\varepsilon_{ijt}$'s are sampled independently from $\text{Uniform}[0,1]$ and ${\cal N}(0,0.5^2),$ respectively.

We use the Epanechnikov kernel with bandwidth values varying on a dense grid.
To evaluate the performance of $\widehat \mu_j(\cdot)$ for $j \in [p]$ and $\widehat \Sigma_{jk}(\cdot,\cdot)$ for $(j,k) \in [p]^2$ given specific bandwidth $h_{\mu,j}$ and $h_{\sSigma,jk},$ we define the corresponding mean integrated squared errors (MISE) as
$$\text{MISE}(\widehat\mu_j,h_{\mu,j})=\int \{\widehat\mu_j(u)-\mu_j(u)\}^2 \mathrm d u,~~\text{MISE}(\widehat\Sigma_{jk},h_{\sSigma,jk})=\int\int \{\widehat\Sigma_{jk}(u,v)-\Sigma_{jk}(u,v)\}^2 \mathrm d u \mathrm d v.$$
We first calculate the elementwise minimal MISEs for the mean and covariance estimators over the grids of candidate bandwidths in prespecified sets $\cal H_{\mu}$ and ${\cal H}_{\sSigma},$
%${\cal H}_{\mu}$ and ${\cal H}_{\sSigma},$
respectively. We then compute their averages and maximums over $j \in [p]$ for the mean functions, that is,
\begin{flalign*}
\text{AveMISE}(\mu)&=\frac{1}{p}\sum_{j}\min_{h_{\mu,j}}\text{MISE}(\widehat\mu_j,h_{\mu,j}),\\
~~~\text{MaxMISE}(\mu)&=\max_j\min_{h_{\mu,j}}\text{MISE}(\widehat\mu_j,h_{\mu,j}),
\end{flalign*}
and over $(j,k) \in [p]^2$ for the covariance functions, that is,
\begin{flalign*}
\text{AveMISE}(\Sigma)&=\frac{1}{p^2}\sum_{j} \sum_k \min_{h_{\sSigma,jk}}\text{MISE}(\widehat\Sigma_{jk},h_{\sSigma,jk}),\\
\text{MaxMISE}(\Sigma)&=\max_{j,k}\min_{h_{\sSigma,jk}}\text{MISE}(\widehat\Sigma_{jk},h_{\sSigma,jk}).
\end{flalign*}
We next use the example of estimating the mean functions to illustrate the rationale of the above measures. While AveMISE($\mu$) presents the averaged elementwise minimal MISEs across $j \in [p],$ some simple calculations in Appendix~\ref{supp_claim} show that the attainable quantity of the minimal elementwise maximum of MISEs, $\min_{(h_{\mu,1}, \dots, h_{\mu,p}) \in {\cal H}_{\mu}^p}\max_{j \in [p]} \text{MISE}(\widehat \mu_j, h_{\mu,j}),$ is equal to MaxMISE($\mu$).
It is worth noting that the optimal selection of bandwidths by minimizing, e.g., for mean functions, each $\text{MISE}(\widehat\mu_j,h_{\mu,j})$ over $h_{\mu,j} \in \cal H_{\mu}$ or $\max_j\text{MISE}(\widehat\mu_j,h_{\mu,j})$ over $(h_{\mu,1}, \dots, h_{\mu,p}) \in {\cal H}_{\mu}^p$ serves to validate our established theoretical results, assuming that $\mu_j(\cdot)$'s are known. In practical scenarios with unknown mean and covariance functions, one can employ the commonly-adopted cross-validation method to estimate the specific MISE as a function of bandwidths, whose minimizer produces the optimal bandwidth selection.

%Such argument applies analogously to estimation of covariance functions.
%It is also noteworthy that the averaged minimal MISE, AveMISE($\mu$)

%o ensure the computational feasibility, we specify a set $\cal H$ of candidate bandwidth values and use $\min_{(h_{\mu,1}, \dots, h_{\mu,p}) \in {\cal H}_{\mu}^p \max_{j} \text{MISE}(\widehat \mu_j, h_{\mu,j})$

%To ensure the feasibility of calculation, we specify a set $\cal H$ of $r$ candidate values for bandwidth, and use $\min_{(i_1,...,i_p) \in [r]^p} \max_{j \in [p]} \text{MISE}(\widehat\mu_j,h_{\mu,i_j}) $ to approximate the $\ell_{\infty}$
%maximum of minimal MISE for the estimated mean functions
%max $l_2$-norm error of mean functions estimates. By some simple calculations, we can find that $\max_j\min_{h_{\mu,j}}\text{MISE}(\widehat\mu_j,h_{\mu,j})$ is exactly what we want. The same applies to covariance case.

\begin{figure}[!htbp]
\centering
\includegraphics[width=7cm,height=6cm]{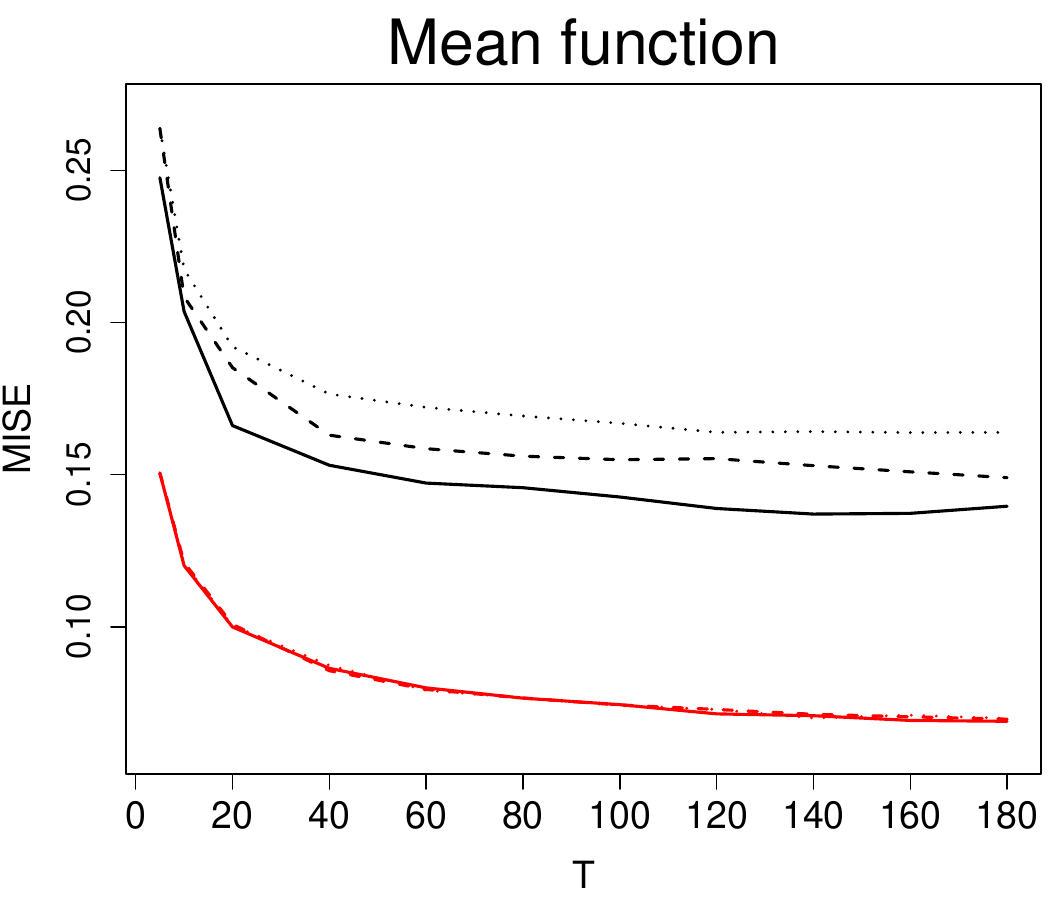}
\hspace{.5cm}
\includegraphics[width=7cm,height=6cm]{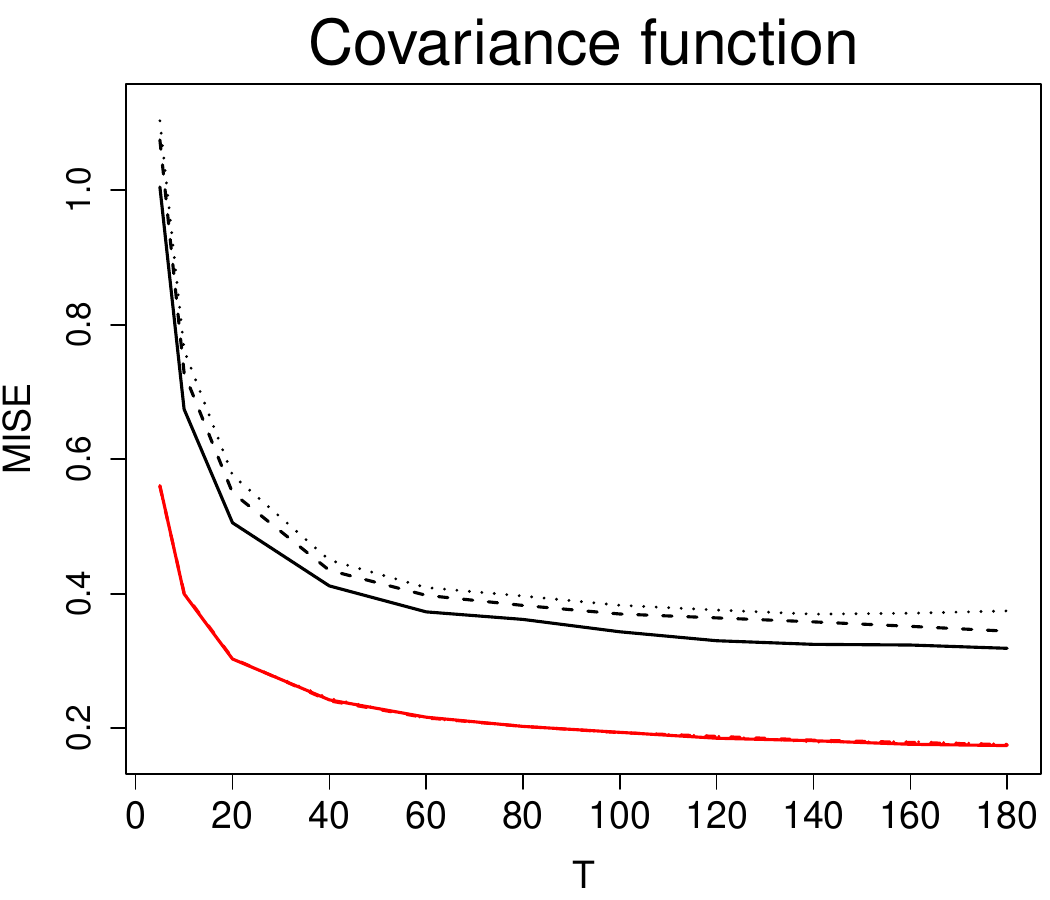}
\caption{\label{plot.sim1} Plots of average MaxMISE (black) and AveMISE (red) against $T$ with $p=50$ (solid), $100$ (dashed) and $150$ (dotted) for mean estimators (left) and covariance estimators (right).
}
\end{figure}

%\begin{figure}[!htbp]
%\centering
%\includegraphics[width=7cm,height=6cm]{Rplotmean.pdf}
%\hspace{.2cm}
%\includegraphics[width=7cm,height=6cm]{Rplotcov.pdf}
%\caption{\label{plot.sim2} \textcolor{red}{Plots of the average $\log(\text{MaxMISE})$s (black) and $\log(\text{AveMISE})$s (red) against $\log(T)$ with $p=50$ (solid), $100$ (dashed) and $150$ (dotted) for the mean estimators (left) and the covariance estimators (right).}
%}
%\end{figure}

%\begin{figure}[!htbp]
%\centering
%\includegraphics[width=7cm,height=6cm]{mean-logmise-t.pdf}
%\hspace{.2cm}
%\includegraphics[width=7cm,height=6cm]{cov-logmise-t.pdf}
%\caption{\label{plot.sim3} \textcolor{red}{Plots of the average $\log(\text{MaxMISE})$s (black) and $\log(\text{AveMISE})$s (red) against $T$ with $p=50$ (solid), $100$ (dashed) and $150$ (dotted) for the mean estimators (left) and the covariance estimators (right).}
%}
%\end{figure}

We firstly consider settings of $n=100,$ $p=50, 100, 150$, and $T=5$, 10, 20, 40, 60, 80, 100, 120, 140, 160, 180.  %.varying from sparse to semi-dense to ultra-dense measurement schedules. 
We ran each simulation 100 times. Figure~\ref{plot.sim1} plots the average AveMISE and MaxMISE as functions of %the sampling frequency 
$T$ for the estimated mean and covariance functions. %Two apparent patterns are observable from Figure~\ref{plot.sim1} 
We observe that both MaxMISE and AveMISE display a similar trend as $T$ increases %from $5$ to $160,$ 
with a steep decline followed by a slight decrease and then a period of stability. Such trend roughly corresponds to the three categories of ``sparse'', ``semi-dense'', and ``ultra-dense'', respectively.
In addition, while AveMISE reflects the performance for univariate functional data, MaxMISE gradually enlarges as $p$ increases from $50$ to $150,$ providing empirical evidence to support that the associated $\log p$-based convergence rates in high-dimensional settings. %. all depend on $\log p$-based multiplicative factors. %Additionally, it is observable that the high-dimensionality causes the transition phase between semi-dense and ultra-dense functional data to slightly shift to the left. 

%Third, compared to the results for the estimated mean functions, an increase in $T$ for sparse and semi-dense functional data leads to an enhanced reduction in relative MISEs for the covariance estimators. 

\begin{figure}[t]
\centering
\includegraphics[width=7.5cm,height=6cm]{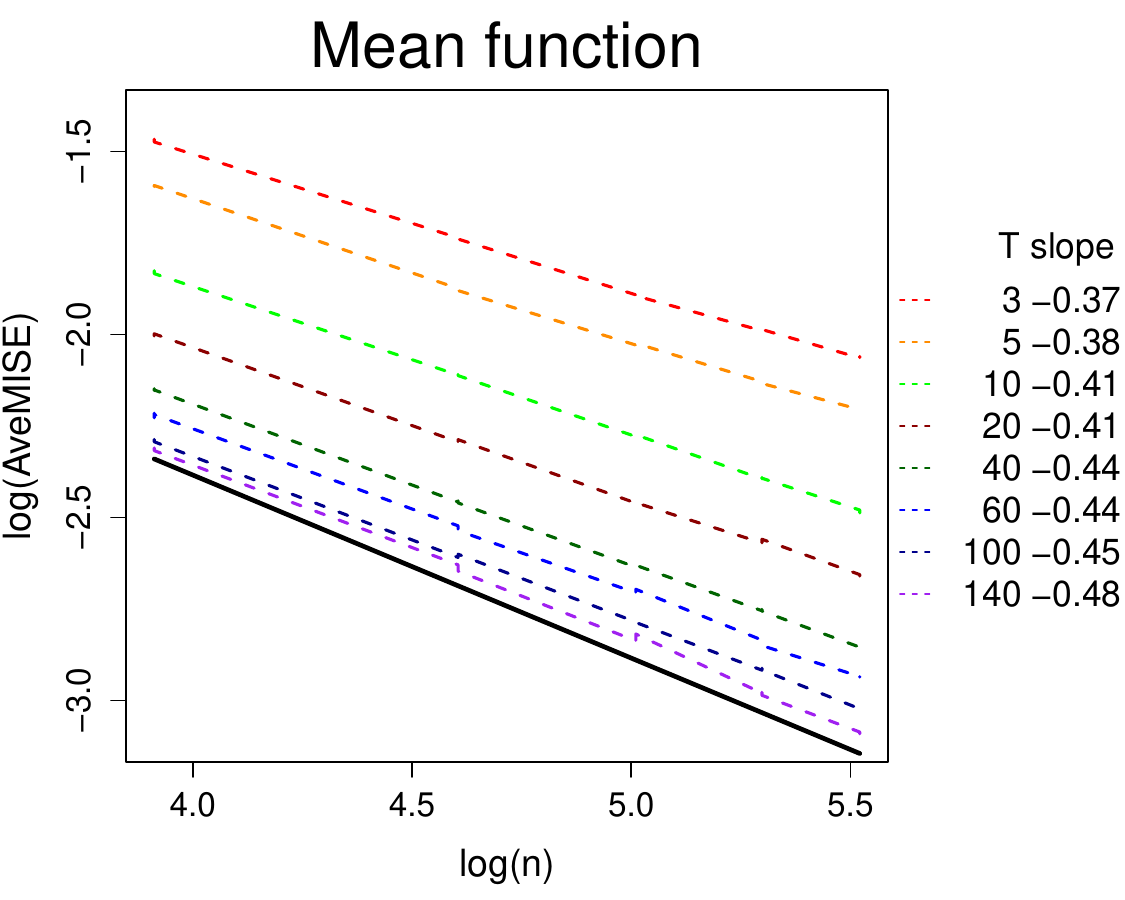}
\hspace{0cm}
\includegraphics[width=7.5cm,height=6cm]{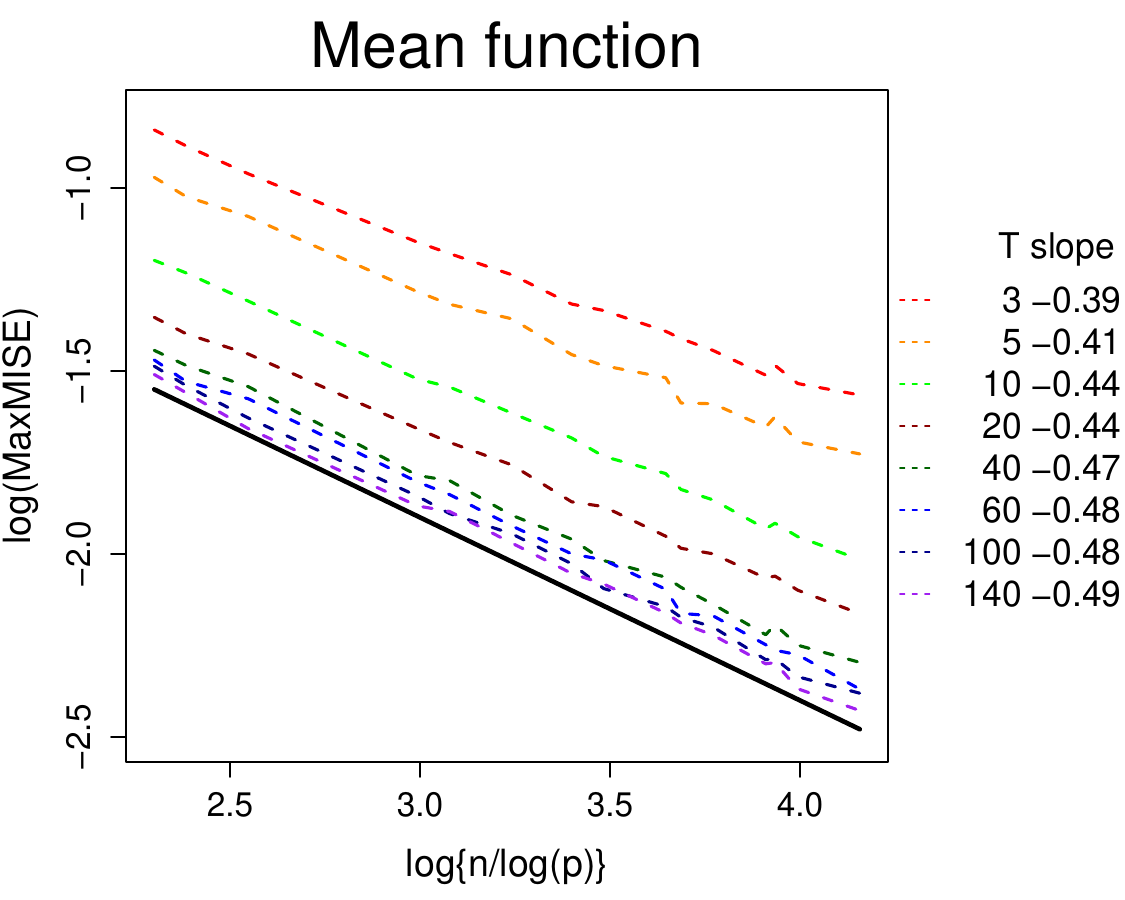}
\includegraphics[width=7.5cm,height=6cm]{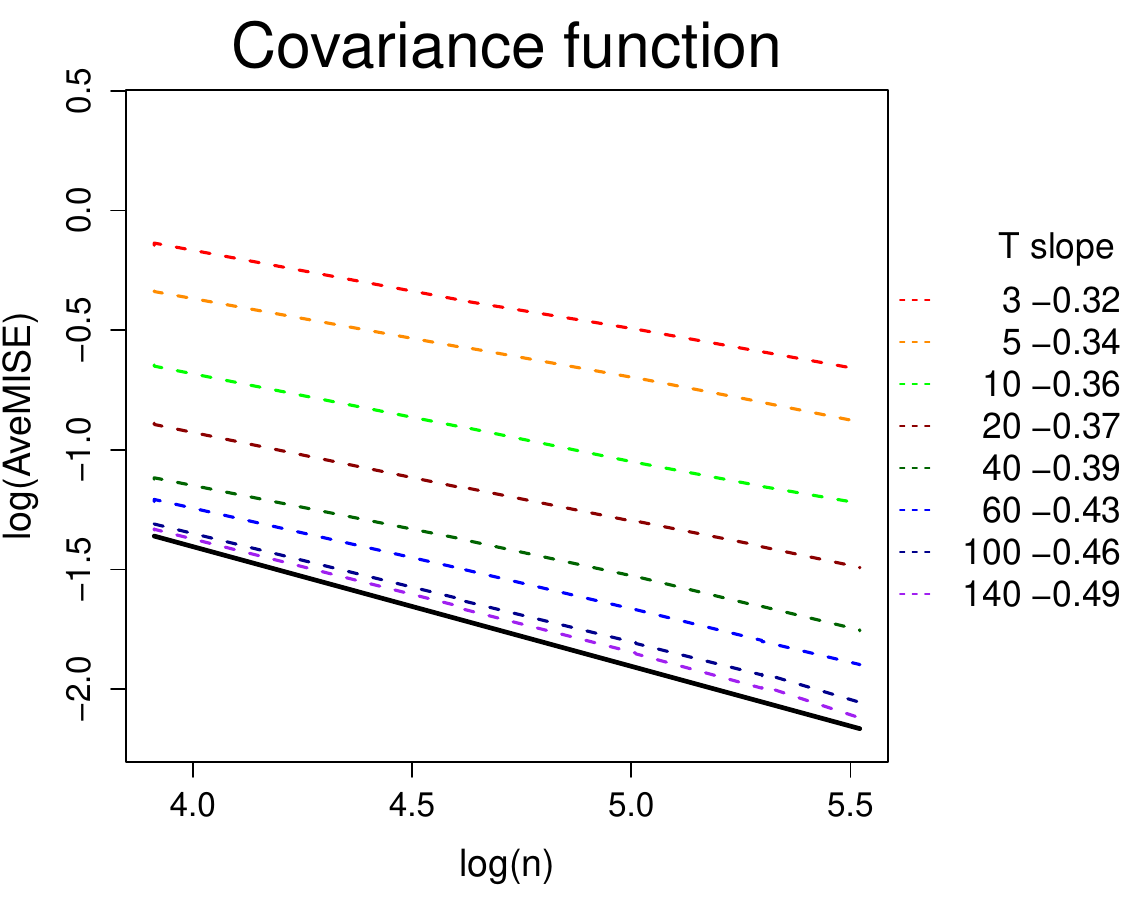}
\hspace{0cm}
\includegraphics[width=7.5cm,height=6cm]{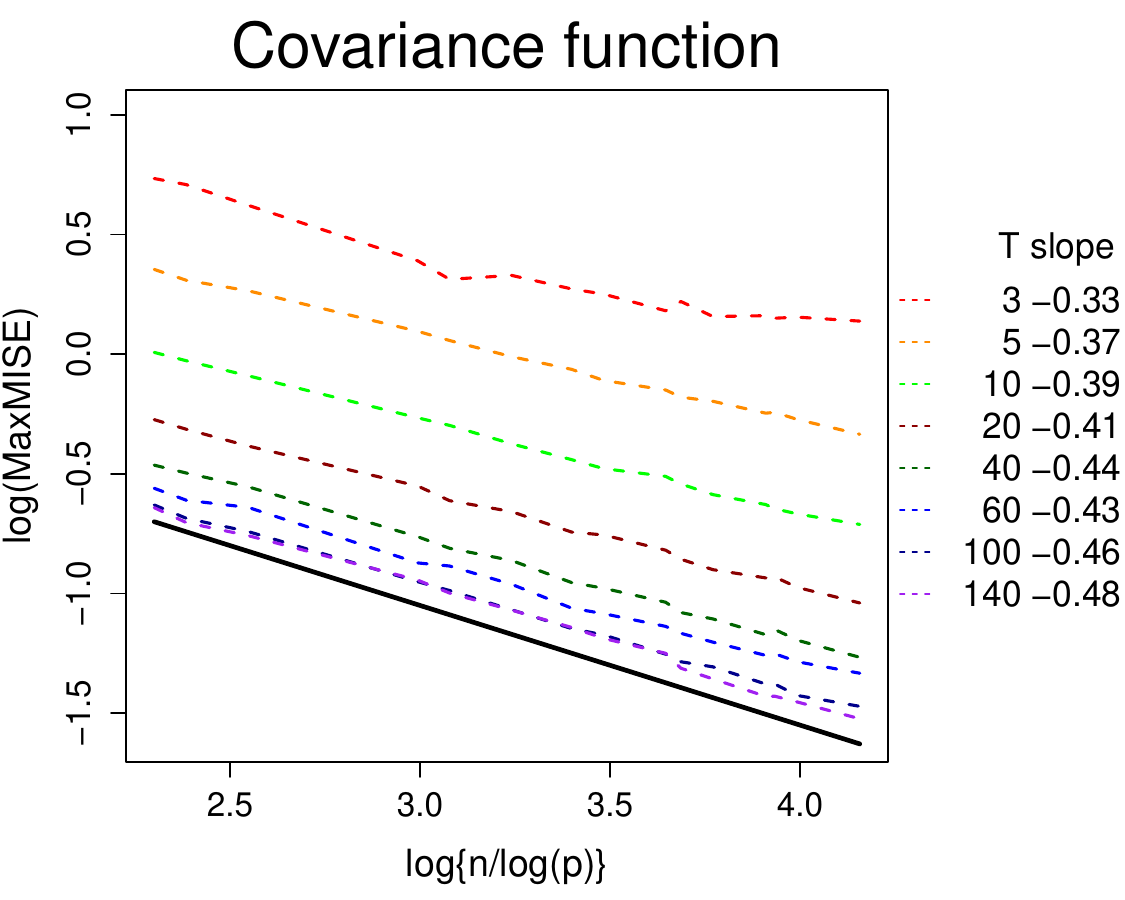}
\caption{\label{plot.sim2} Plots of average $\log(\text{AveMISE})$ against $\log n$ (left) and average $\log(\text{MaxMISE})$ against $\log(n/\log p )$ (right) for mean estimators (top) and covariance estimators (bottom) with $p=50,100,150$ and $n=50,100,150,200,250.$ The colored dashed lines correspond to different values of $T$ ranging from 3 to 140, and the estimated slopes of the corresponding linear fits based on five points for $\log n$ or fifteen points for $\log(n/\log p)$ are also displayed. The slope of the black solid line presents the theoretical value -1/2 (with the intercept being irrelevant here).
}
\end{figure}

%begin{figure}[!htbp]
%\centering
%\includegraphics[width=7.5cm,height=6cm]{cov-ave.pdf}
%\hspace{0cm}
%\includegraphics[width=7.5cm,height=6cm]{cov-max.pdf}
%\caption{\label{plot.sim3} \textcolor{red}{Plots of average $\log(\text{AveMISE})$s against $\log n$ (left) and average $\log(\text{MaxMISE})$s against $\log\{n/\log(p)\}$ (right) for covariance estimators
%with $p=50,100,150$ and $n=50,100,150,200,250$ across different $T$ values ranging from $3$ to $140.$ The estimated slopes in the corresponding linear fits are also displayed.}
%}
%\end{figure}

As suggested by one referee, to further validate the established rates discussed in Remarks~\ref{rmk_asy_phase} and \ref{rmk_nonasy_phase}, we plot the average log(AveMISE) against $\log n$ and the average log(MaxMISE) against $\log(n/\log p)$ across different values of $n,$ $p$ and $T$ for estimated mean and covariance functions in
Figure~\ref{plot.sim2}.
Note that our theory indicates that log(AveMISE) and $\log n,$ as well as log(MaxMISE) and $\log(n/\log p),$ both exhibit linear relationships, e.g., for sparse functional data, with a slope of -2/5 for mean estimators and -1/3 for covariance estimators. Furthermore, when the phase transitions from ``semi-dense" to ``ultra-dense" occur, they all admit linear relationships with a common slope of -1/2. 
Several apparent patterns are observable from Figure~\ref{plot.sim2}.  
First, both plots of log(AveMISE) against five values of $\log n$ and log(MaxMISE) against fifteen values $\log(n/\log p)$ show clear linear patterns across different values of $T.$ This also demonstrates the $\log p$-based high-dimensional effect in elementwise maximum rates.
Second, for sparse functional data (i.e., $T=3, 5$), the estimated slopes of the linear fits are between -0.41 and -0.37 for mean estimators and between -0.37 and -0.32 for covariance estimators, which nicely align with the theoretical slope values of -2/5 and -1/3, respectively.
Third, as $T$ grows from 10 to 140, the estimated slopes under all scenarios gradually increase to some values close to the theoretical slope value of -1/2, and then stabilize, especially when measuring MaxMISE. This suggests that the phase transitions from ``semi-dense" to ``ultra-dense" may occur during the increase of $T.$ For example, the slope for MaxMISE tends to be stable for $T=40, 60, 100, 140,$ indicating that the phase transition may occur around $T=40$ (earlier than the corresponding occurrence of phase transition for AveMISE).  
All of the above observations nicely validate our established theoretical results in Section~\ref{sec.theory}.

\subsection*{Acknowledgments} 
The authors would like to thank the Action Editor and three anonymous reviewers for their insightful comments and suggestions, which have led to significant improvement of our paper. 
Shaojun Guo was partially supported by the National Natural Science Foundation of China (No. 12471270). Dong Li was partially supported by the National Natural Science Foundation of China (No. 72471127). Xinghao Qiao was partially supported by the Seed Fund for Basic Research for New Staff at the University of Hong Kong.

\appendix

\section{Technical Proofs}\label{ap.pf}
The appendix contains proofs of all theorems. Throughout, we use $c, c_1, c_2, \dots$ to denote generic positive finite constants that may be different in different uses. %whose values can be different in different uses but will not cause any confusion.

\setcounter{equation}{0}
\renewcommand{\theequation}{A.\arabic{equation}}

\subsection{Proof of Theorem~\ref{thm_mean_coneq}}
\label{pf.thm.mean}

We organize the proof in four steps. First, we will define $ \widehat \mu(\cdot),$ $ \widetilde\mu(\cdot)$ and obtain the decomposition of $ \widehat\mu(\cdot) - \widetilde\mu(\cdot).$ Second, we will prove the local concentration inequality for fixed interior point $u \in \cU.$ Third, we will prove the concentration inequality in $L_2$ norm. Finally, we will prove the concentration inequality in supremum norm.

\subsubsection{Definition and Decomposition}\label{suppsec:def}
Without loss of generality, we let $ h_{\mu,j} = h$ for $j \in [p]$ and denote $\be_0 = (1, 0)^{\T},$ $\widetilde{\bU}_{ijt} = \big\{1,(U_{ijt}-u)/h\big\}^{\T},$
\begin{equation}
	\begin{split}
		%& \widetilde{U}_{ijt} = \Big(1,\frac{U_{ijt}-u}{h}\Big)^{\T}, \\
		& \widehat{\bS}_j(u) = \sum_{i=1}^{n} v_{ij} \sum_{t=1}^{T_{ij}} \widetilde{\bU}_{ijt} \widetilde{\bU}_{ijt}^{\T} K_h(U_{ijt}-u), \\
		&\widehat{\bR}_j(u) = \sum_{i=1}^{n} v_{ij} \sum_{t=1}^{T_{ij}} \widetilde{\bU}_{ijt} K_h(U_{ijt}-u) Y_{ijk}. \\
	\end{split}
    \nonumber
\end{equation}
%A simple calculation yields $\widehat{\mu}_j(u) = e_0^{\T} \{\widehat{S}_j(u)\}^{-1} \widehat{R}_j(u).$
%\begin{equation}
%	\widehat{\mu}_j(u) = e_0^{\T} \{\widehat{S}_j(u)\}^{-1} \widehat{R}_j(u).
%\end{equation}
%$$
%\widehat{\mu}_j(u) = e_0^{\T} \{\widehat{S}_j(u)\}^{-1} \left\{\widehat{R}_{j1}(u) +  \widehat{R}_{j2}%(u)\right\}
%$$
%where
%$$
%\widehat{R}_{j1}(u) = \sum_{i=1}^{n} v_{ij} \sum_{t=1}^{T_{ij}} \widetilde{U}_{ijt} K_h(U_{ijt}-u) %\mu_j(U_{ijt})
%$$
%,
%$$
%\widehat{R}_{j2}(u) = \sum_{i=1}^{n} v_{ij} \sum_{t=1}^{T_{ij}} \widetilde{U}_{ijt} K_h(U_{ijt}-u) \xi_{ijt}
%$$
%and $ \xi_{ijt} = Y_{ijt} - \mu_j(U_{ijt}) $.
A simple calculation yields that $\widehat{\mu}_j(u) = \be_0^{\T} \{\widehat{\bS}_j(u)\}^{-1} \widehat{\bR}_j(u).$
Let
\begin{equation}
\label{tilde_mu}
\widetilde{\mu}_j(u) = \be_0^{\T} \big[\cE\{\widehat{\bS}_j(u)\}\big]^{-1} \cE\{\widehat{\bR}_j(u)\}.
\end{equation}
%\begin{equation}
%	\widetilde{\mu}_j(u) = e_0^{\T} [E\{\widehat{S}_j(u)\}]^{-1} E\{\widehat{R}_j(u)\},
%	\nonumber
%\end{equation}
We can decompose $\widehat{\mu}_j(u) - \widetilde{\mu}_j(u)$ as
\begin{equation}
	\label{decomp_mean}
	\begin{split}
		\widehat{\mu}_j(u) - \widetilde{\mu}_j(u) = &
   \be_0^\T [\cE\{\widehat{\bS}_j(u)\}]^{-1} [\widehat{\bR}_{j}(u)-\cE\{\widehat{\bR}_j(u)\}]\\ &- \be_0^\T \{\widehat{\bS}_j(u)\}^{-1}[\widehat{\bS}_j(u) - \cE\{\widehat{\bS}_j(u)\}] [\cE\{\widehat{\bS}_j(u)\}]^{-1} \widehat{\bR}_j(u),
  \end{split}
  \nonumber
\end{equation}
which then implies that
\begin{equation}
	\label{ineq_decomp_mean}
	\begin{split}
		|\widehat{\mu}_j(u) - \widetilde{\mu}_j(u)| & \le \|\cE\{\widehat{\bS}_j(u)\}\|_{\min}^{-1} \|\widehat{\bR}_{j}(u)- \cE\{\widehat{\bR}_j(u)\}\| \\
  & \quad + \|\widehat{\bS}_j(u)\|_{\min}^{-1}\|\cE\{\widehat{\bS}_j(u)\}\|_{\min}^{-1}\|\widehat{\bR}_j(u)\|\|\widehat{\bS}_j(u) - \cE\{\widehat{\bS}_j(u)\}\|_{\tF}, \\
\end{split}
\end{equation}
where, for any vector $\bbb=(b_1, \dots, b_p)^{\T},$ we denote its Euclidean norm by $\|\bbb\|=(\sum_{i} b_i^2)^{1/2}$ and, for any matrix $\bB=(B_{ij})_{p \times q}$, we write $\|\bB\|_{\min}=\{\lambda_{\min}(\bB^\T \bB)\}^{1/2}$ and $\|\bB\|_{\tF}=(\sum_{i,j} B_{ij}^2)^{1/2}$ to denote its Frobenius norm.

\subsubsection{Local Concentration Inequality}
%To prove the concentration inequality of $|\widehat{\mu}_j(u) - \widetilde{\mu}_j(u)|.$
%In this subsection, we prove that for every $ u \in \cU $,
%\begin{equation}
%	\label{pointwise_mu}
%	\pr(|\widehat{\mu}_j(u) - \widetilde{\mu}_j(u)| \geq \delta) \leq c \exp\left(- \frac{c \gamma_{n,T,h,j} %\delta^2}{1+ \delta} \right).
%\end{equation}
We will firstly show that there exists some positive constant $c$ (independent of $u$) such that for any $\delta >0$ and $ u \in \cU $,
\begin{equation}
\label{S_ju_bound}
	\pr\Big\{\big\|\widehat{\bS}_{j}(u) - \cE\{\widehat{\bS}_{j}(u)\}\big\|_{\tF} \geq \delta \Big\}\leq 8 \exp\Big(- \frac{c n\widebar{T}_{\mu,j}h \delta^2}{1+ \delta} \Big).
\end{equation}

For $ k,l=1,2 $, let $ \widehat{S}_{jkl}(u) $ be the $ (k,l)$th entry of $\widehat{\bS}_{j}(u).$
Under Assumptions~\ref{cond_T_bd} and \ref{cond_kernel}, we obtain that for any integer $q = 2,3,\dots$ and $s = 0, 1, 2,$
\begin{equation}
\label{iequ.kernel}
\cE\left\{\Big|\Big(\frac{U_{ijt}-u}{h}\Big)^s K_h(U_{ijt}-u)\Big|^q\right\} \leq \int h^{-q} K^q\Big(\frac{t-u}{h}\Big) \Big|\frac{t-u}{h} \Big|^{sq} f_U(t) dt \leq c h^{1-q}.
\end{equation}
Note that Assumption~\ref{cond_weight} implies that the weights $v_{ij}$'s are of the same order $v_{ij} \asymp (n\widebar T_{\mu,j})^{-1}.$ By (\ref{iequ.kernel}), it holds that 
$$\sum_{i=1}^{n} \sum_{t=1}^{T_{ij}} \cE\left\{\Big|\Big(\frac{U_{ijt}-u}{h}\Big)^s K_h(U_{ijt}-u)\Big|^2\right\} \leq c n \widebar T_{\mu,j} h^{-1},
$$
$$
\sum_{i=1}^n \sum_{t=1}^{T_{ij}} \cE\left\{\Big|\Big(\frac{U_{ijt}-u}{h}\Big)^s K_h(U_{ijt}-u)\Big|^q\right\} \leq 2^{-1} q! c n \widebar T_{\mu,j} h^{-1} h^{2-q} \quad\mbox{ for $q \geq 3.$}
$$
By the Bernstein inequality \citep[see Theorem~2.10 and Corollary~2.11 of][]{BLM2014}, we obtain that there exists some positive constant $c$ (independent of $u$) such that for any $\delta >0$ and $ u \in \cU,$
$$
\pr\Big\{\big|\widehat{S}_{jkl}(u) - \cE\{\widehat{S}_{jkl}(u)\}\big| \geq \delta \Big\}\leq 2 \exp\Big(- \frac{c n\widebar{T}_{\mu,j}h \delta^2}{1+ \delta} \Big),
$$
for $k,l=1,2,$ which, by the union bound of probability, implies that (\ref{S_ju_bound}) holds.

For $ k=1,2 $, let $ \widehat{R}_{jk}(u) $ be the $ k$th element of $ \widehat{\bR}_{j}(u) $. We will next show that there exists some positive constant $c$ (independent of $u$) such that for any $ \delta > 0 $ and $ u \in \cU,$
\begin{equation}
\label{coneq.R}
	\pr\Big\{\big|\widehat{R}_{jk}(u) - \cE\{\widehat{R}_{jk}(u)\}\big| \geq \delta\Big\} \leq c \exp\Big(- \frac{c \gamma_{n,T,h,j} \delta^2}{1+ \delta}\Big),
\end{equation}
where $\gamma_{n,T,h,j}=n (1 \wedge \widebar T_{\mu,j} h).$
We only need to consider the case $k=1,$ while the case $k=2$ can be demonstrated in a similar manner. Denote that 
\begin{flalign*}
\xi_{ijt} &= Y_{ijt} - \mu_j(U_{ijt}),\\
\widehat{R}_{j3}(u) & = \sum_{i=1}^{n} v_{ij} \sum_{t=1}^{T_{ij}} K_h(U_{ijt}-u) \mu_{j}(U_{ijt}),\\
\widehat{R}_{j4}(u) &= \sum_{i=1}^{n} v_{ij} \sum_{t=1}^{T_{ij}} K_h(U_{ijt}-u) \xi_{ijt}. 
\end{flalign*} 
Then $\widehat{R}_{j1}(u) - \cE\{\widehat{R}_{j1}(u)\}$ can be rewritten as
\begin{equation}
\label{decompose_R1}
\widehat{R}_{j1}(u) - \cE\{\widehat{R}_{j1}(u)\} = \widehat{R}_{j3}(u) - \cE\{\widehat{R}_{j3}(u)\} + \widehat{R}_{j4}(u).
\end{equation}

Following the same procedure to prove (\ref{S_ju_bound})
%Note that Assumption~\ref{cond_weight} implies that the weights $v_{ij}$'s are of the same order $v_{ij} \asymp (n\widebar T_{\mu,j})^{-1}.$ By (\ref{iequ.kernel}), it holds that $\sum_{i=1}^{n} \sum_{t=1}^{T_{ij}} \cE\{|K_h(U_{ijt}-u) \mu_{j}(U_{ijt})|^2\} \leq c n \widebar T_{\mu,j} h^{-1},$ and
%$\sum_{i=1}^n \sum_{t=1}^{T_{ij}} \cE\{|K_h(U_{ijt}-u) \mu_{j}(U_{ijt})|^q\} \leq 2^{-1} q! c^q n \widebar T_{\mu,j} h^{-1} h^{2-q} $ for $q \geq 3.$
and using the Bernstein inequality, we can obtain that there exists some positive constant $c$ such that for any $\delta >0$ and $ u \in \cU,$
\begin{equation}
\label{concentration_R_3}
	\pr\Big\{\big |\widehat{R}_{j3}(u) - \cE\{\widehat{R}_{j3}(u)\}\big | \geq \delta \Big\}\leq 2 \exp\Big(- \frac{c n\widebar{T}_{\mu,j}h \delta^2}{1+ \delta} \Big).
\end{equation}

Now we consider the tail behavior of $\widehat{R}_{j4}(u).$ Define the event $V_j = \big\{U_{ijt}, t \in [T_{ij}], i \in [n]\big\}.$ %It follows from sub-Gaussianities in Assumption~\ref{cond_subG} and (\ref{model}) that
%\begin{equation}
%\begin{split}
%\cE\left\{\exp(\lambda \xi_{ijt}) \big|V_j\right\} & = \cE\left(\exp\Big[\lambda \big\{X_{ij}(U_{ijt})-\mu_j(U_{ijt})\big\}\Big] \big|V_j\right)+\cE\left\{\exp(\lambda \varepsilon_{ijt}) \big|V_j\right\}\\
%& \le \exp\left[2^{-1} c^2 \lambda^2 \Big\{\sigma_j^2 + \Sigma_{jj}(U_{ijt})\Big\}\right] \le \exp(\lambda^2 c^2). \\
%\end{split} \nonumber
%\end{equation}
%there exists some constant $c>0$ such that, for any $\lambda \in {\mathbb R},$
%$
%\cE\{\exp(\lambda \xi_{ijt}) \big|V_j\big\}\le \exp(\lambda^2 c^2)
%$ for any $\lambda \in {\mathbb R}.$
Rewrite $ \widehat{R}_{j4}(u) = \sum_{i=1}^n v_{ij} \psi_{ij1}(u) $ with $ \psi_{ij1}(u) = \sum_{t=1}^{T_{ij}}  K_h(U_{ijt}-u) \xi_{ijt}.$ %Note that for each $i$ Assumption~\ref{cond_weight} implies that $v_{ij} \asymp (n\widebar T_{\mu,j})^{-1}.$ 
If $\sum_{t=1}^{T_{ij}} K_h(U_{ijt} - u) > 0,$ by Jensen's inequality, we have
\begin{equation}
\begin{split}
& \cE\Big[\exp\{\lambda v_{ij} \psi_{ij1}(u)\} \Big|V_j \Big] \\
= & \cE\Big(\exp\Big[\lambda v_{ij} \frac{\{\sum_{t=1}^{T_{ij}} K_h(U_{ijt} - u)\xi_{ijt}\}\{\sum_{t=1}^{T_{ij}} K_h(U_{ijt} - u)\}}{\sum_{t=1}^{T_{ij}} K_h(U_{ijt} - u)}\Big]\Big|V_j\Big)\\
\leq & \frac{1}{\sum_{t=1}^{T_{ij}} K_h(U_{ijt} - u)} \sum_{t=1}^{T_{ij}} K_h(U_{ijt} - u) \cE\Big[\exp\Big\{\lambda v_{ij} \xi_{ijt}\sum_{t'=1}^{T_{ij}} K_h(U_{ijt'} - u)  \Big\} \Big|V_j\Big] \\
\le & \exp\Big[\lambda^2 c^2 (n\widebar T_{\mu,j})^{-2} \Big\{\sum_{t=1}^{T_{ij}} K_h(U_{ijt} - u)\Big\}^2 \Big].
\end{split}
\nonumber
\end{equation}
Note that the last line comes from the fact that $
\cE\{\exp(\lambda \xi_{ijt}) \big|V_j\big\}\le \exp(\lambda^2 c^2)
$
for any $\lambda \in {\mathbb R},$ which is implied from the sub-Gaussianities in Assumption~\ref{cond_subG} and (\ref{model}).
Clearly,  the above inequality still holds even if $\sum_{t=1}^{T_{ij}} K_h(U_{ijt} - u) = 0.$ Assumption~\ref{cond_timepoints} implies that the number of nonzero terms in $\sum_{t=1}^{T_{ij}} K_h(U_{ijt} - u)$ has an upper bound $c (1\vee \widebar T_{\mu,j}h),$ which yields that $\sum_{t=1}^{T_{ij}} K_h(U_{ijt} - u) \le c h^{-1} (1\vee \widebar T_{\mu,j}h).$ Therefore, for any $\lambda \in \mathbb{R},$ %\textcolor{red}{by Jensen's inequality,}
we obtain that
\begin{equation}
\begin{split}
\cE\left[\exp\Big\{\lambda \sum_{i=1}^nv_{ij} \psi_{ij1}(u) \Big\} \Big|V_j \right]
& = \prod_{i=1}^n \cE\left[\exp\Big\{\lambda v_{ij} \psi_{ij1}(u) \Big\} \Big|V_j \right] \\
&  \le \exp\left[\lambda^2 c^2 \sum_{i=1}^n v_{ij}^2\Big\{\sum_{t=1}^{T_{ij}} K_h(U_{ijt} - u)\Big\}^2 \right]\\
&  \le \exp\Big\{\lambda^2 c^2  (n\widebar T_{\mu,j}h)^{-2} h (1\vee \widebar T_{\mu,j}h)  \sum_{i=1}^n \sum_{t=1}^{T_{ij}} K_h(U_{ijt} - u) \Big\}.
\end{split}
\nonumber
\end{equation}
For any $\delta>0,$ define the event $\Omega_{j,1}(\delta) = \big\{\sum_{i=1}^n \sum_{t=1}^{T_{ij}} K_h(U_{ijt} - u)  \le c (1 + \delta) n\widebar{T}_{\mu,j}\big\}.$  We have
\begin{equation}
\begin{split}
 \cE\left[\exp\Big\{\lambda \sum_{i=1}^nv_{ij} \psi_{ij1}(u) \Big\} \Big|\Omega_{j,1}(\delta) \right] &\le \exp\Big\{\lambda^2 c^2 (1 + \delta) (n\widebar T_{\mu,j}h)^{-2}  n\widebar{T}_{\mu,j}h (1\vee \widebar T_{\mu,j}h) \Big\}\\
 & = \exp\Big\{(1 + \delta) \lambda^2 c^2  (n\widebar T_{\mu,j}h)^{-1}  (1\vee \widebar T_{\mu,j}h) \Big\}.
\end{split}
\nonumber
\end{equation}
As a consequence, we obtain that
\begin{equation}
	\label{sum_exp_psi}
	%\begin{split}
		\pr\Big\{\sum_{i=1}^{n} v_{ij}\psi_{ij1}(u)\geq \delta \Big|\Omega_{j,1}(\delta) \Big\}
		 \leq \exp\Big\{- \lambda \delta + (1 + \delta) \lambda^2 c^2  (n\widebar T_{\mu,j}h)^{-1}  (1\vee \widebar T_{\mu,j}h)\Big\}.
	%\end{split}
\end{equation}
With the choice of $ \lambda = n\widebar{T}_{\mu,j}h \delta/\{2(1 + \delta)c^2 (1 \vee \widebar{T}_{\mu,j}h)\},$ (\ref{sum_exp_psi}) degenerates to
\begin{equation}
	\label{sum_exp_psi2}
		\pr\Big\{\sum_{i=1}^{n} v_{ij} \psi_{ij1}(u) \geq \delta \Big|\Omega_{j,1}(\delta)\Big\}  \leq \exp\Big\{- \frac{cn (1\wedge \widebar{T}_{\mu,j}h)\delta^2}{1 + \delta}\Big\}.
\end{equation}
Note that $\sum_{i=1}^n \sum_{t=1}^{T_{ij}} \cE\{K_h(U_{ijt} - u)\} \le c n \widebar{T}_{\mu,j}.$ By the Bernstein inequality, we obtain that there exists some positive constant $c$ such that for any $\delta >0,$
\begin{equation}
\pr\left( \sum_{i=1}^n \sum_{t=1}^{T_{ij}} \Big[K_h(U_{ijt} - u) - \cE\big\{ K_h(U_{ijt} - u)\big\} \Big] \ge n\widebar{T}_{\mu,j}\delta \right) \le \exp \left( - \frac{c n \widebar{T}_{\mu,j} h \delta^2}{1+\delta}\right),
\nonumber
\end{equation}
which implies that
\begin{equation}
\label{kernel_concentration}
1- \pr\{\Omega_{j,1}(\delta)\} \le \exp \Big(-\frac{c n \widebar{T}_{\mu,j} h \delta^2}{1+\delta}\Big).
\end{equation}
Combining (\ref{sum_exp_psi2}) and (\ref{kernel_concentration}), we obtain that there exists some constant $c>0$ such that for any $\delta >0,$
\begin{equation}
    \pr\Big\{\widehat R_{j4}(u) \geq \delta\Big\} \leq \pr\Big\{\widehat R_{j4}(u) \geq \delta \Big| \Omega_{j,1}(\delta) \Big\} + \pr\Big\{\Omega_{j,1}(\delta)^c\Big\} \leq 2\exp \Big\{- \frac{cn (1\wedge \widebar{T}_{\mu,j}h)\delta^2}{1 + \delta }\Big\},
  \nonumber
\end{equation}
and, consequently,
\begin{equation}
		\pr\Big\{\big|\widehat R_{j4}(u)\big|\geq \delta\Big\} \leq 4\exp \Big\{- \frac{cn (1\wedge \widebar{T}_{\mu,j}h)\delta^2}{1 + \delta }\Big\}.
  \nonumber
\end{equation}
This together with (\ref{decompose_R1}) and (\ref{concentration_R_3}) yields that, there exists some constant $c$ such that
$$
\pr\Big\{\big|\widehat R_{j1}(u) - \cE \{\widehat R_{j1}(u)\}\big|\geq \delta\Big\} \leq 6\exp \Big\{- \frac{cn (1\wedge \widebar{T}_{\mu,j}h)\delta^2}{1 + \delta }\Big\}.
$$
Define the event $\Omega_{j,2}(\delta) = \{\|\widehat{\bS}_j(u) - \cE\{\widehat{\bS}_j(u)\}\|_{\tF} \le \delta/2\}.$
Note that $\cE\{\widehat{\bS}_j(u)\}$  is positive definite.
%It follows from Condition~\ref{cond_kernel} that
%$E\big\{ K_h(U_{ijt}-u)\big\}=1,$
%$E\big\{h^{-1}(U_{ijt}-u) K_h(U_{ijt}-u)\big\}=0$
%and $c<E\big\{h^{-2}(U_{ijt}-u)^2 K_h(U_{ijt}-u)\big\} <\infty,$
%\begin{equation}
%\begin{split}
%& E\big\{ K_h(U_{ijt}-u)\big\} = 1, \\
%& E\Big\{\frac{U_{ijt}-u}{h} K_h(U_{ijt}-u)\Big\} = \int_{\cU} \frac{t-u}{h} K_h(t-u) dt = \int_{-1}^1 x K(x) dx =0, \\
%& E\Big\{(\frac{U_{ijt}-u}{h})^2 K_h(U_{ijt}-u)\Big\} = \int_{\cU} (\frac{t-u}{h})^2 K_h(t-u) dt = \int_{-1}^1 x^2 K(x) dx < \infty, \\
%\end{split} \nonumber
%\end{equation}
%which implies that $E\{\widehat{S}_j(u)\}$  is positive definite.
%means the $(1,1)$th entry of $E\{\widehat{S}_j(u)\}$ is $1$, the $(1,2)$th and $(2,1)$th entry of $E\{\widehat{S}_j(u)\}$ are $0$, the $(2,2)$th entry of $E\{\widehat{S}_j(u)\}$ is positive and finite.From above we prove $E\{\widehat{S}_j(u)\}$ is positive definite with eigenvalues bounded away from $0.$
On the event $\Omega_{j,2}(\delta)$ with $\delta \in (0,1],$ we obtain that
\begin{equation}
\label{Smin.bd}
\|\widehat{\bS}_j(u)\|_{\min} \ge c(1- \delta/2).
\end{equation}
By (\ref{S_ju_bound}), we have
\begin{equation}
\label{pr_Omega2}
1 - \pr\big\{\Omega_{j,2}(\delta)\big\} \leq 8 \exp\Big(- \frac{c n\widebar{T}_{\mu,j}h \delta^2}{1+ \delta} \Big).
\end{equation}
Define the event $\Omega_{j,3}(\delta) = \big\{\big\|\widehat \bR_{j}(u) - \cE \{\widehat \bR_{j}(u)\} \big \| \le \delta\big\}$. Note that, under Assumption~\ref{cond_weight} with $v_{ij} \asymp (n\widebar T_{\mu,j})^{-1},$ $\sum_{i=1}^n v_{ij} \sum_{t=1}^{T_{ij}} \cE\{K_h(U_{ijt} - u)\} \le c$ and $\mu_j(\cdot)$ is uniformly bounded over $\cU$, hence $\| \cE \{\widehat \bR_{j}(u)\} \|$ is uniformly bounded over $\cU$. On the event $\Omega_{j,3}(\delta),$ we have
\begin{equation}
\label{R.bd}
\big\|\widehat \bR_{j}(u) \big\| \le c (1+ \delta).
\end{equation}
On the event $\Omega_{j,2}(\delta) \cap \Omega_{j,3}(\delta)$ with $\delta \in (0,1],$ it follows from
(\ref{ineq_decomp_mean}), (\ref{Smin.bd}) and (\ref{R.bd}) that
\begin{equation}
|\widehat{\mu}_j(u) - \widetilde{\mu}_j(u)| \le c \delta + c (1 - \delta/2)^{-1}(1+ \delta) \delta \le c_3 \delta.
\nonumber
\end{equation}
This together with concentration inequalities in (\ref{coneq.R}) and (\ref{pr_Omega2}) implies that there exist some positive universal constants $c_1$ and $c_2$ such that for any $\delta \in (0,1]$ and $u \in \cU,$
$$
\pr\Big\{|\widehat{\mu}_j(u) - \widetilde{\mu}_j(u)| \ge \delta\Big\}  \le c_2 \exp \left(- c_1 \gamma_{n,T,h,j}\delta^2\right),
$$
which completes the proof of local concentration inequality for the mean estimator.

\subsubsection{Concentration Inequality in $L_2$ Norm}
\label{ap.coneq.mean.L2}
In the proof, we need the following lemma.

\begin{lemma}
\label{l_2 lemma}
Let $X$ be a random variable. If for some constants $c_1, c_2 >0,$ $\pr(|X| > \delta) \leq c_1 \exp\{-c_2^{-1} \min(\delta^2,\delta)\}$ for any $\delta > 0,$ then for any integer $q \geq 1,$
$$
\cE(X^{2q}) \leq q! c_1 (4c_2)^q + (2q)! c_1 (4c_2)^{2q}.
$$
Conversely, if for some positive constants $a_1, a_2,$ $\cE(X^{2q}) \leq q! a_1 a_2^q + (2q)! a_1 a_2^{2q}$ for any integer $q \geq 1,$ then by letting $c_1^*=a_1$ and $c_2^*= 32(a_2+a_2^2),$ we have that
$$
\pr(|X| > \delta) \leq c_1^* \exp\{-c_2^{*-1} \min(\delta^2,\delta)\}
$$
for any $\delta > 0.$
\end{lemma}
\begin{proof}
This lemma can be proved in a similar way to Theorem~2.3 of \cite{BLM2014} and hence the proof is omitted here. In the proof, the following two inequalities are used, i.e., for any $c, \delta >0,$
$$
\frac{1}{2} \min(\delta^2, \delta) \leq \frac{\delta^2}{1+\delta} \leq \min(\delta^2, \delta)
$$
and
$$
\sqrt{\frac{c \delta}{2}} + \frac{c \delta}{2} \leq \frac{c(\delta + \sqrt{\delta^2 + 4\delta / c})}{2} \leq \sqrt{c\delta} + c\delta.
$$
\end{proof}

We are now ready to derive the $L_2$ concentration inequality of $ \|\widehat \mu_{j} - \widetilde \mu_{j}\|_2.$
Let 
\begin{flalign}
\label{add.remark}
\widetilde{\bS}_{j}(u) = (n \widebar{T}_{\mu,j})^{-1}\sum_{i=1}^{n} \sum_{t=1}^{T_{ij}} \widetilde{\bU}_{ijt} \widetilde{\bU}_{ijt}^{\T} K_h(U_{ijt}-u).
\end{flalign}
Then we have that
$\|\widehat{\bS}_{j}(u)\|_{\min} \ge c \|\widetilde{\bS}_{j}(u)\|_{\min}.$ We now give a lower bound on $\|\widetilde{\bS}_{j}(u)\|_{\min}.$ Denote $W = \sup_{u \in \cU} \big\|\widetilde{\bS}_{j}(u) - \cE\{\widetilde{\bS}_{j}(u)\}\big\|_{\tF}.$ 
Let $ \widetilde{S}_{jkl}(u) $ be the $ (k,l)$th entry of $\widetilde{\bS}_{j}(u)$ for $ k,l=1,2.$ Note that 
$\cE\{|(U_{ijt}-u)^a h^{-a} K_h(U_{ijt}-u)|\} \leq c$ and
$\cE(W) \leq 4 \max_{k,l} \cE\{\sup_{u \in \cU} |\widetilde{S}_{jkl}(u)|\}$ for $a=0,1,2.$
In an analogy to Lemma~13.5 of \cite{BLM2014}, we can show that %{\color{red}$\cE\{(n \widebar{T}_{j})^{1/2}W\} \le c,$ which means 
$\cE(W) \leq c (n \widebar{T}_{\mu,j})^{-1/2}.$ %}
Note that Lemma 13.5 of \cite{BLM2014} relies on the results presented in Lemma 13.1 of \cite{BLM2014}. Consequently, Lemma 13.5 assumes that the corresponding index set is countable in order to apply Lemma 13.1, as the supremums of the summation of indicator functions may not be measurable. However, in our specific case, each component of $\widetilde{\bS}_{j}(u) - \cE\{\widetilde{\bS}_{j}(u)\}$ is the sum of continuous functions, for which the supremums over $ \cU$ are measurable. Therefore, when we extend the index set from countable to the uncountable set $\cU$, this lemma, as well as Theorems~11.10 and 12.5 of \cite{BLM2014}, still hold true and can be applicable to our situation.
Moreover, it follows from the facts 
$\var(W) \leq \cE(W^2) \leq 4 \max_{k,l}  \var\{\sup_{u \in \cU} \widetilde{S}_{jkl}(u)\},$
$|(U_{ijt}-u)^a h^{-a}K_h(U_{ijt}-u)| \leq c h^{-1},$
$\cE\{(U_{ijt}-u)^{2a} h^{-2a} K_h^2(U_{ijt}-u)\} \leq c h^{-1}$ for $a=0,1,2,$
%\end{flalign*}
and Theorem~11.10 of \cite{BLM2014} that $\var(n \widebar{T}_{\mu,j} h W) \leq 2 \cE(n \widebar{T}_{\mu,j} h W) + \sum_{i=1}^{n} \sum_{t=1}^{T_{ij}} c h^{-1} h^2 \leq c(n \widebar{T}_{\mu,j})^{1/2} h + c n \widebar{T}_{\mu,j} h,$ which implies that the variance of $W$ is bounded by $ c (n\widebar{T}_{\mu,j}h)^{-1} \leq c \gamma_{n,T,h,j}^{-1}.$ Applying Theorem~12.5 of \cite{BLM2014} yields that there exists some positive constant $c$ such that, for any $\delta >0,$ 
\begin{equation}
\label{Z_bound}
\pr\big\{W - \cE(W) > \delta \big\} % \pr\big\{n \widebar{T}_{j} h W - \cE(n \widebar{T}_{j} h W) > n \widebar{T}_{j} h \delta \big\} 
\leq \exp \Big(- \frac{c \gamma_{n,T,h,j} \delta^2}{1+\delta}\Big).
\end{equation}
Define the event $\Omega_{j,4}(\delta)= \big\{\sup_{u \in \cU} \big\|\widetilde{\bS}_{j}(u) - \cE\{\widetilde{\bS}_{j}(u)\}\big\|_{\tF} \le \delta/2\big\}$ with $\delta \in (0,1].$ By~(\ref{Z_bound}), we obtain that there exists some constant $c>0$ such that, for any $\delta \in (0,1],$
\begin{equation}
\label{Omega4.mean}
1- \pr\left\{\Omega_{j,4}(\delta) \right\} \le 2 \exp \left(- c \gamma_{n,T,h,j} \delta^2\right).
\end{equation}
On the event $\Omega_{j,4} = \Omega_{j,4}(\delta_1)$ with $c\gamma_{n,T,h,j}^{-1/2}<\delta_1\le 1,$  $\|\widehat{\bS}_{j}(u)\|_{\min}\ge c \|\widetilde{\bS}_{j}(u)\|_{\min} \ge c(1 - \delta_1/2) \ge c/2.$ Note that $\cE\{\widehat{\bS}_j(u)\}$ is positive definite and $\| \cE \{\widehat \bR_{j}(u)\} \|$ is uniformly bounded over $\cU$. On the event $\Omega_{j,4}$, it thus follows from (\ref{ineq_decomp_mean}) and $\|\widehat{\bR}_j(u)\| \leq \|\widehat{\bR}_j(u)-\cE\{\widehat{\bR}_j(u)\}\|+\|\cE\{\widehat{\bR}_j(u)\}\| $ that
\begin{equation}
	\label{decomp_of_mean}
	%\begin{split}
	|\widehat{\mu}_j(u) - \widetilde{\mu}_j(u)|  %& \le c \|\widehat{R}_{j}(u)- \cE\{\widehat{R}_{j}(u)\}\| \\
    \leq c \|\widehat{\bR}_{j}(u)- \cE\{\widehat{\bR}_{j}(u)\}\| + c \|\widehat{\bS}_j(u) - \cE\{\widehat{\bS}_j(u)\}\|_{\tF},
%\end{split}
\end{equation}
where the positive constant $c$ does not depend on $u \in \cU.$

Combining (\ref{decomp_of_mean}) with (\ref{S_ju_bound}), (\ref{coneq.R}) and applying the first part of Lemma~\ref{l_2 lemma} yields that, for any $ u \in \cU $ and integer $ q \geq 1 $,
$$
\cE\Big\{|\widehat{\mu}_j(u) - \widetilde{\mu}_j(u)|^{2q} \Big |\Omega_{j,4}\Big\} \leq q! c \Big(\frac{4}{c \gamma_{n,T,h,j}}\Big)^q + (2q)! c \Big(\frac{4}{c \gamma_{n,T,h,j}}\Big)^{2q}.
$$
Applying the second part of Lemma~\ref{l_2 lemma} and (\ref{Omega4.mean}), we can show that, for each $\delta \in (0,1],$
$$
\pr\big(\|\hat \mu_{j}-\tilde \mu_{j}\|_2 \geq \delta\big)  \leq \pr\big(\|\hat \mu_{j}-\tilde \mu_{j}\|_2 \geq \delta \big|\Omega_{j,4}\big) + \pr\big(\Omega_{j,4}^c\big) \leq c_2 \exp\big(-{c_1 \gamma_{n,T,h,j} \delta^2}\big),\\
$$
which means that (\ref{mean_coneq_L2}) in Theorem~\ref{thm_mean_coneq} holds and completes the proof of concentration inequality for the mean estimator in $L_2$ norm.

\subsubsection{Concentration Inequality in Supremum Norm}
We will derive the uniform concentration bound of $ \sup_{u \in \cU} |\widehat \mu_{j}(u) - \widetilde \mu_{j}(u)|.$
We partition the interval $ \cU  =[0,1]$ into $ N $ subintervals $I_s$ for $s \in [N]$ of equal length.
Let $ u_s $ be the center of $I_s,$ then we have
$$
\sup_{u \in \cU} \big|\widehat \mu_{j}(u) - \widetilde \mu_{j}(u)\big| \leq  \max_{s \in [N]} \Big[\big|\widehat \mu_{j}(u_s) - \widetilde \mu_{j}(u_s)\big| +\big|\{\widehat \mu_{j}(u_s)-\widehat \mu_{j}(u)\} - \{\widetilde \mu_{j}(u_s) - \widetilde \mu_{j}(u)\}\big|\Big].
$$
We need to bound the second term.
By some calculations, it suffices to bound $\Big|\{\widehat R_{jk}(u) - \widehat R_{jk}(u_s)\} - \big[\cE\{\widehat R_{jk}(u)\} - \cE\{\widehat R_{jk}(u_s)\}\big]\Big|$ and $\Big|\{\widehat S_{jkl}(u) - \widehat S_{jkl}(u_s)\} - \big[\cE\{\widehat S_{jkl}(u)\} - \cE\{\widehat S_{jkl}(u_s)\}\big]\Big|$ for $ k,l = 1,2 $, which means that we need to bound $ \big|\widehat R_{jk}(u) - \widehat R_{jk}(u_s)\big| $ and $\big|\widehat S_{jkl}(u) - \widehat S_{jkl}(u_s)\big|. $
Let $ u \in I_s $ and consider $ |\widehat R_{j1}(u) - \widehat R_{j1}(u_s)|$ first.
Define the event $\Omega_{R,j1} = \big\{\sum_{i =1}^n v_{ij} \sum_{t = 1}^{T_{ij}} |Y_{ijt}| \leq \cE(\sum_{i=1}^n v_{ij} \sum_{t = 1}^{T_{ij}} |Y_{ijt}|) +1\big\}.$ On this event, it follows from  Assumption~\ref{cond_kernel}(ii) that
\begin{equation}
	\begin{split}
		\big|\widehat R_{j1}(u) - \widehat R_{j1}(u_s)\big| & \leq \Big|\sum_{i =1}^n v_{ij} \sum_{t = 1}^{T_{ij}} Y_{ijt} \big\{K_h(U_{ijt}-u) - K_h(U_{ijt}-u_s)\big\}\Big| \\
		& \leq \frac{c |u-u_s|}{h^2} \sum_{i =1}^n v_{ij} \sum_{t = 1}^{T_{ij}} \Big|Y_{ijt}\Big| \leq \frac{c}{Nh^2} \left\{\cE\Big(\sum_{i=1}^n v_{ij} \sum_{t = 1}^{T_{ij}} \Big|Y_{ijt}\Big|\Big) +1\right\} \leq \frac{c }{Nh^2}.
	\end{split} \nonumber
\end{equation}
Applying similar techniques as above, we can define events $ \Omega_{R,jk} $ and $ \Omega_{S,jkl} $ for $ k,l = 1,2 $.
On the intersection of these events, we can obtain that $|\widehat R_{jk}(u) - \widehat R_{jk}(u_s)| \leq c (Nh^2)^{-1}$ and $|\widehat S_{jkl}(u) - \widehat S_{jkl}(u_s)| \leq c (Nh^2)^{-1}$.
Combing	the above results, we have
\begin{equation}
	\sup_{u \in \cU} \big|\widehat \mu_{j}(u) - \widetilde \mu_{j}(u)\big| \leq  \max_{s \in [N]} \big|\widehat \mu_{j}(u_s) - \widetilde \mu_{j}(u_s)\big| + \frac{c}{Nh^2}.
	\nonumber
\end{equation}
Applying Hoeffding's inequality, we obtain that
%\begin{flalign*} 
$\pr(\Omega_{R,jk}^{\rm c}) \leq \exp\{- (2 \sum_{i=1}^n c^2 v_{ij}^2 T_{ij}^2)^{-1}\} = \exp(-cn) \leq \exp(-c \gamma_{n,T,h,j})$ and
$\pr(\Omega_{S,jkl}^{\rm c}) \leq \exp\{- (2 \sum_{i=1}^n c^2 v_{ij}^2 T_{ij}^2)^{-1}\} = \exp(-cn) \leq \exp(-c \gamma_{n,T,h,j})$ for $k,l=1,2.$
%\end{flalign*}
It follows from the above results and the union bound of probability with the choice of $N = \lfloor c (h^2 \delta)^{-1} \rfloor$ that there exist some positive constants $c_1$ and $c_2$ such that, for any $ \delta \in (0,1],$
\begin{equation}
	\label{1}
	\pr\Big\{\sup_{u \in \cU} |\widehat \mu_{j}(u) - \widetilde \mu_{j}(u)| \geq \delta\Big\} \leq \frac{c_2}{h^2 \delta} \exp(-c_1 \gamma_{n,T,h,j} \delta^2).
\end{equation}
Take arbitrarily small $\epsilon_1>0.$ If $ n^{\epsilon_1}  \gamma_{n,T,h,j} \delta^2 \geq 1,$ then the right side of (\ref{1}) reduces to $ c_2 \{n^{\epsilon_1} \gamma_{n,T,h,j}\}^{1/2} h^{-2} \exp(-c_1 \gamma_{n,T,h,j} \delta^2).$
If $ n^{\epsilon_1} \gamma_{n,T,h,j} \delta^2 \leq 1,$ we can choose $ c_2 $ and $ n^{\epsilon_1} > c $ such that $ c_2 \exp(-c_1 c^{-1}) \geq 1 $ and the same bound $c_1 \{n^{\epsilon_1} \gamma_{n,T,h,j}\}^{1/2} h^{-2} \exp(-c_1 \gamma_{n,T,h,j} \delta^2)$ can still be used.
Hence (\ref{mean_coneq_sup}) in Theorem~\ref{thm_mean_coneq} holds, which completes the proof of concentration inequality for the mean estimator in supremum norm. $\hfill\blacksquare$

\subsection{Proof of Theorem~\ref{thm_cov_coneq}}
\label{pf.thm.cov}
We organize the proof %of Theorem~\ref{thm_cov_coneq}
in four steps. First, we will define $ \widehat{\Sigma}(\cdot,\cdot),$ $ \widetilde{\Sigma}(\cdot,\cdot)$ and obtain the decomposition of $ \widehat{\Sigma}(\cdot,\cdot) - \widetilde{\Sigma}(\cdot,\cdot).$ Second, we will prove the local concentration inequality for fixed $(u,v) \in \cU^2.$ Third, we will prove the concentration inequality in Hilbert--Schmidt norm. Finally, we will prove the concentration inequality in supremum norm.

\subsubsection{Definition and Decomposition}
Without loss of generality, let $h_{\sSigma,jk} = h $ for $(j,k) \in [p]^2$ and denote $\tilde \be_0 = (1,0,0)^{\T},$ $ \widetilde{\bU}_{ijkts} = \big\{1, (U_{ijt} - u)/h, (U_{iks} - v)/h\big\}^{\T}.$ For $ j=k $, let
\begin{equation}
	\begin{split}
		& \widehat{\bXi}_{jj}(u, v) = \sum_{i =1}^{n} w_{ijj} \sum_{1 \leq t \neq s \leq T_{ij}} \widetilde{\bU}_{ijjts} \widetilde{\bU}_{ijjts}^{\T} K_{h}(U_{ijt} - u) K_{h}(U_{ijs} - v), \\
		& \widehat{\bZ}_{jj}(u, v) = \sum_{i =1}^{n} w_{ijj} \sum_{1 \leq t \neq s \leq T_{ij}} \widetilde{\bU}_{ijjts} \Theta_{ijjts} K_{h}(U_{ijt} - u) K_{h}(U_{ijs} - v). \\
	\end{split}
    \nonumber
\end{equation}
For $ j \neq k $, let
\begin{equation}
	\begin{split}
		& \widehat{\bXi}_{jk}(u, v) = \sum_{i =1}^{n}  w_{ijk} \sum_{t=1}^{T_{ij}} \sum_{s=1}^{T_{ik}} \widetilde{\bU}_{ijkts} \widetilde{\bU}_{ijkts}^{\T} K_{h}(U_{ijt} - u) K_{h}(U_{iks} - v), \\
		& \widehat{\bZ}_{jk}(u, v) = \sum_{i =1}^{n}  w_{ijk} \sum_{t=1}^{T_{ij}} \sum_{s=1}^{T_{ik}} \widetilde{\bU}_{ijkts} \Theta_{ijkts} K_{h}(U_{ijt} - u) K_{h}(U_{iks} - v). \\
	\end{split}
 \nonumber
\end{equation}
A simple calculation yields that $ \widehat{\Sigma}_{jk} (u, v) = \tilde \be_0^{\T} \{\widehat{\bXi}_{jk}(u, v)\}^{-1} \widehat{\bZ}_{jk}(u, v).$
Let
\begin{equation}
\label{Sigma_tilde}
\widetilde{\Sigma}_{jk} (u, v) = \tilde \be_0^{\T} \big[\cE \{ \widehat{\bXi}_{jk}(u, v)\}\big]^{-1} \cE \{\widehat{\bZ}_{jk}(u, v)\}.
\end{equation}
We can decompose $\widehat{\Sigma}_{jk} (u, v) - \widetilde{\Sigma}_{jk} (u, v)$ as
\begin{equation}
\label{decomp_cov}
	\begin{split}
		\widehat{\Sigma}_{jk} (u, v) - \widetilde{\Sigma}_{jk} (u, v) = & \tilde \be_0^{\T} (\{\widehat{\bXi}_{jk}(u, v)\}^{-1} - [\cE \{ \widehat{\bXi}_{jk}(u, v)\}]^{-1}) \widehat{\bZ}_{jk}(u, v) \\
	& + \tilde \be_0^{\T} [\cE \{ \widehat{\bXi}_{jk}(u, v)\}]^{-1} [\widehat{\bZ}_{jk}(u, v) - \cE \{\widehat{\bZ}_{jk}(u, v)\}],
	\end{split} \nonumber
\end{equation}
which further implies that
\begin{equation}
\label{ineq_decomp_cov}
	\begin{split}
	&	\big|\widehat{\Sigma}_{jk} (u, v) - \widetilde{\Sigma}_{jk} (u, v)\big|\\
 \le &    \|\cE\{\widehat{\bXi}_{jk}(u, v)\}\|^{-1}_{\min} \|\widehat{\bZ}_{jk}(u, v) - \cE\{\widehat{\bZ}_{jk}(u, v)\}\| \\
		& +  \|\cE\{\widehat{\bXi}_{jk}(u, v)\}\|^{-1}_{\min} \|\widehat{\bXi}_{jk}(u, v)\|^{-1}_{\min} \| \|\widehat{\bZ}_{jk}(u, v)\| \|\widehat{\bXi}_{jk}(u, v) - \cE\{\widehat{\bXi}_{jk}(u, v)\} \|_{\tF}.
	\end{split}
\end{equation}
In the following, we will prove the concentration results for case $j \neq k,$ and the results for the case $j=k$ can be proved in a similar manner.

\subsubsection{Local Concentration Inequality}
We will firstly show that there exists some positive constant $c$ (independent of $u,v$) such that for any $\delta >0$ and $(u,v) \in \cU^2,$
\begin{equation}
\label{Xi_bound}
	\pr\Big\{\big\|\widehat{\bXi}_{jk}(u,v) - \cE\{\widehat{\bXi}_{jk}(u,v)\}\big\|_{\tF} \geq \delta \Big\}\leq 18 \exp\Big(- \frac{c \nu_{n,T,h,jk} \delta^2}{1+ \delta} \Big).
\end{equation}

For $ m,l=1,2,3,$ let $\widehat{\Xi}_{jkml}(u,v) $ be the $ (m,l)$th entry of $\widehat{\bXi}_{jk}(u,v).$
It follows from Assumptions~\ref{cond_T_bd} and \ref{cond_kernel} that for any integer $q = 2,3,\dots$ and $s, s'=0, 1, 2,$
\begin{equation}
	\label{iequ.kernel.cov}
	\begin{split}
	    & \cE\left\{\Big|\Big(\frac{U_{ijt}-u}{h}\Big)^s \Big(\frac{U_{ikt'}-v}{h}\Big)^{s'} K_h(U_{ijt}-u) K_h(U_{ikt'}-v)\Big|^{q}\right\} \\
     & \leq \int h^{-2q} K^q\Big(\frac{t-u}{h}\Big) K^q\Big(\frac{t'-v}{h}\Big) \Big|\frac{t-u}{h} \Big|^{sq} \Big|\frac{t'-v}{h} \Big|^{sq} f_U(t) f_U(t') dt dt' \leq c h^{2-2q}. \\
	\end{split}
\end{equation}
Note that Assumption~\ref{cond_weight} implies that the weights are of the same order $w_{ijk} \asymp (n\widebar T_{\sSigma,jk}^2)^{-1}.$ By (\ref{iequ.kernel.cov}),  $$\sum_{i=1}^{n} \sum_{t=1}^{T_{ij}} \sum_{t'=1}^{T_{ik}} \cE\left\{|(\frac{U_{ijt} - u}{h})^s (\frac{U_{ikt'} - v}{h})^{s'} K_{h}(U_{ijt} - u) K_{h}(U_{ikt'} - v)|^2\right\} \leq c n \widebar T_{\sSigma,jk}^2 h^{-2},$$ $$\sum_{i=1}^{n} \sum_{t=1}^{T_{ij}} \sum_{t'=1}^{T_{ik}} \cE\left\{|(\frac{U_{ijt} - u}{h})^s (\frac{U_{ikt'} - v}{h})^{s'} K_{h}(U_{ijt} - u) K_{h}(U_{ikt'} - v)|^q\right\} \leq 2^{-1} q! c n \widebar T_{\sSigma,jk}^2 h^{2-2q}$$ for $q \geq 3.$ Applying the Bernstein inequality yields that there exists some positive constant $c$ (independent of $u,v$) such that for any $\delta >0$ and $ (u,v) \in \cU^2,$
$$
\pr\Big\{\big|\widehat{\Xi}_{jkml}(u,v) - \cE\{\widehat{\Xi}_{jkml}(u,v)\}\big| \geq \delta \Big\}\leq 2 \exp\Big(- \frac{c \nu_{n,T,h,jk} \delta^2}{1+ \delta} \Big),
$$
for $ m,l=1,2,3,$ which, by the union bound of probability, implies that (\ref{Xi_bound}) holds.

For $ m =1,2,3, $ let $ \widehat{Z}_{jkm}(u,v) $ be the $m$th element of $ \widehat{\bZ}_{jk}(u,v) $. We will next show that, there exits some positive constant $c$ (independent of $u,v$) such that for any $ \delta > 0 $ and $ (u,v) \in \cU^2,$
\begin{equation}
\label{coneq.Z}
	\pr\Big\{\big|\widehat{Z}_{jkm}(u,v) - \cE\{\widehat{Z}_{jkm}(u,v)\}\big| \geq \delta\Big\} \leq c \exp\Big(- \frac{c \nu_{n,T,h,jk} \delta^2}{1+ \delta}\Big),
\end{equation}
where $\nu_{n,T,h,jk}=n (1 \wedge \widebar T_{\sSigma,jk}^2 h^2).$
We only need to consider the case $m=1$, while the results for cases $m=2,3$ can be proved similarly. Denote that 
$$\zeta_{ijkts} = \{Y_{ijt} - \mu_j(U_{ijt})\} \{Y_{iks} - \mu_k(U_{iks})\} - \Sigma_{jk}(U_{ijt},U_{iks}),$$
%{\color{red} $\zeta_{ijkts}$ there should be $\{Y_{ijt} - \widehat \mu_j(U_{ijt})\} \{Y_{iks} - \widehat \mu_k(U_{iks})\} - %\Sigma_{jk}(U_{ijt},U_{ikt'})$ this needs edit by Shaojun}
$$
\widehat{Z}_{jk4}(u,v) = \sum_{i=1}^{n} w_{ijk} \sum_{t=1}^{T_{ij}} \sum_{s=1}^{T_{ik}} K_{h}(U_{ijt} - u) K_{h}(U_{iks} - v) \Sigma_{jk}(U_{ijt},U_{iks}),
$$
$$\widehat{Z}_{jk5}(u,v) = \sum_{i=1}^{n} w_{ijk} \sum_{t=1}^{T_{ij}} \sum_{s=1}^{T_{ik}} K_h(U_{ijt}-u) K_{h}(U_{iks} - v) \zeta_{ijkts}.$$ Then we rewrite $\widehat{Z}_{jk1}(u,v) - \cE\{\widehat{Z}_{jk1}(u,v)\}$ as
\begin{equation}
\label{decompose_Z1}
\widehat{Z}_{jk1}(u,v) - \cE\{\widehat{Z}_{jk1}(u,v)\} = \widehat{Z}_{jk4}(u,v) - \cE\{\widehat{Z}_{jk4}(u,v)\} + \widehat{Z}_{jk5}(u,v).
\end{equation}

Following the same procedure to prove (\ref{Xi_bound}) with the aid of the Bernstein inequality, we can obtain that there exists some positive constant $c$ such that for any $ \delta>0$ and $ (u,v) \in \cU^2,$
\begin{equation}
\label{concentration_Z_4}
	\pr\Big\{\big |\widehat{Z}_{jk4}(u,v) - \cE\{\widehat{Z}_{jk4}(u,v)\}\big| \geq \delta \Big\}\leq 2 \exp\Big(- \frac{c \nu_{n,T,h,jk} \delta^2}{1+ \delta} \Big).
\end{equation}

Now we consider the tail behavior of $\widehat{Z}_{jk5}(u,v).$ Define the event $\widetilde V_{jk} = \{(U_{ijt}, U_{iks}), t \in [T_{ij}], s \in [T_{ik}], i \in [n]\}.$
Note a random variable $X$ is sub-exponential if there exist positive constants $c_1$ and $c_2$ such that  $\cE(\exp[\lambda \{X-\cE(X)\}]) \leq \exp(c_1^2 \lambda^2/2)$ for all $|\lambda| < c_2^{-1}.$ The sub-Gaussianities under Assumption~\ref{cond_subG} implies that, conditional on the event $\widetilde V_{jk},$ $Y_{ijt} - \mu_j(U_{ijt})$ and $Y_{iks} - \mu_k(U_{iks})$ are sub-Gaussian random variables, then $\{Y_{ijt} - \mu_j(U_{ijt})\}\{Y_{iks} - \mu_k(U_{iks})\}$ is a sub-exponential random variable, %$\zeta_{ijkts} = \{Y_{ijt} - \mu_j(U_{ijt})\}\{Y_{iks} - \mu_k(U_{iks})\} - \cE[\{Y_{ijt} - \mu_j(U_{ijt})\}\{Y_{iks} - \mu_k(U_{iks})\}]$,
and hence we can obtain the Bernstein-type bound
%we have $\cE[\{Y_{ijt} - \mu_j(U_{ijt})\}^{2q}] \le (2q)! 2^{-q} (q!)^{-1} c^{2q}$ and $\cE[\{Y_{iks} - \mu_k(U_{iks})\}^{2q}] \le (2q)! 2^{-q} (q!)^{-1} c^{2q},$ which yields that $\cE[\{Y_{ijt} - \mu_j(U_{ijt})\}^{q} \{Y_{iks} - \mu_k(U_{iks})\}^{q}] \le (2q)! 2^{-q} (q!)^{-1} c^{2q}.$
%For $\sup_{q \ge 2} (\cE[\{Y_{ijt} - \mu_j(U_{ijt})\}^{q} \{Y_{iks} - \mu_k(U_{iks})\}^{q}] / q!)^{1/q} \le 2c^2$ is finite, $\zeta_{ijkts}$ is a sub-exponential random variable, so we have
%The sub-Gaussian property in Assumption~\ref{cond_subG} implies that there exists some constant $c>0$ such that
$
\cE\{\exp(\lambda \zeta_{ijkts}) \big|\widetilde V_{jk}\big\}\le \exp\big\{(1-c\lambda)^{-1}{c \lambda^2}\big\}
$
for any $\lambda \in (0, c^{-1}).$
Rewrite $ \widehat{Z}_{jk5}(u,v) = \sum_{i=1}^n w_{ijk} \phi_{ijk1}(u,v) $ with $ \phi_{ijk1}(u,v) = \sum_{t=1}^{T_{ij}} \sum_{s=1}^{T_{ik}}  K_h(U_{ijt}-u) K_h(U_{iks}-v) \zeta_{ijkts}.$
Note that, for each $i,$ Assumption~\ref{cond_weight} implies that $w_{ijk} \asymp (n\widebar T_{\sSigma,jk}^2)^{-1}.$ If $\sum_{t=1}^{T_{ij}} \sum_{s=1}^{T_{ik}} K_h(U_{ijt} - u) K_h(U_{iks}-v)> 0$ holds, it follows from Jensen's inequality and the above result that
\begin{equation}
	\begin{split}
		& \cE\Big[\exp\{\lambda w_{ijk} \phi_{ijk1}(u,v)\} \Big|\widetilde V_{jk} \Big] \\
		\leq & \frac{1}{\sum_{t=1}^{T_{ij}} \sum_{s=1}^{T_{ik}}  K_h(U_{ijt}-u) K_h(U_{iks}-v)} \sum_{t=1}^{T_{ij}} \sum_{s=1}^{T_{ik}} \Big( K_h(U_{ijt}-u) K_h(U_{iks}-v) \\
  & \times \cE\Big[\exp\Big\{\lambda w_{ijk} \zeta_{ijkts} \sum_{t'=1}^{T_{ij}} \sum_{s'=1}^{T_{ik}}  K_h(U_{ijt'}-u) K_h(U_{iks'}-v)  \Big\} \Big|\widetilde V_{jk}\Big] \Big)\\
		\le & \exp\left[\frac{c \lambda^2  (n\widebar T_{\sSigma,jk}^2)^{-2}\Big\{\sum_{t,s} K_h(U_{ijt} - u) K_h(U_{iks}-v)\Big\}^2}{1 - c\lambda  (n\widebar T_{\sSigma,jk}^2)^{-1} \sum_{t,s} K_h(U_{ijt} - u) K_h(U_{iks}-v)}\right]. \\
	\end{split}
	\nonumber
\end{equation}
where $0<\lambda (n\widebar T_{\sSigma,jk}^2)^{-1} \sum_{t=1}^{T_{ij}} \sum_{s=1}^{T_{ik}}  K_h(U_{ijt}-u) K_h(U_{iks}-v) < c^{-1}. $
It is obvious that the above inequality still holds even if $\sum_{t=1}^{T_{ij}} \sum_{s=1}^{T_{ik}} K_h(U_{ijt} - u) K_h(U_{iks}-v) = 0.$ Assumption~\ref{cond_timepoints} implies that the number of nonzero terms in $\sum_{t=1}^{T_{ij}} \sum_{s=1}^{T_{ik}} K_h(U_{ijt} - u) K_h(U_{iks}-v)$ has an upper bound $c (1\vee \widebar T_{\sSigma,jk}^2h^2),$ which yields that $\sum_{t=1}^{T_{ij}} \sum_{s=1}^{T_{ik}} K_h(U_{ijt} - u) K_h(U_{iks}-v) \le c h^{-2} (1\vee \widebar T_{\sSigma,jk}^2 h^2).$ Therefore, for any $\lambda$ satisfying $0< \lambda (n \widebar{T}_{\sSigma,jk}^2 h^2)^{-1} (1\vee \widebar T_{\sSigma,jk}^2h^2)  < c^{-1}$ for some constant $c>0$, we obtain that
\begin{equation}
\begin{split}
& \cE\left[\exp\Big\{\lambda \sum_{i=1}^n w_{ijk} \phi_{ijk1}(u,v) \Big\} \Big|\widetilde V_{jk} \right]\\
%&  \le \exp\left[\lambda^2 c^2 \sum_{i=1}^n v_{ij}^2\Big\{\sum_{t=1}^{T_{ij}} K_h(U_{ijt} - u)\Big\}^2 \right]\\
&  \le \exp\left\{\frac{c \lambda^2(n\widebar T_{\sSigma,jk}^2)^{-2} h^{-2} (1\vee \widebar T_{\sSigma,jk}^2h^2)  \sum_{i=1}^n \sum_{t=1}^{T_{ij}} \sum_{s=1}^{T_{ik}} K_h(U_{ijt} - u)K_h(U_{iks} - v)}{1 - c \lambda (n \widebar{T}_{\sSigma,jk}^2 h^2)^{-1} (1\vee \widebar T_{\sSigma,jk}^2h^2)} \right\}.
\end{split}
\nonumber
\end{equation}
For any $\delta>0,$ define the event
$$
\Lambda_{jk,1}(\delta) = \left\{\sum_{i=1}^n \sum_{t=1}^{T_{ij}} \sum_{s=1}^{T_{ik}} K_h(U_{ijt} - u) K_h(U_{iks} - v) \le c (1 + \delta) n\widebar{T}_{\sSigma,jk}^2\right\}.
$$
We have
\begin{equation}
\begin{split}
 \cE\left[\exp\Big\{\lambda \sum_{i=1}^n w_{ijk} \phi_{ijk1}(u,v) \Big\} \Big|\Lambda_{jk,1}(\delta) \right]
 \le \exp\left\{\frac{c \lambda^2(1+\delta) (n\widebar T_{\sSigma,jk}^2h^2)^{-1} (1\vee \widebar T_{\sSigma,jk}^2h^2)}{1 - c \lambda (n \widebar{T}_{\sSigma,jk}^2 h^2)^{-1} (1\vee \widebar T_{\sSigma,jk}^2h^2)} \right\}.
\end{split}
\nonumber
\end{equation}
Consequently, we obtain that
\begin{equation}
	\label{sum_exp_phi}
	%\begin{split}
		\pr\Big\{\sum_{i=1}^{n} w_{ijk} \phi_{ijk1}(u)\geq \delta \Big|\Lambda_{jk,1}(\delta) \Big\}
		 \leq \exp\left\{-\lambda \delta + \frac{c \lambda^2(1+\delta) (n\widebar T_{\sSigma,jk}^2h^2)^{-1} (1\vee \widebar T_{\sSigma,jk}^2h^2)}{1 - c \lambda (n \widebar{T}_{\sSigma,jk}^2 h^2)^{-1} (1\vee \widebar T_{\sSigma,jk}^2h^2)} \right\}.
	%\end{split}
\end{equation}
With the choice of $ \lambda = n \widebar{T}_{\sSigma,jk}^2 h^2 \delta \big\{2 c (1 + \delta)(1 \vee \widebar{T}_{\sSigma,jk}^2 h^2)+c \delta (1 \vee \widebar{T}_{\sSigma,jk}^2 h^2)\big\}^{-1},$ (\ref{sum_exp_phi}) reduces to
\begin{equation}
	\label{sum_exp_phi2}
		\pr\Big\{\sum_{i=1}^{n} w_{ijk} \phi_{ijk1}(u,v) \geq \delta \Big|\Lambda_{jk,1}(\delta)\Big\}  \leq \exp\Big\{- \frac{c\nu_{n,T,h,jk}\delta^2}{1 + \delta}\Big\},
\end{equation}
where the constant $c$ is chosen to satisfy $c \lambda (n \widebar{T}_{\sSigma,jk}^2 h^2)^{-1} (1\vee \widebar T_{\sSigma,jk}^2h^2) \le 1/2.$
Note that 
\begin{flalign*}
\sum_{i=1}^n \sum_{t=1}^{T_{ij}} \sum_{s=1}^{T_{ik}} \cE\big\{K_h(U_{ijt} - u)K_h(U_{iks} - v)\big\} \le c n \widebar{T}_{\sSigma,jk}^2.
\end{flalign*}
By the Bernstein inequality, we obtain that there exists some positive constant $c$ such that for any $\delta >0$
\begin{eqnarray*}
& & \pr\left( \sum_{i=1}^n \sum_{t=1}^{T_{ij}} \sum_{s=1}^{T_{ik}} \Big[ K_h(U_{ijt} - u) K_h(U_{iks} - v) - \cE\big\{K_h(U_{ijt} - u)K_h(U_{iks} - v)\big\} \Big] \ge n\widebar{T}_{\sSigma,jk}^2 \delta \right)\\
&\le& \exp \left( \frac{-c \nu_{n,T,h,jk} \delta^2}{1+\delta}\right),
\end{eqnarray*}
which implies that
\begin{equation}
\label{kernel_concentration2}
1- \pr\{\Lambda_{jk,1}(\delta)\} \le \exp \Big(-\frac{c \nu_{n,T,h,jk} \delta^2}{1+\delta}\Big).
\end{equation}
Combining (\ref{sum_exp_phi2}) and (\ref{kernel_concentration2}), we obtain that there exists some constant $c>0$ such that for any $\delta >0,$
\begin{equation}
    \pr\Big\{\widehat Z_{jk5}(u,v) \geq \delta\Big\}  \leq \pr\Big\{\widehat Z_{jk5}(u,v) \geq \delta \Big| \Lambda_{jk,1}(\delta)\Big\} + \pr\Big\{\Lambda_{jk,1}(\delta)^c\Big\} \leq 2\exp \Big(- \frac{c\nu_{n,T,h,jk}\delta^2}{1 + \delta }\Big),
  \nonumber
\end{equation}
which leads to
\begin{equation}
		\pr\Big\{\big|\widehat Z_{jk5}(u,v)\big|\geq \delta\Big\} \leq 4\exp \Big(- \frac{c\nu_{n,T,h,jk}\delta^2}{1 + \delta }\Big).
  \nonumber
\end{equation}
It follows from the above, (\ref{decompose_Z1}) and (\ref{concentration_Z_4}) that for each $\delta >0$ and $(u,v)\in \cU^2,$ there exists some positive constant $c$ such that
$$
\pr\Big\{\big|\widehat Z_{jk1}(u,v) - \cE \{\widehat Z_{jk1}(u,v)\}\big|\geq \delta\Big\} \leq 6\exp \Big(- \frac{c\nu_{n,T,h,jk}\delta^2}{1 + \delta}\Big).
$$
%and as a result,
%\begin{equation}
%\label{conseq_Z}
%\pr\Big\{\big\|\widehat Z_{jk}(u) - E \{\widehat Z_{jk}(u)\}\big\|\geq \delta\Big\} \leq 18 \exp \Big\{- \frac{cn (1\wedge \widebar{T}_{jk}^2 h^2)\delta^2}{1 + \delta}\Big\}.
%\end{equation}
Define the event $\Lambda_{jk,2}(\delta) = \big\{\|\widehat{\bXi}_{jk}(u,v) - \cE\{\widehat{\bXi}_{jk}(u,v)\}\|_{\tF} \le \delta/2\big\}.$ Note that $\cE\{\widehat{\bXi}_{jk}(u,v)\}$ is positive definite.
%It follows from Condition~\ref{cond_kernel} that
%$E\big\{ K_h(U_{ijt}-u)K_h(U_{iks}-v)\big\}=1,$ $E\big\{ h^{-1}(U_{ijt}-u)K_h(U_{ijt}-u)K_h(U_{iks}-v)\big\}=0,$ $E\big\{h^{-1}(U_{iks}-v)K_h(U_{ijt}-u)K_h(U_{iks}-v)\big\}=0,$ $E\big\{ h^{-1}(U_{ijt}-u)h^{-1}(U_{iks}-v)K_h(U_{ijt}-u)K_h(U_{iks}-v)\big\}=0,$ $c <E\big\{h^{-2}(U_{ijt}-u)^2 K_h(U_{ijt}-u) K_h(U_{iks}-v)\big\} <\infty,$
%and $c <E\big\{h^{-2}(U_{iks}-v)^2 K_h(U_{ijt}-u) K_h(U_{iks}-v)\big\} <\infty,$
%which implies that $E\{\widehat{\Xi}_{jk}(u,v)\}$ is diagonal and positive definite with eigenvalues bounded away from $0.$
On the event $\Lambda_{jk,2}(\delta)$ with $\delta \in (0,1],$ we obtain that
\begin{equation}
\label{Ximin.bd}
\|\widehat{\bXi}_{jk}(u,v)\|_{\min} \ge c(1- \delta/2).
\end{equation}
By (\ref{Xi_bound}), we have
\begin{equation}
\label{pr_Omega2-2}
1 - \pr\big\{\Lambda_{jk,2}(\delta)\big\} \leq 18 \exp\Big(- \frac{c \nu_{n,T,h,jk} \delta^2}{1+ \delta} \Big).
\end{equation}
Define the event $\Lambda_{jk,3}(\delta) = \big\{\big\|\widehat \bZ_{jk}(u,v) - \cE \{\widehat \bZ_{jk}(u,v)\} \big \| \le \delta\big\}$. Note that, under Assumption~\ref{cond_weight} with $w_{ijk} \asymp (n\widebar T_{\sSigma,jk}^2)^{-1},$ $\sum_{i=1}^n w_{ijk} \sum_{t=1}^{T_{ij}} \sum_{s=1}^{T_{ik}} \cE\{K_h(U_{ijt} - u)K_{h}(U_{iks}-v)\} \le c,$ hence $\| \cE \{\widehat \bZ_{jk}(u,v)\} \|$ is uniformly bounded over $\cU^2$. On the event $\Omega_{jk,3}(\delta),$ we have
\begin{equation}
\label{Z.bd}
\big\|\widehat \bZ_{jk}(u,v) \big\| \le c (1+ \delta).
\end{equation}
On the event $\Lambda_{jk,2}(\delta) \cap \Lambda_{jk,3}(\delta)$ with $\delta \in (0,1],$ it follows from
(\ref{ineq_decomp_cov}), (\ref{Ximin.bd}) and (\ref{Z.bd}) that
\begin{equation}
\big|\widehat{\Sigma}_{jk}(u,v) - \widetilde{\Sigma}_{jk}(u,v)\big| \le c \delta + c (1 - \delta/2)^{-1}(1+ \delta) \delta \le c_3 \delta.
\nonumber
\end{equation}
This together with concentration inequalities in (\ref{coneq.Z}) and (\ref{pr_Omega2-2}) implies that there exist some positive universal constants $c_1$ and $c_2$ such that for any $\delta \in (0,1]$ and $(u,v)\in \cU^2,$
$$
\pr\Big\{\big|\widehat\Sigma_{jk}(u,v) - \widetilde\Sigma_{jk}(u,v)\big| \ge \delta\Big\}  \le c_2 \exp \left(- c_1 \nu_{n,T,h,jk}\delta^2\right),
$$
which completes the proof of local concentration inequality for the covariance estimator.

\subsubsection{Concentration Inequality in Hilbert--Schmidt Norm}
We will derive the $L_2$ concentration inequality of $\|\widehat\Sigma_{jk} - \widetilde \Sigma_{jk}\|_{\cS}.$
Let $$\widetilde{\bXi}_{jk}(u,v) = (n \widebar{T}_{\sSigma,jk}^2)^{-1}\sum_{i=1}^{n} \sum_{t=1}^{T_{ij}} \sum_{s = 1}^{T_{ik}} \widetilde{\bU}_{ijkts} \widetilde{\bU}_{ijkts}^{\T} K_h(U_{ijt}-u)K_h(U_{iks}-v).$$ Then we have that
$\|\widehat{\bXi}_{jk}(u,v)\|_{\min} \ge c \|\widetilde{\bXi}_{jk}(u,v)\|_{\min}.$ Similar to Appendix~\ref{ap.coneq.mean.L2}, we will give a lower bound on $\|\widetilde{\bXi}_{jk}(u,v)\|_{\min}.$ Denote $\widetilde W = \sup_{(u,v) \in \cU^2} \big\|\widetilde{\bXi}_{jk}(u,v) - \cE\{\widetilde{\bXi}_{jk}(u,v)\}\big\|_{\tF}.$ 
For $ t,s=1,2,3,$ let $ \widetilde{\bXi}_{jkts}(u,v) $ be the $ (t,s)$th entry of $\widetilde{\bXi}_{jk}(u,v).$ Note that 
$
\cE\{|(U_{ijt}-u)^{a} h^{-a} (U_{iks}-v)^{b} h^{-b} K_h(U_{ijt}-u) K_h(U_{iks}-v)|\} \leq c
$ for $a, b=0,1,2,$
and, moreover,
$\cE(\widetilde W) \leq 6\max_{t,s} \cE\{ \sup_{(u,v) \in \cU^2} |\widetilde{\bXi}_{jkts}(u,v)|\}.$ 
%\end{flalign*}
In an analogy to Lemma~13.5 of \cite{BLM2014} and by the similar arguments below (\ref{add.remark}) in Appendix~\ref{ap.coneq.mean.L2}, we can show that %{\color{red}$\cE\{ (n \widebar{T}_{jk}^2)^{1/2} \widetilde W\} \leq c,$ which means 
$\cE(\widetilde W) \leq c (n \widebar{T}_{\sSigma,jk}^2)^{-1/2}.$ %}
In addition, it follows from the facts 
%\begin{flalign*}
$\var(W) \leq \cE(W^2) \leq 9 \max_{t,s}  \var\{\sup_{(u,v) \in \cU^2} \widetilde{\bXi}_{jkts}(u,v)\},$
$|(U_{ijt}-u)^{a} h^{-a} (U_{iks}-v)^{b} h^{-b} K_h(U_{ijt}-u) K_h(U_{iks}-v)| \leq c h^{-2},$
$\cE\{(U_{ijt}-u)^{2a} h^{-2a} (U_{iks}-v)^{2b} h^{-2b} K_h^2(U_{ijt}-u) K_h^2(U_{iks}-v)\} \leq c h^{-2}$ for $a,b=0,1,2,$
%\end{flalign*}
and Theorem~11.10 of \cite{BLM2014} that $\var(n \widebar{T}_{\sSigma,jk}^2 h^2 \widetilde W) \leq 2 \cE(n \widebar{T}_{\sSigma,jk}^2 h^2 \widetilde W) + \sum_{i=1}^{n} \sum_{t=1}^{T_{ij}} \sum_{s = 1}^{T_{ik}} c h^{-2} h^4 \leq c (n \widebar{T}_{\sSigma,jk}^2)^{1/2} h^2 + c n \widebar{T}_{\sSigma,jk}^2 h^2,$ which implies that the variance of $\widetilde W$ is bounded by $ c (n \widebar{T}_{\sSigma,jk}^2 h^2)^{-1} \leq c (\nu_{n,T,h,jk})^{-1}.$
Noting the similar arguments below (\ref{add.remark}) and applying Theorem~12.5 of \cite{BLM2014}, we obtain that there exists some positive constant $c$ such that, for any $\delta >0,$ 
\begin{equation}
\label{Z_bound2}
 \pr\Big\{\widetilde W - \cE(\widetilde W) > \delta \Big\} %= \pr\Big\{n \widebar{T}_{jk}^2 h^2 \widetilde W - \cE(n \widebar{T}_{jk}^2 h^2 \widetilde W) > n \widebar{T}_{jk}^2 h^2 \delta \Big\} 
 \leq \exp \Big(- \frac{c \nu_{n,T,h,jk}\delta^2}{1+\delta}\Big).
\end{equation} 
Define the event $\Lambda_{jk,4}(\delta)= \big\{\sup_{(u,v) \in \cU^2} \big\|\widetilde{\bXi}_{jk}(u,v) - \cE\{\widetilde{\bXi}_{jk}(u,v)\}\big\|_{\tF} \le \delta/2\big\}$ with $\delta \in (0,1].$ By~(\ref{Z_bound2}), we obtain that there exists some constant $c>0$ such that, for any $\delta \in (0,1],$
$$
1- \pr\left\{\Lambda_{jk,4}(\delta) \right\} \le 2 \exp \left(- c \nu_{n,T,h,jk}\delta^2\right).
$$
On the event $\Lambda_{jk,4} = \Lambda_{jk,4}(\tilde \delta_1)$ with $c(\nu_{n,T,h,jk})^{-1/2}<\tilde \delta_1\le 1,$ we have $\|\widehat{\bXi}_{jk}(u,v)\|_{\min}\ge c \|\widetilde{\bXi}_{jk}(u,v)\|_{\min} \ge c(1 - \tilde \delta_1/2) \ge c/2.$ Notice that $\cE\{\widehat{\bXi}_{jk}(u,v)\}$ is positive definite and $\| \cE \{\widehat \bZ_{jk}(u,v)\} \|$ is uniformly bounded over $\cU^2$. On the event $\Lambda_{jk,4}$, it thus follows from (\ref{ineq_decomp_cov}) and $\|\widehat{\bZ}_{jk}(u,v)\| \leq \|\widehat{\bZ}_{jk}(u,v)-\cE\{\widehat{\bZ}_{jk}(u,v)\}\|+\|\cE\{\widehat{\bZ}_{jk}(u,v)\}\| $ that
\begin{equation}
	\label{decomp_of_cov}
	%\begin{split}
	|\widehat{\Sigma}_{jk}(u,v) - \widetilde{\Sigma}_{jk}(u,v)|  %& \le c \|\widehat{R}_{j}(u)- E\{\widehat{R}_{j}(u)\}\| \\
    \leq c \|\widehat{\bZ}_{jk}(u,v)- \cE\{\widehat{\bZ}_{jk}(u,v)\}\| + c \|\widehat{\bXi}_{jk}(u,v) - \cE\{\widehat{\bXi}_{jk}(u,v)\}\|_{\tF}\\
%\end{split}
\end{equation}
and the positive constant $c$ does not depend on $(u,v)\in \cU^2.$

Combining (\ref{decomp_of_cov}) with (\ref{Xi_bound}), (\ref{coneq.Z}) and applying the first part of Lemma~\ref{l_2 lemma} yields that, for any $ (u,v) \in \cU^2$ and integer $ q \geq 1 $,
$$
\cE\Big\{|\widehat{\Sigma}_{jk}(u,v) - \widetilde{\Sigma}_{jk}(u,v)|^{2q} \Big |\Lambda_{jk,4}\Big\} \leq q! c \Big(\frac{4}{c \nu_{n,T,h,jk}}\Big)^q + (2q)! c \Big(\frac{4}{c \nu_{n,T,h,jk}}\Big)^{2q}.
$$
Applying the second part of Lemma~\ref{l_2 lemma}, we can show that, for each $\delta \in (0,1],$ 
%\begin{flalign*}
$$\pr\Big(\|\widehat \Sigma_{jk}-\widetilde \Sigma_{jk}\|_{\cS} \geq \delta\Big) \leq \pr\Big(\|\widehat \Sigma_{jk}-\widetilde \Sigma_{jk}\|_{\cS} \geq \delta \big|\Lambda_{jk,4}\Big) + \pr\big(\Lambda_{jk,4}^c\big)  \leq c_4 \exp\big(-{c_3\nu_{n,T,h,jk}\delta^2}\big),$$%\\
%\end{flalign*}
which means that (\ref{cov_coneq_L2}) in Theorem~\ref{thm_cov_coneq} holds and completes the proof of concentration inequality for the covariance estimator in Hilbert--Schmidt norm.

\subsubsection{Concentration Inequality in Supremum Norm}
We will derive the uniform concentration bound of $ \sup_{(u,v) \in \cU^2} |\widehat \Sigma_{jk}(u,v)-\widetilde \Sigma_{jk}(u,v)|.$
We partition the interval $ \cU  =[0,1]$ into $ N $ subintervals $I_s$ for $s \in [N]$ of equal length.
Let $ u_s$ and $v_{s'} $ be the centers of $I_s$ and $I_{s'},$ respectively, then we have
\begin{equation}
	\begin{split}
		\sup_{(u,v) \in \cU^2} |\widehat \Sigma_{jk}(u,v) - \widetilde \Sigma_{jk}(u,v)| \leq
		& \max_{s, s' \in [N]} \Big[\big|\widehat \Sigma_{jk}(u_s, v_{s'}) - \widetilde \Sigma_{jk}(u_s, v_{s'})\big|\Big. \\
		& +\Big. \big|\{\widehat \Sigma_{jk}(u, v)-\widehat \Sigma_{jk}(u_s, v_{s'})\} - \{\widetilde \Sigma_{jk}(u, v) - \widetilde \Sigma_{jk}(u_s, v_{s'})\}\big|\Big].
	\end{split} \nonumber
\end{equation}
We need to bound the second term.
By some calculations, it suffices to bound $\Big|\{\widehat{Z}_{jkm}(u,v) - \widehat{Z}_{jkm}(u_s,v_{s'})\} - \big[\cE\{\widehat{Z}_{jkm}(u,v)\} - \cE\{\widehat{Z}_{jkm}(u_s,v_{s'})\}\big]\Big|$ and $\Big|\{\widehat{\Xi}_{jkml}(u,v) - \widehat{\Xi}_{jkml}(u_s,v_{s'})\} - \big[\cE\{\widehat{\Xi}_{jkml}(u,v)\} - \cE\{\widehat{\Xi}_{jkml}(u_s,v_{s'})\}\big]\Big|$ for $ m,l = 1,2,3,$ which means that we need to bound $ \big|\widehat{Z}_{jkm}(u,v) - \widehat{Z}_{jkm}(u_s,v_{s'})\big| $ and $ \big|\widehat{\Xi}_{jkml}(u,v) - \widehat{\Xi}_{jkml}(u_s,v_{s'})\big|.$
Let $ (u,v) \in I_s \times I_{s'} $ and consider $ |\widehat{Z}_{jk1}(u,v) - \widehat{Z}_{jk1}(u_s,v_{s'})| $ for the case of $j\neq k$ first. The results for the case of $j=k$ can be proved in a similar fashion.
Define the event $ \Lambda_{Z,jk1} = \big\{\sum_{i =1}^n w_{ijk} \sum_{t=1}^{T_{ij}} \sum_{t'=1}^{T_{ik}} |\Theta_{ijktt'}| \leq \cE(\sum_{i =1}^n w_{ijk} \sum_{t=1}^{T_{ij}} \sum_{t'=1}^{T_{ik}} |\Theta_{ijktt'}|) +1\big\}.$ On this event, it follows from  Assumption~\ref{cond_kernel}(ii) that
\begin{equation}
	\begin{split}
		& \Big|\widehat{Z}_{jk1}(u,v) - \widehat{Z}_{jk1}(u_s,v_{s'})\Big| \\
		\leq & \Big|\sum_{i=1}^n w_{ijk} \sum_{t=1}^{T_{ij}} \sum_{t'=1}^{T_{ik}} \Theta_{ijktt'}\Big[\{K_h(U_{ijt} - u) - K_h(U_{ijt} - u_s)\} K_h(U_{ikt'} - v) \\
        & + \big\{K_h(U_{ikt'} - v) - K_h(U_{ikt'} - v_{s'})\big\} K_h(U_{ijt} - u_s)\Big]\Big| \\
	\leq	&  \frac{c(|u-u_s| \vee |v-v_{s'}|)}{h^2}  \sum_{i=1}^n w_{ijk} \sum_{t=1}^{T_{ij}} \sum_{t'=1}^{T_{ik}} \Big|\Theta_{ijktt'}\Big| \big\{K_h(U_{ikt'} - v) + K_h(U_{ijt} - u_s)\big\}\\
       \leq  & \frac{c}{Nh^3} \left\{\cE\Big(\sum_{i =1}^n w_{ijk} \sum_{t=1}^{T_{ij}} \sum_{t'=1}^{T_{ik}} \Big| \Theta_{ijktt'}\Big|\Big) +1\right\} \leq \frac{c}{Nh^3}. \\
	\end{split} \nonumber
\end{equation}
Applying similar techniques as above, we can define events $ \Lambda_{Z,jkm} $ and $ \Lambda_{\Xi,jkml} $ for $ m,l = 1,2,3.$
On the intersection of these events, we can obtain that $|\widehat{Z}_{jkm}(u,v) - \widehat{Z}_{jkm}(u_s,v_{s'})| \leq c (Nh^3)^{-1}$ and $|\widehat{\Xi}_{jkml}(u,v) - \widehat{\Xi}_{jkml}(u_s,v_{s'})| \leq c (Nh^3)^{-1}.$
Combing	the above results, we have
\begin{equation}
	\sup_{(u,v) \in \cU^2} |\widehat \Sigma_{jk}(u,v) - \widetilde \Sigma_{jk}(u,v)| \leq  \max_{s, s'\in [N]} \big|\widehat \Sigma_{jk}(u_s, v_{s'}) - \widetilde \Sigma_{jk}(u_s, v_{s'})\big| + \frac{c}{Nh^3}.
	\nonumber
\end{equation}
By the Bernstein inequality, we have $\pr(\Lambda_{Z,jkm}^{\rm c}) \leq \exp\{-\lambda + c_1^2 \lambda^2 (n-c_2 \lambda)^{-1}\}$ for $\lambda \in (0,n c^{-1}].$ With $\lambda=n(2c_1^2+c_2)^{-1},$ the right-side reduces to $ \exp(-cn) \leq \exp(-c \nu_{n,T,h,jk})$ for $m,l=1,2,3.$
Similarly, $ \pr(\Lambda_{\Xi,jkml}^{\rm c}) \leq \exp(-cn) \leq \exp(-c \nu_{n,T,h,jk})$ for $m,l=1,2,3.$
It follows from the above results and the union bound of probability with the choice of $N = \lfloor c (h^3 \delta)^{-1} \rfloor$ that there exist some positive constants $c_3$ and $c_4$ such that, for any $ \delta \in (0,1],$
\begin{equation}
	\label{2}
	\pr\Big\{\sup_{(u,v) \in \cU^2} |\widehat \Sigma_{jk}(u,v) - \widetilde \Sigma_{jk}(u,v)| \geq \delta\Big\} \leq \frac{c_4}{h^6 \delta^2} \exp(-c_3\nu_{n,T,h,jk} \delta^2).
\end{equation}
Take arbitrarily small $\epsilon_2>0.$ If $ n^{\epsilon_2} \nu_{n,T,h,jk} \delta^2 \geq 1,$ the right side of (\ref{2}) reduces to $ c_4 n^{\epsilon_2} \nu_{n,T,h,jk} h^{-6} \exp(-c_3 \nu_{n,T,h,jk} \delta^2).$
If $ n^{\epsilon_2} \nu_{n,T,h,jk} \delta^2 \leq 1,$ we can choose $ c_4 $ and $ n^{\epsilon_2} > c $ such that $ c_4 \exp(-c_3 c^{-1}) \geq 1 $ and hence the same bound $c_4 n^{\epsilon_2} \nu_{n,T,h,jk} h^{-6} \exp(-c_3 \nu_{n,T,h,jk} \delta^2)$ can still be used.
We complete the proof of (\ref{cov_coneq_sup}) in Theorem~\ref{thm_cov_coneq}, the concentration inequality for the covariance estimator in supremum norm. $\hfill\blacksquare$

\subsection{Proof of Theorem~\ref{thm_mean_maxrate}}

Note that $ \|\widehat \mu_{j}-\mu_{j}\|_2 \leq \|\widehat \mu_{j}- \widetilde \mu_{j}\|_2 + \|\widetilde \mu_{j}-\mu_{j}\|_2,$ it suffices to bound $\|\widetilde \mu_j -\mu_j\|_2.$
By (\ref{tilde_mu}), for any $u \in \cU,$ $$\widetilde \mu_{j}(u)-\mu_{j}(u) = \be_0^{\T} \Big[\cE\big\{\widehat{\bS}_j(u)\big\}\Big]^{-1} \cE\Big[\widehat{\bR}_j(u) - \widehat{\bS}_j(u) \big\{\mu_j(u), 0\big\}^{\T}\Big].$$
By the Taylor expansion, we have
$$
\cE_{\varepsilon}\Big[\widehat \bR_j(u) -\widehat \bS_j(u) \{\mu_j(u),0\}^{\T}\Big|V_j\Big] = \sum_{i=1}^{n} v_{ij} \sum_{t=1}^{T_{ij}} \widetilde{\bU}_{ijt} K_{h_{\mu,j}}(U_{ijt} - u) \big\{\mu_j(U_{ijt}) - \mu_j(u)\big\}:= \bJ_1 + \bJ_2,
$$
with
\begin{equation}
\begin{split}
& \bJ_1 = \sum_{i=1}^{n} v_{ij} \sum_{t=1}^{T_{ij}} \widetilde{\bU}_{ijt} K_{h_{\mu,j}}(U_{ijt} - u) \frac{U_{ijt} - u}{h_{\mu,j}} h_{\mu,j} \frac{\partial \mu_j(u)}{ \partial u}, \\
& \bJ_2 = \frac{1}{2} \sum_{i=1}^{n} v_{ij} \sum_{t=1}^{T_{ij}} \widetilde{\bU}_{ijt} K_{h_{\mu,j}}(U_{ijt} - u)  \Big(\frac{U_{ijt} - u}{h_{\mu,j}}\Big)^2 h_{\mu,j}^2 \frac{\partial \mu_j^2(\delta_{ijt})}{\partial u^2},\\
\end{split} \nonumber
\end{equation}
where $\delta_{ijt} \in [u-h_{\mu,j}, u+ h_{\mu,j}]$ and $\cE_{\varepsilon}$ denotes the expectation over $\{Y_{ijt}\}$ in (\ref{model}) conditional on the event $V_j = \{U_{ijt}, t \in [T_{ij}], i \in [n]\}.$
First consider $\bJ_1,$ which equals to the second column of $\widehat{\bS}_j(u)$ multiplied by $ h_{\mu,j} \partial \mu_j(u)/\partial u,$ hence
$$\cE_U\big[\be_0^{\T}[\cE\{\widehat{\bS}_j(u)\}]^{-1} \bJ_1\big] = h_{\mu,j} \frac{\partial \mu_j(u)}{\partial u} \be_0^{\T} [\cE\{\widehat{\bS}_j(u)\}]^{-1} \cE\{\widehat{\bS}_j(u)\} (0,1)^{\T} = 0,$$
where $\cE_U$ denotes the expectation over $V_j.$
Consider $\bJ_2$ next. Under Assumption~\ref{cond_cov_der}, we have $K_1 = \sup_{j \in [p], \xi \in \cU} |\partial \mu_j^2(\xi)/\partial u^2| < \infty.$
Each entry of $|\bJ_2|$ is bounded by the $(1,1)$th entry of $\widehat{\bS}_j(u)$ multiplied by $K_1 h_{\mu,j}^2/2,$ and by $\cE\{ K_h(U_{ijt}-u)\} \leq 1$ we have that the $(1,1)$th entry of $\cE\{\widehat{\bS}_j(u)\}$ is bounded by $1.$ Note that $\cE\{\widehat{\bS}_j(u)\}$ is positive definite. These results together yield that
$$
\big|\cE_U(\be_0^{\T} [\cE\{\widehat{\bS}_j(u)\}]^{-1} \bJ_2)\big| \le \|\cE\{\widehat{\bS}_j(u)\}\|_{\min}^{-1} \|\cE_U(|\bJ_2|)\| \leq c K_1 h_{\mu,j}^2,
$$
%$$\big|E_U(e_0^{\T} [E\{\hat{S}_j(u)\}]^{-1} J_2)\big| \leq e_0^{\T} [E\{\hat{S}_j(u)\}]^{-1} E_U(|J_2|) \leq \frac{1}{2} K_1 h_{\mu,j}^2,$$
which implies that $ |\widetilde \mu_{j}(u)-\mu_{j}(u)| \leq c K_1 h_{\mu,j}^2$ for any $u \in \cU.$ Hence
\begin{equation}
\label{mu.bias}
\sup_{u \in \cU}\big|\widetilde \mu_{j}(u)-\mu_{j}(u)\big| = O(h_{\mu,j}^2)~~\text{and}~~\big\|\widetilde \mu_{j}-\mu_{j}\|_2 = O(h_{\mu,j}^2).
\end{equation}

%Define $ h_{\mu,min} =\min_{1 \leq j \leq p} h_{\mu,j} $, we have $ \|\widehat \mu_{j}-\mu_{j}\| \leq \|\widehat \mu_{j}- \widetilde \mu_{j}\| + \|\widetilde \mu_{j}-\mu_{j}\| $.
%For $ u - h_{\mu,j} \leq t \leq u + h_{\mu,j} $,  by Taylor expansion, $ \mu_j(t) - \mu_j(u) = (t - u) \mu^{(1)}_j(u) + \frac{1}{2} (u - t)^2 \mu^{(2)}_j(\xi) $ where $ u - h_{\mu,j} \leq \xi \leq u + h_{\mu,j} $.
%For any $ u \in \cU $, $ \widetilde \mu_{j}(u)-\mu_{j}(u) = e_0^T[E\{\hat{S}_j(u)\}]E\{\hat{R}_j(u) - \hat{S}_j(u) (\mu_j(u), 0)^T\} $.
%The expectation here can be rewritten as $ E\{\hat{R}_j(u) - \hat{S}_j(u) (\mu_j(u), 0)^T\} = E_{U_{ijk}}[E_{e_{ijk}}\{\hat{R}_j(u) - \hat{S}_j(u) (\mu_j(u), 0)^T | \{U_{ijk}\}\}] $.
%With fixed $ \{U_{ijk}\} $, $ E_{e_{ijk}}\{\hat{R}_j(u) - \hat{S}_j(u) (\mu_j(u), 0)^T | \{U_{ijk}\}\} = \sum_{i=1}^{n} v_{ij} \sum_{k=1}^{T_{ij}} \tilde{U}_{ijk} K_{h_{\mu,j}}(U_{ijk} - u) \{\mu_j(U_{ijk}) - \mu_j(u)\} = \sum_{i=1}^{n} v_{ij} \sum_{k=1}^{T_{ij}} \tilde{U}_{ijk} K_{h_{\mu,j}}(U_{ijk} - u) \{(U_{ijk} - u) \mu_j^{(1)}(u) + \frac{1}{2} (U_{ijk} - u)^2 \mu_j^{(2)}(\xi_{ijk})\} $.
%Consider that $ K_2 = sup_{u - h_{\mu,min} \leq \xi \leq u + h_{\mu,min}} \mu^{(2)}(\xi) < \infty $,
%so we have $ E_{e_{ijk}}\{\hat{R}_j(u) - \hat{S}_j(u) (\mu_j(u), 0)^T | \{U_{ijk}\}\} = O(h_{\mu,j}^2) $,
%then $ E\{\hat{R}_j(u) - \hat{S}_j(u) (\mu_j(u), 0)^T\} = O(h_{\mu,j}^2) $, it shows that $ |\widetilde \mu_{j}(u)-\mu_{j}(u)| = O(h_{\mu,j}^2) $ and $ \|\widetilde \mu_{j}-\mu_{j}\| = O(h_{\mu,j}^2) $.

For $M > 0$ with the choice of $\delta = (\log p)^{1/2} (\min_{j}\gamma_{n,T,h,j})^{-1/2} M \leq 1,$ %which requires $n (\widebar T_{\mu,j} h_{\mu,j} \wedge 1) \geq M^2 \log p,$
it follows from the union bound of probability and (\ref{mean_coneq_L2}) in Theorem~\ref{thm_mean_coneq} that
\begin{eqnarray}
	\nonumber
		& &\pr\Big\{\frac{\max_{j \in [p]}\|\widehat \mu_{j}- \widetilde \mu_{j}\|_2}{(\log p)^{1/2} (\min_{j}\gamma_{n,T,h,j})^{-1/2}} \geq M\Big\}\\ \nonumber
		&\leq& \sum_{j = 1}^{p} \pr\Big\{\|\widehat \mu_{j}- \widetilde  \mu_{j}\|_2\geq M (\log p)^{1/2} (\min_{j}\gamma_{n,T,h,j})^{-1/2} \Big\} \\
		&\leq& \sum_{j = 1}^{p} c \exp\Big(-c \gamma_{n,T,h,j} M^2 \frac{\log p}{\min_{j}\gamma_{n,T,h,j}} \Big)
		\leq  %\sum_{j = 1}^{p} c \exp(-c M^2 \log p)
		c \exp\{(1 - c M^2)\log p\}.
\label{maxrate_mu}
\end{eqnarray}
%When $ p $ goes to infinity, for $n^{-1} \log p$ goes to $0$ and so $n$ tends to infinity,
We can choose a large $M$ such that $1-cM^2<0,$ the right-side of (\ref{maxrate_mu}) tends to 0. Hence
$\max_{j \in [p]} \|\widehat \mu_j - \widetilde \mu_j\|_2=O_\p\{(\log p)^{1/2} (\min_j \gamma_{n, T, h,j})^{-1/2}\}.$ Combing this with (\ref{mu.bias})  yields that
\begin{equation}
	\max_{j\in [p]}\|\hat \mu_{j}-\mu_{j}\|_2 = O_\p\left\{\Big(\frac{\log p}{\min_{j}\gamma_{n,T,h,j}}\Big)^{1/2} + \max_j h_{\mu,j}^2\right\}
	\nonumber,
\end{equation}
which completes the proof of (\ref{mean_max_L2_dense}).

The rate of convergence in (\ref{mean_max_sup_dense}) can be proved following a similar procedure.
Let $h_{\mu,\min} = \min_j h_{\mu,j},$ we assume that $h_{\mu,\min} \asymp \{\log (p\vee n)/n\}^{\kappa_1}$ where $\kappa_1 \in (0, 1/2].$
Consider $h_{\mu,\min}$ with some $\kappa_1^* > 1/2$, the corresponding rate is not faster than that with $\kappa_1 \in (0, 1/2].$ To be specific,
under the sparse design, the rate of $\max_{j \in [p]} \sup_{u \in \cU} |\hat \mu_{j}(u)-\mu_{j}(u)|$ is $\{\log (p\vee n)/n\}^{1/2-\kappa_1^*/2},$ which is slower than $\{\log (p\vee n)/n\}^{2/5}$ with $\kappa_1=1/5.$
Under the dense design with $\widebar T_{\mu} \{\log (p \vee n)/n\}^{\kappa_1^*} \rightarrow 0$ and $\widebar T_{\mu} \{\log (p \vee n)/n\}^{3/2} \rightarrow 0,$ the rate is $\{\log (p\vee n)/n\}^{1/2-\kappa_1^*/2} \widebar T_{\mu}^{-1/2},$ which is slower than $\{\log (p\vee n)/n\}^{2/5} \widebar T_{\mu}^{-2/5}$ with $\kappa_1 \in (1/5, 1/2].$
Under the dense design with $\widebar T_{\mu} \{\log (p \vee n)/n\}^{\kappa_1^*} \rightarrow 0$ and $\widebar T_{\mu} \{\log (p \vee n)/n\}^{3/2} \rightarrow \tilde c$ or $\infty,$ the rate is $\{\log (p\vee n)/n\}^{1/2-\kappa_1^*/2} \widebar T_{\mu}^{-1/2},$ which is slower than $\{\log (p\vee n)/n\}^{1/2}$ with $\kappa_1 = 1/4.$
Under the dense design with $\widebar T_{\mu} \{\log (p \vee n)/n\}^{\kappa_1^*} \rightarrow \tilde c$ or $\infty,$ the rate is $\{\log (p\vee n)/n\}^{1/2},$ which is the same as $\{\log (p\vee n)/n\}^{1/2}$ with $\kappa_1 = 1/4.$
Based on the above four cases, if $\kappa_1^* > 1/2,$ the corresponding rate is not faster than that with some $\kappa_1 \in (0, 1/2]$ and hence the $\kappa_1$ that corresponds with the optimal bandwidth under sparse or dense design is in $(0,1/2].$
For $M > 0$ with the choice of $\delta = \{\log (p \vee n)\}^{1/2} (\min_{j}\gamma_{n,T,h,j})^{-1/2} M \leq 1,$  by the union bound of probability and (\ref{mean_coneq_sup}) in Theorem~\ref{thm_mean_coneq},  we have
\begin{eqnarray}
& &\pr\Big\{\frac{\max_{j \in [p]} \sup_{u \in \cU} |\widehat \mu_{j}(u)- \widetilde \mu_{j}(u)|}{\{\log (p \vee n)\}^{1/2} ( \min_{j} \gamma_{n,T,h,j} )^{-1/2}} \geq M\Big\} \nonumber\\
&\leq & \sum_{j = 1}^{p} \pr\Big\{\sup_{u \in \cU} |\widehat \mu_{j}(u)- \widetilde  \mu_{j}(u)| \geq M \{\log (p \vee n)\}^{1/2} ( \min_{j} \gamma_{n,T,h,j} )^{-1/2}\Big\} \nonumber\\
&\leq & \sum_{j = 1}^{p} \frac{c (n^{\epsilon_1}  \gamma_{n,T,h,j})^{1/2}}{h_{\mu,\min}^2} \exp\Big\{-c \gamma_{n,T,h,j} M^2 \frac{\log(p \vee n)}{\min_{j}\gamma_{n,T,h,j}}\Big\} \nonumber\\
&\leq& c \exp\Big\{ \Big(\frac{3+\epsilon_1}{2} + 2 \kappa_1 - c M^2\Big)\log(p \vee n)\Big\}. \label{maxrate_mu1}
\end{eqnarray}
We can choose a large $M$ such that $3+\epsilon_1+4\kappa_1-2cM^2<0,$ the right-side of (\ref{maxrate_mu1}) tends to 0.
Hence $$\max_{j \in [p]} \sup_{u \in \cU} \big|\hat \mu_{j}(u)-\widetilde \mu_{j}(u)\big| = O_\p\left[\Big\{\frac{\log (p \vee n)}{\min_{j}\gamma_{n,T,h,j}}\Big\}^{1/2}\right].$$
Combing this with (\ref{mu.bias}) yields that
$$
\max_{j \in [p]} \sup_{u \in \cU} \big|\hat \mu_{j}(u)-\mu_{j}(u)\big| = O_\p\left[\Big\{\frac{\log (p \vee n)}{\min_{j}\gamma_{n,T,h,j}}\Big\}^{1/2} + \max_j h_{\mu,j}^2\right],
$$
which completes the proof of (\ref{mean_max_sup_dense}). $\hfill\blacksquare$

\subsection{Proof of Theorem~\ref{thm_cov_maxrate}}
Note $ \|\widehat \Sigma_{jk}- \Sigma_{jk}\|_{\cS} \leq \|\widehat \Sigma_{jk}- \widetilde \Sigma_{jk}\|_{\cS} + \|\widetilde \Sigma_{jk}- \Sigma_{jk}\|_{\cS},$ it suffices to bound $\|\widetilde \Sigma_{jk}- \Sigma_{jk}\|_{\cS}.$
By (\ref{Sigma_tilde}), for any $(u,v) \in \cU^2,$ $$\widetilde{\Sigma}_{jk} (u, v) - \Sigma_{jk} (u, v) = \tilde \be_0^{\T}	\Big[\cE \big\{ \widehat{\bXi}_{jk}(u, v)\big\}\Big]^{-1} \cE\Big[\widehat{\bZ}_{jk}(u, v) - \widehat{\bXi}_{jk}(u, v) \big\{\Sigma_{jk}(u, v),0,0\big\}^{\T}\Big].$$
By the Taylor expansion, we have
\begin{equation}
\begin{split}
    & \cE_{\varepsilon}\Big[\widehat{\bZ}_{jk}(u, v) - \widehat{\bXi}_{jk}(u, v) \big\{\Sigma_{jk}(u, v), 0, 0 \big\}^{\T} | \widetilde V_{jk}\Big] \\
    & = \sum_{i=1}^{n} w_{ijk} \sum_{t=1}^{T_{ij}} \sum_{s=1}^{T_{ik}} \widetilde{\bU}_{ijkts} K_{h_{\sSigma,jk}}(U_{ijt} - u) K_{h_{\sSigma,jk}}(U_{iks} - v) \big\{\Sigma_{jk} (U_{ijt}, U_{iks}) - \Sigma_{jk} (u, v)\big\} \\
    & := \bL_1 + \bL_2 + \bL_3,
\end{split} \nonumber
\end{equation}
where
\begin{equation}
\begin{split}
    \bL_1 &= \sum_{i=1}^{n} w_{ijk} \sum_{t=1}^{T_{ij}} \sum_{s=1}^{T_{ik}} \widetilde{\bU}_{ijkts} K_{h_{\sSigma,jk}}(U_{ijt} - u) K_{h_{\sSigma,jk}}(U_{iks} - v) \frac{U_{ijt} - u}{h_{\sSigma,jk}} h_{\sSigma,jk} \frac{\partial \Sigma_{jk}}{\partial u} (u,v), \\
    \bL_2 &= \sum_{i=1}^{n} w_{ijk} \sum_{t=1}^{T_{ij}} \sum_{s=1}^{T_{ik}} \widetilde{\bU}_{ijkts} K_{h_{\sSigma,jk}}(U_{ijt} - u) K_{h_{\sSigma,jk}}(U_{iks} - v) \frac{U_{iks} - v}{h_{\sSigma,jk}} h_{\sSigma,jk} \frac{\partial \Sigma_{jk}}{\partial v} (u,v), \\
    \bL_3 &= \sum_{i=1}^{n} w_{ijk} \sum_{t=1}^{T_{ij}} \sum_{s=1}^{T_{ik}} \widetilde{\bU}_{ijkts} K_{h_{\sSigma,jk}}(U_{ijt} - u) K_{h_{\sSigma,jk}}(U_{iks} - v) \widetilde L_{ijk}, \\
    \widetilde L_{ijk} &= \frac{1}{2} h_{\sSigma,jk}^2  \Big(\frac{U_{ijt} - u}{h_{\sSigma,jk}}\Big)^2 \frac{\partial^2 \Sigma_{jk}}{\partial u^2} (\delta_{ijkts1},\delta_{ijkts2}) + \frac{1}{2} h_{\sSigma,jk}^2  \Big(\frac{U_{iks} - v}{h_{\sSigma,jk}}\Big)^2 \frac{\partial^2 \Sigma_{jk}}{\partial v^2} (\delta_{ijkts1},\delta_{ijkts2}) \\
    & +  h_{\sSigma,jk}^2  \frac{U_{ijt} - u}{h_{\sSigma,jk}} \frac{U_{iks} - v}{h_{\sSigma,jk}} \frac{\partial^2 \Sigma_{jk}}{\partial u \partial v} (\delta_{ijkts1},\delta_{ijkts2}),
\end{split} \nonumber
\end{equation}
$(\delta_{ijkts1},\delta_{ijkts2}) \in [u-h_{\sSigma,jk},u+h_{\sSigma,jk}] \times [v-h_{\sSigma,jk},v+h_{\sSigma,jk}]$ and the event $\widetilde V_{jk}=\{(U_{ijt}, U_{iks}), t \in [T_{ij}], s \in [T_{ik}], i \in [n]\}.$
First consider $L_1,$ which equals to the second column of $ \widehat{\bXi}_{jk}(u, v)$ multiplied by $ h_{\sSigma,jk} \partial \Sigma_{jk}(u,v) /\partial u,$ hence $\cE_U\Big(\tilde \be_0^{\T} \Big[\cE \big\{ \widehat{\bXi}_{jk}(u, v)\big\}\Big]^{-1} \bL_1\Big)$ equals to
$$h_{\sSigma,jk} \frac{\partial \Sigma_{jk}}{\partial u} (u,v) \tilde \be_0^{\T} \Big[\cE \big\{ \widehat{\bXi}_{jk}(u, v)\big\}\Big]^{-1} \cE \Big\{ \widehat{\bXi}_{jk}(u, v)\Big\} (0,1,0)^{\T} =0,$$ where $\cE_U$ denotes the expectation over $\widetilde V_{jk}.$
Following the similar procedure, we can show that
$$\cE_U\Big(\tilde \be_0^{\T} \Big[\cE \big\{ \widehat{\bXi}_{jk}(u, v)\big\}\Big]^{-1} \bL_2\Big)=0.$$
%Similarly, $N_2$ equals to the third column of $ \widehat{\Xi}_{jk}(u, v)$ multiplied by $h_{\sSigma,jk} \partial \Sigma_{jk}(u,v) /\partial v,$ hence
%$$E_U\Big(e_0^{\T} \Big[E \big\{ \widehat{\Xi}_{jk}(u, v)\big\}\Big]^{-1} N_2\Big) = h_{\sSigma,jk} \frac{\partial \Sigma_{jk}}{\partial v} (u,v) e_0^{\T} \Big[E \big\{ \widehat{\Xi}_{jk}(u, v)\big\}\Big]^{-1} E \Big\{ \widehat{\Xi}_{jk}(u, v)\Big\} (0,1,0)^{\T} =0.$$
Then consider $\bL_3.$ By Assumption~\ref{cond_cov_der}, we have $$ K_2= \sup_{(j,k) \in [p]^2,(u,v) \in \cU^2} \{\partial^2 \Sigma_{jk}(u,v)/\partial v^2,\partial^2 \Sigma_{jk}(u,v) /\partial u^2 , \partial^2 \Sigma_{jk}(u,v)/\partial u \partial v\} < \infty.$$
Each entry of $|\bL_3|$ is bounded by the $(1,1)$th entry of $\widehat{\bXi}_{jk}(u, v)$ multiplied by $2 K_2 h_{\sSigma,jk}^2,$
and by $\cE\{ K_h(U_{ijt}-u)K_h(U_{iks}-v)\}\leq 1,$ the $(1,1)$th entry of $\cE\{\widehat{\bXi}_{jk}(u, v)\}$ is bounded by $1.$
Note that $\cE\{ \widehat{\bXi}_{jk}(u, v)\}$ is positive definite. Combining these results, we have
$$
\Big|\cE_U\Big(\tilde \be_0^{\T} \Big[\cE\big\{ \widehat{\bXi}_{jk}(u, v)\big\}\Big]^{-1} \bL_3\Big)\Big| \leq \|\cE\{\widehat{\bXi}_{jk}(u, v)\}\|^{-1}_{\min} \|\cE_U(|\bL_3|)\| \leq c K_2 h_{\sSigma,jk}^2,
$$
%$$
%\Big|E_U(e_0^{\T} \Big[E\big\{ \widehat{\Xi}_{jk}(u, v)\big\}\Big]^{-1} L_3)\Big| \leq e_0^{\T} \Big[E\big\{ \widehat{\Xi}_{jk}(u, v)\big\}\Big]^{-1} E_U(|L_3|) \leq \frac{1}{2} K_2 h_{\sSigma,jk}^2,
%$$
which implies that $|\widetilde \Sigma_{jk}(u, v)- \Sigma_{jk}(u, v)|\leq c K_2 h_{\sSigma,jk}^2$ for any $(u,v) \in \cU^2.$
Hence
\begin{equation}
	\label{cov.bias}
	\sup_{(u,v) \in \cU^2}\big|\widetilde \Sigma_{jk}(u, v)- \Sigma_{jk}(u, v)\big| = O(h_{\sSigma,jk}^2)~~\text{and}~~\big\|\widetilde \Sigma_{jk}-\Sigma_{jk}\|_\cS = O(h_{\sSigma,jk}^2).
\end{equation}
For $ M >0 $ with the choice of $\delta = (\log p)^{1/2} (\min_{j,k} \nu_{n,T,h,jk})^{-1/2} M \leq 1$, it follows from the union bound of probability and (\ref{cov_coneq_L2}) in Theorem~\ref{thm_cov_coneq} that
\begin{flalign}	\label{maxrate_cov}
\begin{split}
& \pr\Big\{\frac{\max_{j,k \in [p]}\|\widehat \Sigma_{jk}- \widetilde \Sigma_{jk}\|_{\cS}}{(\log p)^{1/2} (\min_{j,k} \nu_{n,T,h,jk})^{-1/2}} \geq M\Big\} \\
\leq & \sum_{j=1}^p \sum_{k=1}^p \pr\Big\{\|\widehat \Sigma_{jk}- \widetilde  \Sigma_{jk}\|_{\cS} \geq M (\log p)^{1/2} (\min_{j,k} \nu_{n,T,h,jk})^{-1/2}\Big\} \\
\leq & \sum_{j=1}^p \sum_{k=1}^p c \exp\Big(-c \nu_{n,T,h,jk} M^2 \frac{\log p}{\min_{j,k} \nu_{n,T,h,jk}}\Big) \leq  c \exp\{ (2 - c M^2)\log p\}. 
\end{split}
\end{flalign}
We can choose a large $M$ such that $2-cM^2<0,$ the right-side of (\ref{maxrate_cov}) tends to 0.
Hence $\max_{j,k \in [p]} \|\widehat \Sigma_{jk}- \widetilde \Sigma_{jk}\|_{\cS}=O_\p\{(\log p)^{1/2} (\min_{j,k} \nu_{n,T,h,jk})^{-1/2}\}.$ Combing this with (\ref{cov.bias}) yields that
\begin{equation}
	\max_{j,k \in [p]}\|\widehat \Sigma_{jk}-\Sigma_{jk}\|_{\cS} = O_\p\left\{\Big(\frac{\log p}{\min_{j,k} \nu_{n,T,h,jk}}\Big)^{1/2} + \max_{j,k} h_{\Sigma,jk}^2\right\}, \nonumber
\end{equation}
which completes the proof of (\ref{cov_max_L2_dense}).

The rate of convergence in (\ref{cov_max_sup_dense}) can be proved following a similar procedure.
Let $h_{\sSigma,\min} = \min_{j,k} h_{\sSigma,jk},$ we can assume that $h_{\sSigma,\min} \asymp \{\log (p\vee n)/n\}^{\kappa_2}$ where $\kappa_2 \in (0, 1/2].$
Consider $h_{\sSigma,\min}$ with some $\kappa_2^* > 1/2,$ the corresponding rate is not faster than that with some $\kappa_2 \in (0, 1/2].$
Specifically, under the sparse design, the rate of $\max_{j,k \in [p]} \sup_{(u,v) \in \cU^2} |\widehat \Sigma_{jk}(u,v)-\Sigma_{jk}(u,v)|$ is $\{\log (p\vee n)/n\}^{1/2-\kappa_2^*},$ which is slower than $\{\log (p\vee n)/n\}^{1/3}$ with $\kappa_2=1/6.$
Under the dense design with $\widebar T_{\sSigma} \{\log (p \vee n)/n\}^{\kappa_2^*} \rightarrow 0$ and $\widebar T_{\sSigma} \{\log (p \vee n)/n\} \rightarrow 0,$ the rate %of $\max_{j,k \in [p]} \sup_{(u,v) \in \cU^2} |\widehat \Sigma_{jk}(u,v)-\Sigma_{jk}(u,v)|$
is $\{\log (p\vee n)/n\}^{1/2-\kappa_2^*} \widebar T_{\sSigma}^{-1},$ which is slower than $\{\log (p\vee n)/n\}^{1/3} \widebar T_{\sSigma}^{-2/3}$ with $\kappa_2 \in (1/6, 1/2].$
Under the dense design with $\widebar T_{\sSigma} \{\log (p \vee n)/n\}^{\kappa_2^*} \rightarrow 0$ and $\widebar T_{\sSigma} \{\log (p \vee n)/n\} \rightarrow \widetilde c$ or $\infty,$ the rate %of $\max_{j,k \in [p]} \sup_{(u,v) \in \cU^2} |\widehat \Sigma_{jk}(u,v)-\Sigma_{jk}(u,v)|$
is $\{\log (p\vee n)/n\}^{1/2-\kappa_2^*} \widebar T_{\sSigma}^{-1},$ which is slower than $\{\log (p\vee n)/n\}^{1/2}$ with $\kappa_2 = 1/4.$
Under the dense design with $\widebar T_{\sSigma} \{\log (p \vee n)/n\}^{\kappa_2^*} \rightarrow \widetilde c$ or $\infty,$ the rate %$\max_{j,k \in [p]} \sup_{(u,v) \in \cU^2} |\widehat \Sigma_{jk}(u,v)-\Sigma_{jk}(u,v)|$
is $\{\log (p\vee n)/n\}^{1/2},$ which is the same as $\{\log (p\vee n)/n\}^{1/2}$ with $\kappa_2 = 1/4.$
Based on the above four cases, if $\kappa_2^* > 1/2,$ the corresponding rate is not faster than that with some $\kappa_2 \in (0, 1/2]$ and hence $\kappa_2$ that corresponds with the optimal bandwidth under sparse or dense design is in $(0,1/2].$
For $M > 0$ with the choice of $\delta=\{\log (p \vee n)\}^{1/2} (\min_{j}\gamma_{n,T,h,j})^{-1/2} M \leq 1,$ by the union bound of probability and (\ref{cov_coneq_sup}) in Theorem~\ref{thm_cov_coneq}, we have
\begin{eqnarray*}
	\label{maxsuprate_cov}
	%\begin{split}
		& & \pr\Big\{\frac{\max_{j,k \in [p]} \sup_{(u,v) \in \cU^2} |\widehat \Sigma_{jk}(u,v)- \widetilde \Sigma_{jk}(u,v)|}{{\log(p \vee n)}^{1/2} (\min_{j,k} \nu_{n,T,h,jk})^{-1/2}} \geq M\Big\} \nonumber\\
		&\leq & \sum_{j=1}^p \sum_{k=1}^p\pr\Big\{\sup_{(u,v) \in \cU^2} |\widehat \Sigma_{jk}(u,v)- \widetilde \Sigma_{jk}(u,v)| \geq M \{\log(p \vee n)\}^{1/2} (\min_{j,k} \nu_{n,T,h,jk})^{-1/2}\Big\} \nonumber\\
		&\leq & \sum_{j=1}^p \sum_{k=1}^p \frac{c n^{\epsilon_2} \nu_{n,T,h,jk}}{h_{\sSigma,\min}^6} \exp\Big\{-c \nu_{n,T,h,jk} M^2 \frac{\log(p \vee n)}{\min_{j,k} \nu_{n,T,h,jk}} \Big\} \nonumber\\
		&\leq & c \exp\{ (3 +\epsilon_2 + 6\kappa_2 - c M^2)\log(p \vee n)\}. %\\
	%\end{split}
\end{eqnarray*}
We can choose a large $M$ such that $3+\epsilon_2+6\kappa_2-cM^2<0,$ the right-side of the above inequality %(\ref{maxsuprate_cov})
tends to 0.
Hence $$\max_{j,k \in [p]} \sup_{(u,v) \in \cU^2} |\widehat \Sigma_{jk}(u,v)- \widetilde \Sigma_{jk}(u,v)| = O_\p\left[\Big\{\frac{\log (p \vee n)}{\min_{j,k} \nu_{n,T,h,jk}}\Big\}^{1/2}\right].$$
Combing this with (\ref{cov.bias}) yields that
\begin{equation}
	\max_{j,k \in [p]} \sup_{(u,v) \in \cU^2} \big|\widehat \Sigma_{jk}(u,v)-\Sigma_{jk}(u,v)\big| = O_\p\left[\Big\{\frac{\log (p \vee n)}{\min_{j,k} \nu_{n,T,h,jk}}\Big\}^{1/2} + \max_{j,k} h_{\Sigma,jk}^2\right], \nonumber
\end{equation}
which completes the proof of (\ref{cov_max_sup_dense}). $\hfill\blacksquare$

\subsection{Proof of Proposition~\ref{thm_fpca_maxrate}}

%Let $\delta_{jl} = {\min}_{1 \le k \le l}\{\lambda_{jk} - \lambda_{j(k+1)}\}$ and $\widehat \Delta_{jj} = \widehat \Sigma_{jj} - \Sigma_{jj}$ for $j=1, \dots,p$ and $l=1,2 \dots.$ 
It follows from (4.43) and Lemma~4.3 of \cite{Bbosq1} that
\begin{equation*}
\label{eigen.bd}
\underset{l \in [d_j]}{\sup}~ |\widehat \lambda_{jl} - \lambda_{jl}| \le \|\widehat \Sigma_{jj} - \Sigma_{jj}\|_{\cS}~~\mbox{and}~~\underset{l \in [d_j]}{\sup} ~\delta_{jl}\|\widehat \phi_{jl} - \phi_{jl}\|_2 \le 2\sqrt{2} \|\widehat \Sigma_{jj} - \Sigma_{jj}\|_{\cS}.
\end{equation*}
Combining the above with Theorem~\ref{thm_cov_maxrate} yields the result in Proposition~\ref{thm_fpca_maxrate}. $\hfill\blacksquare$

\section{Verification of the Claim in Section~\ref{sec.sim}}
\label{supp_claim}

For the mean estimator, define the set of $r$ candidate bandwidths ${\cal H}_{\mu} = \{h_{\mu}^{(1)},\dots,h_{\mu}^{(r)}\}.$
In our simulations, the bandwidth for each dimension can be chosen from ${\cal H}_{\mu}$ freely, and hence there are $ r^p $ possible outcomes. The targeted evaluation metric is $ \text{global}_{\text{opt}} = \min_{(m_1,...,m_p) \in [r]^p} \max_{j\in [p]} \text{MISE}(\widehat{\mu}_j, h_{\mu}^{(m_j)}).$
We will show that $ \text{global}_{\text{opt}} = \text{MaxMISE}(\mu),$ the right side of which is much easier to calculate as it only takes into account $pr$ cases.
On one hand, it is obvious that $\text{global}_{\text{opt}} \leq \text{MaxMISE}(\mu)$.
On the other hand, for fixed $ (j, m_j) \in [p] \times [r],$ $ \text{MISE}(\widehat{\mu}_j, h_{\mu}^{(m_j)}) \geq \min_{m \in [r]} \text{MISE}(\widehat{\mu}_j, h_{\mu}^{(m)}),$
%For any $ (i_1,...,i_p) \in [r]^p $, $ \max_{j \in [p]} MISE(\widehat{\mu}_j, h_{i_j}) \geq \max_{j \in [p]} \min_{i \in [r]} MISE(\widehat{\mu}_j, h_{i}) $.
and hence $\text{global}_{\text{opt}} \geq$
$$\min_{(m_1,...,m_p) \in [r]^p} \max_{j \in [p]} \min_{m \in [r]}  \text{MISE}(\widehat{\mu}_j, h_{\mu}^{(m)})  = \max_{j \in [p]} \min_{m \in [r]} \text{MISE}(\widehat{\mu}_j, h_{\mu}^{(m)}) = \text{MaxMISE}(\mu).
$$
%\begin{equation}
%	\begin{split}
%		global_{opt} & \leq \min_{(i_1,...,i_p) \in [r]^p} \max_{j \in [p]} \min_{i \in [r]} MISE(\widehat{\mu}_j, h_{i}) \\
%		& = \max_{j \in [p]} \min_{i \in [r]} MISE(\widehat{\mu}_j, h_{i}) \\
%		& = MaxMISE(\mu), \\
%	\end{split} \nonumber
%\end{equation}
Combining the above results yields $ \text{global}_{\text{opt}} = \text{MaxMISE}(\mu).$
The corresponding claim for the covariance estimator can be verified in the same way. %$\hfill\blacksquare$

\vskip 0.2in
\bibliography{paperbib}
\end{document}